\documentclass{amsart}
\pdfoutput=1
\usepackage{float}
\floatstyle{boxed}
\restylefloat{figure}
\usepackage{mathtools}
\usepackage{amssymb}
\usepackage{amsthm}
\usepackage{mathrsfs}
\usepackage[all]{xy}
\UseComputerModernTips
\CompileMatrices
\usepackage{graphicx} 
\usepackage[colorlinks,allcolors=blue]{hyperref}
\usepackage{ifthen}
\usepackage{latexsym}
\allowdisplaybreaks
\numberwithin{equation}{section}
\newtheorem{theorem}[equation]{Theorem}
\newtheorem{lemma}[equation]{Lemma}
\newtheorem{proposition}[equation]{Proposition}
\newtheorem{corollary}[equation]{Corollary}
\theoremstyle{definition}
\newtheorem{definition}[equation]{Definition}

\newtheorem{notation}[equation]{Notation}
\theoremstyle{remark}
\newtheorem{remark}[equation]{Remark}

\renewcommand{\phi}{\varphi}

\DeclareMathSymbol{\boxprod}{\mathbin}{AMSa}{"03} 
\DeclareMathSymbol{\mixprod}{\mathbin}{AMSa}{"4F} 

\newcommand{\convto}{\Rightarrow}
\newcommand{\dirsum}{\oplus}
\newcommand{\Dirsum}{\bigoplus}
\newcommand{\disjunion}{\sqcup}
\newcommand{\Disjunion}{\coprod}
\newcommand{\dual}{^\vee}

\newcommand{\hmtpc}{\simeq}

\newcommand{\includesin}{\hookrightarrow}

\newcommand{\iso}{\cong}

\newcommand{\Mackey}[1]{{\underline {#1}}}

\newcommand{\smsh}{\wedge}

\newcommand{\susp}{\Sigma}
\newcommand{\tensor}{\otimes}
\newcommand{\union}{\cup}

\newcommand{\C}{{\mathbb C}}

\newcommand{\R}{{\mathbb R}}

\newcommand{\V}{{\mathscr V}}
\newcommand{\Vrep}{\mathbb{V}}

\newcommand{\Wrep}{\mathbb{W}}

\newcommand{\Z}{\mathbb{Z}}
\newcommand{\GG}{{C_2}}
\newcommand{\CP}{\C P}
\newcommand{\tE}{\tilde E}

\newcommand{\orb}[1]{{\mathscr{O}_{#1}}}
\newcommand{\sorb}[1]{{\widehat{\mathscr{O}}_{#1}}}
\newcommand{\conc}[1]{\langle #1 \rangle}

\DeclareMathOperator{\Hom}{Hom}

\newcommand{\tensorS}{\tensor_{\Mackey H_\GG^{RO(\GG)}(S^0)}}
\begin{document}
\title{The $RO(\Pi B)$-graded $\GG$-equivariant ordinary cohomology of $B_\GG U(1)$}

\author{Steven R. Costenoble}
\address{Department of Mathematics\\Hofstra University\\
   Hempstead, NY 11549}
\email{Steven.R.Costenoble@Hofstra.edu}



\date{August 19, 2023}

\abstract
We calculate the ordinary $\GG$-cohomology, with Burnside ring coefficients, of 
$\CP_\GG^\infty = B_\GG U(1)$,
the complex projective space, a model for the classifying space for $\GG$-equivariant
complex line bundles. The $RO(\GG)$-graded Bredon ordinary cohomology was
calculated by Gaunce Lewis, but here we extend to a larger grading
in order to capture a more natural set of generators.
These generators include the Euler class
of the tautological bundle, which lies outside of the $RO(\GG)$-graded theory.
\endabstract
\maketitle
\tableofcontents


\section*{Introduction}

There has been a recent resurgence of interest in calculations in Bredon's equivariant ordinary cohomology.
See, for example, Dugger \cite{Dug:AHSSforKR} and \cite{Dug:Grassmannians}, Kronholm \cite{Kro:SerreSS},
Hogle \cite{Hog:Grassmannians}, Hazel \cite{Haz:FundClasses},
and particularly the use of such calculations by Hill, Hopkins, and Ravenel
to solve the Kervaire invariant problem \cite{HHR:Kervaire}.
A large influence on the present paper was the earlier
calculation of the $\GG$-equivariant cohomology of
projective spaces by Lewis in \cite{Le:projectivespaces}.
All of these calculations used $RO(G)$ grading.

A general problem with $RO(G)$-graded ordinary cohomology is that it
does not have a Thom isomorphism, hence does not have Euler classes, for arbitrary vector bundles.
It does have such for vector bundles all of whose fibers are modeled on a single representation $V$
of $G$. In that case, the Thom isomorphism shifts grading by $V$ and the Euler class
of the bundle lives in grading $V$. But many naturally occurring bundles, like the tautological
bundle on a projective space, have varying representations over the space,
and there is no Thom isomorphism or Euler class for such a bundle using $RO(G)$ grading.
Recalling that the nonequivariant cohomology of a projective space is generated by the Euler
class of its tautological bundle, 
it is clear that not having such a class equivariantly is a serious lacuna.

It was with this in mind that Stefan Waner and I wrote \cite{CW:ordinaryhomology},
defining and exploring an ordinary cohomology theory with a grading expanded
beyond $RO(G)$, in which there are natural gradings to use as the dimensions
of arbitrary vector bundles. This theory does possess a Thom isomorphism theorem
for any vector bundle, which allows us to define Euler classes of
arbitrary vector bundles.

The main purpose of this paper is to calculate,
using the expanded grading from \cite{CW:ordinaryhomology},
the ordinary cohomology of
$B_\GG U(1)$, the classifying space for 
complex $\GG$-line bundles, which we can model as an infinite complex projective space.
The result is Theorem~\ref{thm:multstructure} in which,
as expected, the Euler class of the tautological bundle plays a central role
as one of the multiplicative generators.

Applications of these results have already been made in \cite{CHT:projspaces},
where Thomas Hudson, Sean Tilson, and I used this calculation as input
to the calculation of the cohomologies of finite projective spaces.
One interesting point is that the cohomologies of the finite projective spaces
are not just quotient rings of the infinite case, but involve some new elements.
The phenomenon that leads to these new elements is discussed in \S\ref{sec:components}.
We used those calculations to give an equivariant
version of B\'ezout's theorem, an important part of which was the
calculation of the Euler classes of line bundles associated with hypersurfaces.
Once again, those Euler classes live in gradings outside the $RO(\GG)$ grading,
so use of the expanded grading was crucial.
We are working on generalizing those results further, including the possibility
of developing an equivariant Schubert calculus.

This paper is also meant to be a first step toward a useful theory of characteristic classes for
equivariant vector bundles. The next steps will require similar calculations
for the classifying spaces $B_G U(n)$ for larger $n$.

This paper is organized as follows.
Part~\ref{part:cohomology} reviews
some necessary background material and describes ordinary cohomology with the expanded grading.
The equivariant cohomology of a $\GG$-space is a module over the
$RO(\GG)$-graded equivariant cohomology of a point, which is distinctly nontrivial
away from the integer-graded part. The $\GG$-cohomology of a point was
calculated by Stong, in an unpublished manuscript, and first published by
Lewis in \cite{Le:projectivespaces}. We summarize the calculation
in \S\ref{sec:pointSummary}, where we also give the cohomologies
of $E\GG$ and $\tE\GG$; part of our calculation of the cohomology of $B_\GG U(1)$
is based on isotropy separation,
that is, on the cofibration sequence $(E\GG)_+\to S^0 \to \tE\GG$.

Part~\ref{part:BU1} of the paper gives the calculation of the cohomology of $B_\GG U(1)$.
The main result is Theorem~\ref{thm:multstructure}, which gives a simple
description of this cohomology as an algebra over the cohomology of a point,
in terms of multiplicative generators and relations.
We also prove directly that it is a free module
over the cohomology of a point, which was expected,
given general freeness results like \cite{FerLew:Freeness} in the $RO(\GG)$-graded case.
The expanded grading includes $RO(\GG)$ as a subgroup and we compare
our results to Lewis's $RO(\GG)$-graded calculation in \S\ref{sec:Lewis},
where we will be able to see his result as a restriction
of ours, but a restriction that does not see the most natural generators.
We also describe the algebra obtained when using other coefficient systems,
including constant $\Z$ coefficients.

Part~\ref{part:basespace} calculates the cohomology of a related space that carries information
about the component structure only.
This calculation helps explain the presence of some of the generators of the
cohomology of $B_\GG U(1)$.
It also helps explain a divisibility phenomenon that shows up in the cohomology
of finite projective spaces.
We use the calculation as well to explain why, equivariantly,
$B_\GG U(1)$ does not also represent an ordinary cohomology group,
as opposed to the convenient fact that, nonequivariantly, $BU(1)$ represents
$H^2(-;\Z)$.

This work owes a large debt to Gaunce Lewis, of course.
I would like to think that he would have enjoyed this
paper, if he hadn't written it himself first.
This paper is also founded on my long and continuing collaboration with Stefan Waner that produced
\cite{CW:ordinaryhomology}, which in turn came out of joint work with Peter May
on equivariant orientation theory, which culminated in \cite{CMW:orientation}.
Finally, a word of praise for serendipity and the MathOverflow website,
where I stumbled across a discussion of Bergman's diamond lemma \cite{Berg:diamondLemma}, which greatly
simplified and improved the line of argument used in Sections~\ref{sec:ringstructure}
and~\ref{sec:structure}.

\part{Equivariant ordinary cohomology}\label{part:cohomology}

\section{The representation ring and the Burnside ring}\label{sec:RepBurnsideRings}

Throughout this paper our ambient group will be $\GG = \{1, t\}$ (written multiplicatively), 
unless we say explicitly that we are stating results for more general groups.
We begin by introducing notations for some representations of $\GG$.

\begin{definition}\label{def:irrRepresentations}
\hspace{2em}
\begin{enumerate}
\item
Let $\R$ denote the one-dimensional (real) representation of $\GG$ with trivial action.

\item
Let $\C$ denote the one-dimensional complex representation of $\GG$ with trivial action.
We fix once and for all an identification $\C = \R^2$.

\item
Let $\R^\sigma$ denote the nontrivial irreducible real representation of $\GG$,
that is, $\R$ with $t$ acting by multiplication by $-1$.
We fix once and for all a nonequivariant identification $\R = \R^\sigma$.

\item
Let $\C^\sigma$ denote the nontrivial complex representation of $\GG$,
that is, $\C$ with $t$ acting by multiplication by $-1$.
With our fixed identification $\C = \R^2$, we also think of $\C^\sigma = (\R^\sigma)^2$
as a real representation.

\end{enumerate}
\end{definition}

Thus, the real representation ring $RO(\GG)$ is a free abelian group on two generators, which
we call $1$ (the class of $\R$) and $\sigma$ (the class of $\R^\sigma$).

\begin{definition}\label{def:ro0G}
If $\alpha\in RO(\GG)$, let $|\alpha|\in \Z$ denote the (nonequivariant) dimension of $\alpha$
and let $\alpha^\GG\in\Z$ denote the dimension of its fixed set.
\end{definition}

The Burnside ring $A(\GG)$ is also free abelian on two generators,
which we call $1$ (the class of the orbit $\GG/\GG$) and $g$ (the class of the orbit $\GG/e$),
with $g^2 = 2g$.
We let $\kappa = 2 - g$, so $\kappa^2 = 2\kappa$ and
$\{1,\kappa\}$ is another basis of $A(\GG)$.
We let $\epsilon\colon A(\GG)\to\Z$ denote the augmentation map,
with $\epsilon(1) = 1$, $\epsilon(g) = 2$,
and $\epsilon(\kappa) = 0$.

Segal \cite{Seg:equivariantstable} showed that, for finite $G$, $A(G)$ is
isomorphic to the ring of stable $G$-maps from $S^0$ to itself.
When $G = \GG$, such a stable map $f$ is determined by its nonequivariant degree
and the degree of its fixed-set map $f^\GG$.
This allows us to explicitly identify maps of spheres with elements of $A(\GG)$.
With $S^\sigma$ denoting the one-point compactification of $\R^\sigma$,
an important example is the map $f\colon S^\sigma\to S^\sigma$
induced by negation on $\R^\sigma$;
$f$ has nonequivariant degree $-1$ and $f^\GG$ has degree 1 (being the identity map on $S^0$).
Thus, $f$ must represent the element $1-g = \kappa-1 \in A(\GG)$, because $\epsilon(1-g) = -1$ and
$(1-g)^\GG = 1$. Notice that $(1-g)^2 = 1$, so $1-g$ is a unit in $A(\GG)$.
In fact, the units in $A(\GG)$ are precisely $\pm 1$ and $\pm (1-g) = \mp (1-\kappa)$.

\section{Some Mackey functors}\label{sec:MackeyFunctors}

Equivariant ordinary cohomology uses Mackey functors as coefficients, and
we will view it as Mackey functor--valued, so we review some basic
facts about such functors and establish some notation.

\begin{definition}
Let $\orb{\GG}$ denote the orbit category of $\GG$ and let $\sorb{\GG}$ denote the stable orbit category,
i.e., the category of orbits of $\GG$ and stable $\GG$-maps between them.
\end{definition}

$\orb\GG$ and $\sorb\GG$ each have two objects, $\GG/\GG$ and $\GG/e$.
We picture the maps as follows:
\[
 \xymatrix{
  \GG/\GG \\
  \GG/e \ar[u]^\rho \ar@(dl,dr)[]_t \\
  \orb{\GG}
 }
\qquad\qquad
 \xymatrix{
  \GG/\GG \ar@(ur,ul)[]_{A(\GG)} \ar@/^/[d]^{\tau} \\
  \GG/e \ar@/^/[u]^{\rho} \ar@(dl,dr)[]_{\Z[\GG]} \\
  \sorb{\GG}
 }
\]
That is, in the stable orbit category, the ring of self maps of $\GG/\GG$ is
$A(\GG)$ while the ring of self maps of $\GG/e$ is isomorphic to 
the group ring $\Z[\GG] \iso \Z[t]/\langle t^2 \rangle$.
The group of maps $\GG/e\to \GG/\GG$ is free abelian on the projection $\rho$ while the group
of maps $\GG/\GG\to \GG/e$ is free abelian on the transfer map $\tau$.
We have $\rho \circ \tau = g$
and $\tau \circ \rho = 1+t$.
Finally, $\rho t = \rho$, $t\tau = \tau$, $g\rho = 2\rho$, and $\tau g = 2\tau$.

\begin{definition}
A {\em Mackey functor} is a contravariant additive functor from the stable orbit category to the
category of abelian groups.
\end{definition}

Following common practice, we will denote Mackey functors using an underline, for example, $\Mackey T$.
In \cite{CW:ordinaryhomology} we used overlines for contravariant functors and underlines for covariant functors.
In the case of finite groups, as here, we need not make the distinction because the stable orbit category
is self-dual.

We will generally picture a Mackey functor $\Mackey T$ using a diagram of the following form:
\[
 \xymatrix{
  \Mackey T(\GG/\GG) \ar@/_/[d]_{\rho} \\
  \Mackey T(\GG/e) \ar@/_/[u]_{\tau} \ar@(dl,dr)[]_{t^*}
 }
\]
Here, $\rho$ and $\tau$ are the maps induced by the maps of the same name in $\sorb \GG$.
$\Mackey T(\GG/\GG)$ should be a module over the Burnside ring; the action is specified by this diagram
because the action of $g$ is given by $\tau\circ\rho$.

We now review and give names to the Mackey functors that will appear in our calculations, beginning
with the following two:
\[
 \Mackey A_{\GG/\GG} = \sorb{\GG}(-,\GG/\GG)\colon \xymatrix{
		A(\GG) \ar@/_/[d]_{\epsilon} \\
		\Z \ar@/_/[u]_{\cdot g} \ar@(dl,dr)[]_{1}
	}
\qquad
 \Mackey A_{\GG/e} = \sorb{\GG}(-,\GG/e)\colon \xymatrix{
		\Z \ar@/_/[d]_{\cdot (1+t)} \\
		\Z[\GG] \ar@/_/[u]_{\epsilon} \ar@(dl,dr)[]_{\cdot t}
	  }
\]
We call $\Mackey A_{\GG/\GG}$ the {\em Burnside ring Mackey functor}
and often write $\Mackey A = \Mackey A_{\GG/\GG}$ for brevity.
In $\Mackey A_{\GG/e}$, 
$\epsilon\colon\Z[\GG]\to\Z$ is
\[
 \epsilon(a_0 + a_1 t) = a_0 + a_1.
\]
$\Mackey A$ and $\Mackey A_{\GG/e}$, being represented functors,
are both projective, with
\[
 \Hom_{\sorb{\GG}}(\Mackey A, \Mackey T) \iso \Mackey T(\GG/\GG) \qquad\text{and}\qquad
 \Hom_{\sorb{\GG}}(\Mackey A_{\GG/e}, \Mackey T) \iso \Mackey T(\GG/e).
\]

The next two are examples of a general construction $\conc C$
for any abelian group $C$, but these are the two cases that will occur in our calculations:
\[
 \conc\Z \colon \xymatrix{
		\Z \ar@/_/[d] \\
		0 \ar@/_/[u] \ar@(dl,dr)[]_{}
	   }
\qquad\qquad
 \conc{\Z/2}\colon \xymatrix{
			\Z/2 \ar@/_/[d] \\
			0 \ar@/_/[u] \ar@(dl,dr)[]
	     }
\]

The last four Mackey functors that will occur are the following.
\[
\begin{array}{rcrc}
 \Mackey\Z\colon & \xymatrix{
		\Z \ar@/_/[d]_{1} \\
		\Z \ar@/_/[u]_{2} \ar@(dl,dr)[]_{1}
	}
&\qquad\qquad
 \Mackey\Z_{-}\colon & \xymatrix{
		0 \ar@/_/[d] \\
		\Z \ar@/_/[u] \ar@(dl,dr)[]_{-1}
	 }
\\ \\
 \Mackey\Z'\colon & \xymatrix{
		\Z \ar@/_/[d]_{2} \\
		\Z \ar@/_/[u]_{1} \ar@(dl,dr)[]_{1}
	}
&\qquad\qquad
 \Mackey\Z'_{-}\colon & \xymatrix{
		\Z/2 \ar@/_/[d]_{0} \\
		\Z \ar@/_/[u]_{\pi} \ar@(dl,dr)[]_{-1}
	 }
\end{array}
\]
$\Mackey\Z$ is the functor usually referred to as {\em constant $\Z$ coefficients}.

\subsection*{Multiplicative structures}

Let $\Mackey S\boxprod \Mackey T$ denote the box product of Mackey functors
as described, for example, in \cite{Le:projectivespaces}.
For us it suffices to know that a map $\Mackey S\boxprod \Mackey T \to \Mackey U$
is equivalent to a pair of maps
\begin{align*}
 \Mackey S(\GG/\GG)\tensor \Mackey T(\GG/\GG) &\to \Mackey U(\GG/\GG) \quad\text{and} \\
 \Mackey S(\GG/e)\tensor \Mackey T(\GG/e) &\to \Mackey U(\GG/e)
\end{align*}
satisfying the following conditions, where we write
$xy$ for the image of $x\tensor y$ under the appropriate one of these maps:
\begin{align*}
 t(xy) &= (tx)(ty), \\
 \rho(xy) &= \rho(x)\rho(y), \\
 \tau(x\rho(y)) &= \tau(x)y, \quad\text{and} \\
 \tau(\rho(x)y) &= x\tau(y).
\end{align*}
The last two conditions are called the {\em Frobenius relations.}
By convention, if $x\in \Mackey S(\GG/\GG)$ and $y\in\Mackey T(\GG/e)$,
we will write $xy$ for $\rho(x)y \in \Mackey U(\GG/e)$.

The functor $\Mackey A$ has a self-pairing 
$\Mackey A\boxprod \Mackey A\to \Mackey A$ using the usual
ring structures on $A(\GG)$ and $\Z$.
A {\em unital ring} is a Mackey functor $\Mackey T$ with an associative and unital
pairing $\Mackey T\boxprod \Mackey T\to \Mackey T$,
where the unit is given by a map $\Mackey A\to \Mackey T$.
(Here, we use that $\Mackey A$ is the unit for $\boxprod$,
meaning that $\Mackey A\boxprod\Mackey T \iso \Mackey T$ for
any Mackey functor $\Mackey T$.)
The conditions above say that this is equivalent to $\Mackey T(\GG/\GG)$ being
a unital ring, $\Mackey T(\GG/e)$ being a unital ring (with the action of $t$ being a ring map), 
$\rho\colon \Mackey T(\GG/\GG)\to \Mackey T(\GG/e)$ being a ring map,
and $\tau\colon\Mackey T(\GG/e)\to\Mackey T(\GG/\GG)$ being a left and right
$\Mackey T(\GG/\GG)$-module map.
Clearly, $\Mackey A$ is itself a unital ring.
Every Mackey functor is a {\em module} over $\Mackey A$ in the obvious sense.

$\Mackey\Z$ is a unital ring with the usual ring structure on $\Z$.
A module over $\Mackey\Z$ is precisely a Mackey functor such that $\tau\circ\rho$
is multiplication by $2$.
The functors
$\Mackey A_{G/e}$, $\Mackey\Z'$, $\Mackey\Z_-$, $\Mackey\Z'_-$, and $\conc{\Z/2}$ 
are all modules over $\Mackey\Z$.
$\Mackey\Z$ and $\Mackey A_{G/e}$ are projective $\Mackey\Z$-modules.

\subsection*{Generators and relations}

We will want to describe the results of our calculations in terms of generators and relations.
When doing so, we identify elements of $\Mackey T(\GG/\GG)$ with maps $\Mackey A\to \Mackey T$
and elements of $\Mackey T(\GG/e)$ with maps $\Mackey A_{\GG/e}\to \Mackey T$.
For example, we can say that $\Mackey\Z$ is generated by an element $\xi$ at level $\GG/\GG$
subject to the relation $\kappa\xi = 0$. By this we mean that 
the following sequence is exact:
\[
 \Mackey A \xrightarrow{\kappa} \Mackey A \xrightarrow{\xi} \Mackey\Z \to 0.
\]
Here, $\kappa$ is the map 
corresponding to $\kappa\in \Mackey A(\GG/\GG)$ and
$\xi$ is the map
corresponding to $1\in \Mackey\Z(\GG/\GG)$, which we are also calling $\xi$
while thinking of it as an abstract generator.

Generators may occur at either level $\GG/\GG$ or level $\GG/e$, and similarly for relations.
Here are descriptions of the other examples in terms of generators and relations:

\begin{itemize}
\item $\conc\Z$: Generated by an element $e$ at level $\GG/\GG$ subject to $\rho(e) = 0$.
That is, there is an exact sequence
\[
 \Mackey A_{G/e} \to \Mackey A \to \conc\Z \to 0
\]
where the first map is specified at level $\GG/e$ by $\epsilon\colon\Z[G]\to\Z$.

\item $\conc{\Z/2}$: Generated by an element $e$ at level $\GG/\GG$ subject to
$\rho(e) = 0$ and $2e = 0$.

\item $\Mackey\Z_-$: Generated by an element $\iota$ at level $\GG/e$
such that $\tau\iota = 0$.

\item $\Mackey\Z'$: Generated by an element $\iota$ at level $\GG/e$
such that $t\iota = \iota$.

\item $\Mackey\Z'_-$: Generated by an element $\iota$ at level $\GG/e$
such that $t\iota = -\iota$.

\end{itemize}

We noted that several of these Mackey functors are modules over the ring $\Mackey\Z$.
We can describe modules over $\Mackey\Z$ in terms of generators and relations as well,
where a generator at level $\GG/\GG$ gives a copy of $\Mackey\Z$ while a generator
at level $\GG/e$ gives a copy of $\Mackey A_{\GG/e}$.
The modules $\Mackey\Z_-$, $\Mackey\Z'$, and $\Mackey\Z'_-$ are described
by the same generators and relations as above.
However, we can simplify the description of $\conc{\Z/2}$:
As an $\Mackey\Z$-module, it is generated by an element $e$ at level $\GG/\GG$
such that $\rho(e) = 0$.

\section{Equivariant ordinary $RO(\Pi B)$-graded cohomology}\label{sec:genCohomology}

In \cite{CW:ordinaryhomology}, Stefan Waner and the author gave a detailed exposition of equivariant
ordinary cohomology graded on ``representations of the fundamental groupoid.''
In this section we review some of the basic definitions and properties.
We assume that $G$ is an arbitrary finite group throughout this section, though
\cite{CW:ordinaryhomology} is written in the more general context of compact Lie groups.

\subsection*{The equivariant fundamental groupoid and its representations}

When $X$ is a $G$-space, we have the following definition,
given originally by tom Dieck \cite{tD:transfgroups}
and used extensively in \cite{CMW:orientation} and \cite{CW:ordinaryhomology}.

\begin{definition}
The {\em equivariant fundamental groupoid of $X$}, denoted $\Pi_G X$
(or $\Pi X$ when $G$ is understood), is the category
whose objects are the $G$-maps $x\colon G/H\to X$ for all the orbits $G/H$ of $G$, 
and whose maps from
$x$ to $y\colon G/K\to X$ are pairs $(\omega, \alpha)$, where
$\alpha\colon G/H\to G/K$ is a $G$-map and $\omega$ is a $G$-homotopy class of paths,
rel endpoints, from $x$ to $y\circ\alpha$.
Composition is induced by composition of maps of orbits and the usual composition
of path classes.
\end{definition}

$\Pi_G$ is a 2-functor, taking $G$-maps to functors
and homotopies to natural isomorphisms.
There is a functor $\pi\colon \Pi_G X\to \orb G$, with
$\pi(x\colon G/H\to X) = G/H$ and $\pi(\omega,\alpha) = \alpha$.
This makes $\Pi_G X$ a {\em bundle of groupoids} over $\orb G$
in the language of~\cite{CMW:orientation}.

\begin{definition}\label{def:virtualcat}
For $n$ an integer, let $v\V_G(n)$ denote the category of 
{\em virtual $n$-dimensional orthogonal bundles over orbits}
as defined in~\cite[\S1.3]{CW:ordinaryhomology}.
Its objects are classes of formal differences $G\times_H V - G\times_H W$
under an appropriate stabilization relation.
Maps are given by maps of virtual bundles covering given maps of orbits.
\end{definition}

We have a functor $\pi\colon v\V_G(n)\to \orb G$ that gives the underlying
orbits and maps of orbits.
This makes $v\V_G(n)$ into a bundle of groupoids over $\orb G$.

\begin{definition}
A {\em virtual $n$-dimensional orthogonal representation of $\Pi_G X$} is a functor
$\Pi_G X\to v\V_G(n)$ over $\orb G$.
These form the {\em category of virtual $n$-dimensional representations of $\Pi_G X$}
when we take morphisms to be natural isomorphisms.
\end{definition}

For each representation $V$ of $G$, there is a representation of $\Pi_G X$
we denote by $\Vrep$, given by taking each $x\colon G/H\to X$ to $G/H\times V$
and each $(\omega,\alpha)$ to $\alpha\times 1$.
This generalizes to virtual representations of $G$ to give
representations $\Vrep\ominus \Wrep$.

Less trivially, if $\xi\colon E\to X$ is a $G$-vector bundle over $X$, there is an associated
representation of $\Pi_G X$, denoted $\xi^*$, given by
$\xi^*(x\colon G/H\to X) = x^*(\xi)$.

\begin{definition}
$RO(\Pi_G X)$, the {\em orthogonal representation ring of $\Pi_G X$}, 
is the ring whose elements are the isomorphism classes of
the virtual orthogonal representations of $\Pi_G X$ of all dimensions.
Addition is given by direct sum of bundles and multiplication by tensor product.
\end{definition}

In particular, $RO(\Pi_G(*)) \iso RO(G)$ when $*$ denotes the one-point $G$-space.
In general, the map $X\to *$ induces a map $RO(G)\to RO(\Pi_G X)$, taking $V\ominus W$
to $\Vrep \ominus \Wrep$ as above.

\subsection*{Ordinary cohomology}

Let $B$ be a $G$-space and let $(X,q,\sigma)$ be an ex-space over $B$, which is to say that
$q\colon X\to B$ is a $G$-map and $\sigma\colon B\to X$ is a section of $q$.
Suppose also given a virtual representation $\gamma$ of $\Pi_G B$ and a Mackey functor
$\Mackey T$. 
In \cite{CW:ordinaryhomology} we defined the
{\em $\gamma$th reduced ordinary cohomology group of $X$ with coefficients in $\Mackey T$,}
$H_G^\gamma(X;\Mackey T)$,
a contravariant functor of
$X$ and a covariant functor of $\gamma$ and of $\Mackey T$.
(Because we will be using only reduced cohomology, we write $H$ rather than $\tilde H$.)

If we have
a parametrized space $q\colon X\to B$, rather than an ex-space, we can form the ex-space
$(X,q)_+$, which we will write $X_+$ when $q$ is understood, given by
$(X\disjunion B, q\disjunion 1, \sigma)$, where
$\sigma\colon B\to X\disjunion B$ is the evident inclusion.
In particular, we may consider $H_G^\gamma(B_+;\Mackey T)$.

The collection of groups $H_G^\gamma(X;\Mackey T)$ as $\gamma$ varies gives the
{\em $RO(\Pi_G B)$-graded ordinary cohomology} of the ex-space $X$ over $B$.
(In \cite{CW:ordinaryhomology} we allow more general coefficient systems than just
Mackey functors, but in this paper we will stick to the simpler case.)
In particular, when we restrict the virtual representation $\gamma$ to
be of the form $\Vrep\ominus\Wrep$, the resulting groups are exactly the
$RO(G)$-graded ordinary cohomology of $X$
discussed in \cite{LMM:roghomology}, \cite{May:alaska}, and \cite{CW:ordinaryhomology},
which generalizes Bredon's original integer-graded theory.
In particular, this theory obeys a dimension axiom that takes the following form:
For $x\colon G/H\to B$ and integers $n$, we have
\[
 H_G^n((G/H,x)_+;\Mackey T) \iso
  \begin{cases}
    \Mackey T(G/H) & \text{if $n=0$} \\
    0 & \text{if $n\neq 0$,}
  \end{cases}
\]
naturally in $x$.

The $RO(\Pi_G B)$-graded theory has many nice properties, discussed in detail
in \cite{CW:ordinaryhomology}. One of the main reasons for introducing the enlarged
grading is to get general Thom isomorphism and Poincar\'e duality theorems.
In particular, the Thom isomorphism takes the following form:
If $\xi\colon E\to B$ is a $G$-vector bundle, let
$T(\xi)$ denote the ex-space over $B$ given by taking the one-point
compactification of each fiber, with the section given by the compactification points.
A {\em Thom class} for $\xi$ is a class $t\in H_G^{\xi^*}(T(\xi);\Mackey A)$
such that, for every $x\colon G/K\to B$, the element
\begin{align*}
 x^*(t) &\in H_G^{x^*\xi^*}(T(x^*\xi);\Mackey A) \\
  &\iso H_G^{x^*\xi^*}(G\times_K S^V;\Mackey A) \\
  &\iso H_K^{V}(S^V;\Mackey A_{K/K}) \\
  &\iso A(K)
\end{align*}
is a generator. 
Here, $x^*\xi$ is the pullback of $\xi$ along $x$, so $x^*\xi = G\times_K V$
for some representation $V$ of $K$; $x^*\xi^*$ denotes the representation of 
$RO(\Pi G/K) \iso RO(K)$ associated to $x^*\xi$.

\begin{theorem}[Thom Isomorphism {\cite[3.11.3]{CW:ordinaryhomology}}]
If $\xi\colon E\to B$ is a $G$-vector bundle, then there exists
a Thom class $t\in H_G^{\xi^*}(T(\xi);\Mackey A)$.
For every Thom class $t$, the map
\[
 t\cup - \colon H_G^\gamma(B_+;\Mackey T) \to H_G^{\gamma+\xi^*}(T(\xi);\Mackey T)
\]
is an isomorphism for every representation $\gamma$ and Mackey functor $\Mackey T$.
\qed
\end{theorem}

As usual, given a Thom class $t(\xi)$ for a $G$-bundle $\xi$, we define
the corresponding {\em Euler class} $e(\xi) \in H_G^{\xi^*}(B_+;\Mackey A)$
to be the restriction of $t(\xi)$ along the zero section $B_+\to T(\xi)$.

For computational purposes, it is useful to view cohomology not as group valued,
but as Mackey functor valued. We let
$\Mackey H_G^\gamma(X;\Mackey T)$ be the Mackey functor defined by
\[
 \Mackey H_G^\gamma(X;\Mackey T)(G/K) = H_G^\gamma(G/K_+\smsh_B X; \Mackey T).
\]
Because $H_G^*(-)$ is stable (it has suspension isomorphisms for all representations),
this is actually a functor on stable maps between orbits, so does define
a Mackey functor. 
On the other hand, the Wirthm\"uller isomorphism
(\cite[3.9.5]{CW:ordinaryhomology}) allows us to write this as
\[
 \Mackey H_G^\gamma(X;\Mackey T)(G/K) \iso H_K^{\gamma|K}(X;\Mackey T|K),
\]
where $\gamma|K$ and $\Mackey T|K$ are the restrictions from $G$ to $K$
defined in the most obvious ways.
Thus, treating ordinary cohomology as a Mackey functor amounts to considering
the cohomologies of $X$ for all subgroups of $G$ simultaneously,
along with the associated restriction and transfer maps.
As often happens, the more structure present, the more limited
the possibilities, hence the easier the computations.

It will be useful to distinguish between, say, cohomology graded on $RO(\Pi_G B)$
and cohomology graded on $RO(G)$.
There are traditional notations like $E^*$ and various {\it ad hoc}
variations like $E^\bullet$ or $E^{\textstyle\star}$ to denote grading on different groups.
We adopt instead the following.

\begin{notation}
If $E$ is an object graded on a group $R$, and $Q\subset R$ is any subset, we write $E^Q$ for the collection
of those parts of $E$ graded on elements in $Q$. In particular, we write $E^R$ for the whole
graded object.
\end{notation}

For example, we will write $\Mackey H_G^{RO(\Pi B)}$ for cohomology graded on $RO(\Pi_G B)$
(for a specified base space $B$) and
$\Mackey H_G^{RO(G)}$ for the part graded only on $RO(G)$.
If $\alpha\in RO(\Pi_G B)$, it is often useful to look at
$\Mackey H_G^{\alpha+RO(G)}$, the part of the cohomology graded on the coset
$\alpha+RO(G)\subset RO(\Pi_G B)$.

We also need change-of-grading, defined as follows.

\begin{definition}
If $E^R$ is an object graded on a group $R$ and $\phi\colon Q\to R$ is a group homomorphism,
we define the object $E^Q$ graded on $Q$ by
$E^q = E^{\phi(q)}$ for $q\in Q$.
\end{definition}

In particular, we will see the following occur in several places in our calculations.

\begin{proposition}\label{prop:kerunits}
If $E^R$ is a graded ring with identity and $\phi\colon Q\to R$ is a group homomorphism,
then $E^Q$ is again a graded ring with identity and, for every $q\in\ker\phi$, there is a
distinguished multiplicative unit $\zeta_q \in E^q$, with $\zeta_0 = 1$ and
$\zeta_q\zeta_{q'} = \zeta_{q+q'}$.
\end{proposition}

\begin{proof}
If $x\in E^q$ and $y\in E^{q'}$ with $q, q'\in Q$, then $x\in E^{\phi(q)}$ and $y\in E^{\phi(q')}$, so
\[
 xy \in E^{\phi(q) + \phi(q')} = E^{\phi(q+q')} = E^{q+q'}.
\]
It is straightforward to check that, with this product, $E^Q$ is a graded unital ring,
with identity $1\in E^0 = E^{\phi(0)}$.

For each $q\in \ker\phi$, we let $\zeta_q \in E^q = E^{\phi(q)} = E^0$ be the element corresponding to
the identity $1\in E^0$. The equality $\zeta_q\zeta_{q'} = \zeta_{q+q'}$ then follows
by definition. In particular, $\zeta_q\zeta_{-q} = \zeta_0 = 1$, so each $\zeta_q$ is a unit.
\end{proof}

As an example, if $B$ is a $G$-space, then
the integer-graded nonequivariant cohomology $H^\Z(B_+;\Z)$ can
be regraded on $RO(G)$ via the dimension homomorphism $RO(G)\to \Z$.
This is precisely $\Mackey H_G^{RO(G)}(B_+;\Mackey A)(G/e)$,
the Mackey-functor valued 
$RO(G)$-graded equivariant cohomology of $B$ evaluated at the orbit $G/e$.
There will be a unit in $H^{RO(G)}(B_+;\Z)$ for every
virtual representation $V-W$ of dimension zero, encoding the fact that
nonequivariant cohomology cannot tell the difference between $V$ and $W$.
Similar units appear in $\Mackey H_G^{RO(\Pi B)}(B_+;\Mackey A)(G/e)$.

\section{The cohomologies of a point, $E\GG$, and $\tE\GG$}\label{sec:pointSummary}

We now return to the specific case of $G = \GG$.
From this point on, all cohomology will be assumed to have coefficients in
$\Mackey A$ unless stated otherwise, and we write
$\Mackey H_\GG^\alpha(X)$ for $\Mackey H_\GG^\alpha(X;\Mackey A)$.

In Part~\ref{part:BU1} we shall calculate the cohomology of $B_\GG U(1)$
as an algebra over the cohomology of a point, 
the latter being the reduced
cohomology of $S^0$. The calculation will use
the cofibration sequence $(E\GG)_+\to S^0 \to \tE\GG$, where $E\GG$ is a nonequivariantly contractible 
free $\GG$-space. To that end, we need explicit descriptions of
the $RO(\GG)$-graded cohomologies of $(E\GG)_+$, $S^0$, and $\tE\GG$, and
the maps in the long exact sequence induced by the cofibration sequence.
The calculation of the $\GG$-cohomology of a point
with Burnside ring coefficients is due to Stong
but first appeared in the literature in \cite{Le:projectivespaces}.
The description we give here is essentially Lewis's, but with some changes in notation.

We first introduce some elements.

\begin{definition}\label{def:iota}
Let
\[
 \iota\in \Mackey H_\GG^{\sigma-1}(S^0)(\GG/e)
\]
be the unit given by Proposition~\ref{prop:kerunits} (and the comments following it) by the fact that
$\sigma-1$ is in the kernel of the dimension homomorphism $RO(\GG)\to\Z$.
Its inverse is an element $\iota^{-1}\in\Mackey H_\GG^{1-\sigma}(S^0)(\GG/e)$
such that
$\iota\cdot \iota^{-1}$ is the identity in $\Mackey H_\GG^0(S^0)(\GG/e)$.\end{definition}

\begin{definition}\label{def:EulerClasses}
Let
\[
 e \in \Mackey H_\GG^\sigma(S^0)(\GG/\GG)
\]
be the Euler class of $\R^\sigma$, that is,
the image of the unit under the map
\[
 H_\GG^0(S^0) \iso H_\GG^\sigma(S^\sigma) \to H_\GG^\sigma(S^0),
\]
where the last map is restriction along the inclusion.
\end{definition}

\begin{theorem}[{\cite[2.1 \& 4.3, with proofs in the Appendix]{Le:projectivespaces}}]\label{thm:evenCohomPoint}
Additively, for $\alpha\in RO(\GG)$,
\[
 \Mackey H_\GG^\alpha(S^0) \iso
  \begin{cases}
   \Mackey A & \text{if $\alpha = 0$} \\
   \Mackey\Z & \text{if $|\alpha| = 0$ and $\alpha^\GG < 0$ is even} \\
   \Mackey\Z_- & \text{if $|\alpha| = 0$ and $\alpha^\GG \leq 1$ is odd} \\
   \Mackey\Z' & \text{if $|\alpha| = 0$ and $\alpha^\GG > 0$ is even} \\
   \Mackey\Z'_- & \text{if $|\alpha| = 0$ and $\alpha^\GG \geq 3$ is odd} \\
   \conc\Z & \text{if $|\alpha| \neq 0$ and $\alpha^\GG = 0$} \\
   \conc{\Z/2} & \text{if $|\alpha| > 0$ and $\alpha^\GG < 0$ is even} \\
   \conc{\Z/2} & \text{if $|\alpha| < 0$ and $\alpha^\GG \geq 3$ is odd} \\
   0 & \text{otherwise.}
  \end{cases}
\]
Multiplicatively, $\Mackey H_\GG^{RO(\GG)}(S^0)$ is a strictly commutative $RO(\GG)$-graded ring,
generated by elements
\begin{align*}
 \iota &\in \Mackey H_\GG^{\sigma-1}(S^0)(\GG/e) \\
 \iota^{-1} &\in \Mackey H_\GG^{1-\sigma}(S^0)(\GG/e) \\
 \xi &\in \Mackey H_\GG^{2(\sigma-1)}(S^0)(\GG/\GG) \\
 e &\in \Mackey H_\GG^{\sigma}(S^0)(\GG/\GG) \\
 e^{-m}\kappa &\in \Mackey H_\GG^{-m\sigma}(S^0)(\GG/\GG) & & m\geq 1 \\
 e^{-m}\delta\xi^{-n} &\in \Mackey H_\GG^{1 - m\sigma - 2n(\sigma-1)}(S^0)(\GG/\GG)
   & & m, n \geq 1.
\end{align*}
These generators satisfy the following {\em structural} relations:
\begin{align*}
 \tau(\iota^{-1}) &= 0 \\
 \tau(\iota^{-2n-1}) &= e^{-1}\delta\xi^{-n} & & \text{for $n\geq 1$} \\
 \kappa\xi &= 0 \\
 \rho(\xi) &= \iota^2 \\
 \rho(e) &= 0 \\
 \rho(e^{-m}\kappa) &= 0 & & \text{for $m\geq 1$}\\
 \rho(e^{-m}\delta\xi^{-n}) &= 0 & & \text{for $m\geq 2$ and $n\geq 1$}\\
 2e^{-m}\delta\xi^{-n} &= 0 & & \text{for $m\geq 2$ and $n\geq 1$}
\intertext{and the following {\em multiplicative} relations:}
  \iota\cdot \iota^{-1} &= \rho(1) \\
  e\cdot e^{-m}\kappa &= e^{-m+1}\kappa & &\text{for $m\geq 1$} \\
 \xi\cdot e^{-m}\kappa &= 0 & & \text{for $m\geq 1$} \\
 e^{-m}\kappa\cdot e^{-n}\kappa &= 2e^{-m-n}\kappa & & \text{for $m\geq 0$ and $n\geq 0$} \\
 e\cdot e^{-m}\delta\xi^{-n} &= e^{-m+1}\delta\xi^{-n} & &\text{for $m\geq 2$ and $n\geq 1$} \\
 \xi\cdot e^{-m}\delta\xi^{-n} &= e^{-m}\delta\xi^{-n+1} & & 
   \text{for $m\geq 1$ and $n\geq 2$} \\
 \xi\cdot e^{-m}\delta\xi^{-1} &= 0 & & \text{for $m\geq 2$} \\
 e^{-m}\kappa \cdot e^{-n}\delta\xi^{-k} &= 0 & & \text{if $m\geq 0$, $n\geq 1$, and $k\geq 1$}\\
 e^{-m}\delta\xi^{-k}\cdot e^{-n}\delta\xi^{-\ell} &= 0 & & 
   \text{if $m, n, k, \ell\geq 1$}
\end{align*}
The following relations are implied by the preceding ones:
\begin{align*}
 \kappa e &= 2e \\
 2e^m\xi^n &= 0 & & \text{if $m>0$ and $n>0$} \\
 t\iota^k &= (-1)^k\iota^k & & \text{for all $k$} \\
 \xi\cdot\tau(\iota^k) &= \tau(\iota^{k+2}) & &\text{for all $k$} \\
 e\cdot\tau(\iota^k) &= 0 & &\text{for all $k$} \\
 e^{-m}\kappa \cdot \tau(\iota^k) &= 0 & & \text{for all $m\geq 1$ and $k$} \\
 e^{-m}\delta\xi^{-n}\cdot \tau(\iota^k) &= 0 & & \text{for all $m, n\geq 1$ and $k$} \\
 \tau(\iota^k)\cdot\tau(\iota^\ell) &= 0 & &\text{if $k$ or $\ell$ is odd} \\
 \tau(\iota^{2k})\cdot\tau(\iota^{2\ell}) &= 2\tau(\iota^{2(k+\ell)}) & &\text{for all $k$ and $\ell$}\\
 \tau(\iota^{2k+1}) &= 0 & & \text{if $k\geq 0$} \\
 e\cdot e^{-1}\delta\xi^{-n} &= 0 & & \text{if $n\geq 1$} 
\end{align*}
\qed
\end{theorem}

The notation $e^{-m}\kappa$ comes from the fact that 
$e^m\cdot e^{-m}\kappa = \kappa \in H_\GG^0(S^0) = A(\GG)$.
The reason for the notation $e^{-m}\delta\xi^{-n}$ should become clearer shortly.

It helps to have a way to visualize these calculations, and the common
way of doing so is to plot the groups or Mackey functors on a grid.
Different authors, however, have used different axes.
Lewis, in \cite{Le:projectivespaces}, uses
$\alpha^\GG$ as the horizontal axis and $|\alpha|$ as the vertical axis.
Dugger \cite{Dug:AHSSforKR} and Kronholm \cite{Kro:SerreSS} use the
so-called motivic grading, plotting $|\alpha|$ vs $|\alpha| - \alpha^\GG$,
although Dugger uses $|\alpha|$ as the vertical axis while
Kronholm uses it as the horizontal axis.
Thinking of $\alpha = a + b\sigma$, we will use $a = \alpha^\GG$ for the horizontal axis
and $b = |\alpha| - \alpha^\GG$ for the vertical axis,
consistent with the diagrams given in \cite{CHT:projspaces}.

Figure~\ref{fig:EvenCohomPoint} shows the Mackey functors $\Mackey H_\GG^\alpha(S^0)$, 
with dots representing zero functors.
Figure~\ref{fig:EvenCohomPointGens} gives the elements that generate the corresponding Mackey functors.
Generators shown in parentheses represent elements at level $\GG/e$.
\begin{figure}
\[\def\objectstyle{\scriptstyle}
 \xymatrix@!0@R=5ex@C=2.5em{
  & & & & &\phantom{a}\ b\sigma & & & & & \\
  & \conc{\Z/2} &\cdot& \conc{\Z/2} &\cdot& \conc\Z &\cdot& \cdot &\cdot& \cdot & \cdot \\
  & \Mackey\Z &\cdot& \conc{\Z/2} &\cdot& \conc\Z &\cdot& \cdot &\cdot& \cdot & \cdot \\
  & \cdot &\Mackey\Z_{\mathrlap{-}}& \conc{\Z/2} &\cdot& \conc\Z &\cdot& \cdot &\cdot& \cdot & \cdot \\
  & \cdot &\cdot& \Mackey\Z &\cdot& \conc\Z &\cdot& \cdot &\cdot& \cdot & \cdot \\
  & \cdot &\cdot& \cdot &\Z_{\mathrlap{-}}& \conc\Z &\cdot& \cdot &\cdot& \cdot & \cdot \\
  \ar@{-}'[rrrrr][rrrrrrrrrrr]
   &  &  &  &  & \Mackey A &  &  &  &  &  & a\\
  & \cdot &\cdot& \cdot &\cdot& \conc\Z &\Mackey \Z_{\mathrlap{-}}& \cdot & \cdot & \cdot & \cdot  \\
  & \cdot &\cdot& \cdot &\cdot& \conc\Z &\cdot& \Mackey \Z' & \cdot & \cdot & \cdot  \\
  & \cdot &\cdot& \cdot &\cdot& \conc\Z &\cdot& \cdot & \Mackey \Z'_{\mathrlap{-}} & \cdot & \cdot  \\
  & \cdot &\cdot& \cdot &\cdot& \conc\Z &\cdot& \cdot & \conc{\Z/2} & \Mackey \Z' &  \cdot  \\
  & \cdot &\cdot& \cdot &\cdot& \conc\Z &\cdot& \cdot & \conc{\Z/2} & \cdot &  \Mackey \Z'_{\mathrlap{-}}  \\
  & \cdot &\cdot& \cdot &\cdot& \conc\Z &\cdot& \cdot & \conc{\Z/2} & \cdot &  \conc{\Z/2}  \\
  & & & & & \ar@{-}'[u]'[uu]'[uuu]'[uuuu]'[uuuuu]'[uuuuuu]'[uuuuuuu]'[uuuuuuuu]'[uuuuuuuuu]'[uuuuuuuuuu]'[uuuuuuuuuuu]'[uuuuuuuuuuuu][uuuuuuuuuuuuu]
 }
\]
\caption{$\protect\Mackey H_\GG^{RO(\GG)}(S^0)$ at $a+b\sigma$}\label{fig:EvenCohomPoint}
\end{figure}
\begin{figure}
\[\def\objectstyle{\scriptstyle}
 \xymatrix@!0@R=5ex@C=3em{
  & & & & & \phantom{a}\  b\sigma & & & & & \\
  & e\xi^2 &\cdot& e^3\xi &\cdot& e^5 &\cdot& \cdot &\cdot& \cdot & \cdot \\
  & \xi^2 &\cdot& e^2\xi &\cdot& e^4 &\cdot& \cdot &\cdot& \cdot & \cdot \\
  & \cdot &(\iota^3) & e\xi &\cdot& e^3 &\cdot& \cdot &\cdot& \cdot & \cdot \\
  & \cdot &\cdot& \xi &\cdot& e^2 &\cdot& \cdot &\cdot& \cdot & \cdot \\
  & \cdot &\cdot& \cdot &(\iota)& e &\cdot& \cdot &\cdot& \cdot & \cdot \\
  \ar@{-}'[rrrrr][rrrrrrrrrrr]
   &  &  &  &  & 1 &  &  &  &  &  & a\\
  & \cdot &\cdot& \cdot &\cdot& e^{-1}\kappa &(\iota^{-1})& \cdot & \cdot & \cdot &  \cdot  \\
  & \cdot &\cdot& \cdot &\cdot& e^{-2}\kappa &\cdot& (\iota^{-2}) & \cdot & \cdot & \cdot  \\
  & \cdot &\cdot& \cdot &\cdot& e^{-3}\kappa &\cdot& \cdot & (\iota^{-3}) & \cdot &  \cdot  \\
  & \cdot &\cdot& \cdot &\cdot& e^{-4}\kappa &\cdot& \cdot & e^{-2}\delta\xi^{-1} & (\iota^{-4}) &  \cdot  \\
  & \cdot &\cdot& \cdot &\cdot& e^{-5}\kappa &\cdot& \cdot & e^{-3}\delta\xi^{-1} & \cdot &  (\iota^{-5})  \\
  & \cdot &\cdot& \cdot &\cdot& e^{-6}\kappa &\cdot& \cdot & e^{-4}\delta\xi^{-1} & \cdot &  e^{-2}\delta\xi^{-2}  \\
  & & & & & \ar@{-}'[u]'[uu]'[uuu]'[uuuu]'[uuuuu]'[uuuuuu]'[uuuuuuu]'[uuuuuuuu]'[uuuuuuuuu]'[uuuuuuuuuu]'[uuuuuuuuuuu]'[uuuuuuuuuuuu]'[uuuuuuuuuuuuu]
 }
\]
\caption{The generators of $\protect\Mackey H_\GG^{RO(\GG)}(S^0)$}\label{fig:EvenCohomPointGens}
\end{figure}

\begin{remark}[Relationship with other common notations]\label{rem:othernotations}

The related and somewhat simpler ring $\Mackey H_{C_2}^{RO(C_2)}(S^0;\Mackey \Z)$
has been used extensively in recent literature. 
The map $\Mackey A\to \Mackey \Z$ that is $\epsilon$ at level $C_2/C_2$
induces a map
\[
	\epsilon_*\colon \Mackey H_{C_2}^{RO(C_2)}(S^0;\Mackey A) \to \Mackey H_{C_2}^{RO(C_2)}(S^0;\Mackey \Z).
\]
Elements of $\Mackey H_{C_2}^{RO(C_2)}(S^0;\Mackey \Z)$ have been given
a multitude of different names:
Under $\epsilon_*$, $e$ maps to the element called
$a$ in \cite{HuKr:adamsnovikov}, $a_\sigma$ in \cite{HHR:Kervaire}, 
and $\rho$ in the motivic literature, for example in \cite{Dug:Grassmannians}
(in which the coefficients are further reduced to the constant $\Z/2$ functor).
The element $\xi$ maps to the element called $u_{2\sigma}$ in \cite{HHR:Kervaire} and
the element $\tau^2$ in \cite{Dug:Grassmannians},
whereas $\tau(\iota^{-2})$ maps to the element $\theta$ in \cite{Dug:Grassmannians}.
In \cite{Dug:AHSSforKR}, Dugger uses constant $\Z$ coefficients but a completely different
set of names. The element he calls $\theta$ there is the image of our $\tau(\iota^{-3})$,
while the image of $\tau(\iota^{-2})$ is called $\alpha$.
\end{remark}

We will also need to know the $RO(\GG)$-graded cohomology of $E\GG$.

\begin{theorem}\label{thm:EvenEG}
Additively, for $\alpha\in RO(\GG)$,
\[
 \Mackey H_\GG^{\alpha}((E\GG)_+) \iso
  \begin{cases}
   \Mackey\Z & \text{if $|\alpha| = 0$ and $\alpha^\GG$ is even} \\
   \Mackey\Z_- & \text{if $|\alpha| = 0$ and $\alpha^\GG$ is odd} \\
   \conc{\Z/2} & \text{if $|\alpha|>0$ and $\alpha^\GG$ is even} \\
   0 & \text{otherwise.}
  \end{cases}
\]
Multiplicatively, $\Mackey H_\GG^{RO(\GG)}((E\GG)_+)$ is a strictly commutative $RO(\GG)$-graded
$\Mackey\Z$-algebra generated by
\begin{align*}
 e &\in \Mackey H_\GG^{\sigma}((E\GG)_+)(\GG/\GG), \\
 \iota &\in \Mackey H_\GG^{\sigma-1}((E\GG)_+)(\GG/e), \\
 \iota^{-1} &\in \Mackey H_\GG^{1-\sigma}((E\GG)_+)(\GG/e), \\
 \xi &\in \Mackey H_\GG^{2(\sigma-1)}((E\GG)_+)(\GG/\GG), \quad\text{and} \\
 \xi^{-1} &\in \Mackey H_\GG^{2(1-\sigma)}((E\GG)_+)(\GG/\GG),
\end{align*}
subject to the relations
\begin{align*}
 \rho(e) &= 0 \\
 \tau(\iota) &= 0 \\
 \rho(\xi) &= \iota^2 \\
 \iota\cdot \iota^{-1} &= \rho(1) \quad\text{and} \\
 \xi\cdot\xi^{-1} &= 1.
\end{align*}
\end{theorem}

Figures~\ref{fig:EvenEG} and~\ref{fig:EvenEGGens} show $\Mackey H_\GG^{RO(\GG)}((E\GG)_+)$ and its generators.
It follows from the calculation, and
the action of $\Mackey H_\GG^{RO(\GG)}(S^0)$ given by
Proposition~\ref{prop:EvenMaps}, that
$
 \Mackey H_\GG^{RO(\GG)}((E\GG)_+) \iso \Mackey H_\GG^{RO(\GG)}(S^0)[\xi^{-1}]
$.
\begin{figure}
\[\def\objectstyle{\scriptstyle}
 \xymatrix@!0@R=5ex@C=2.5em{
  & & & & &\phantom{a}\ b\sigma & & & & \\
  & \conc{\Z/2} &\cdot& \conc{\Z/2} &\cdot& \conc{\Z/2} &\cdot& \conc{\Z/2} &\cdot& \conc{\Z/2}  \\
  & \Mackey\Z &\cdot& \conc{\Z/2} &\cdot& \conc{\Z/2} &\cdot& \conc{\Z/2} &\cdot& \conc{\Z/2} \\
  & \cdot &\Mackey\Z_{\mathrlap{-}}& \conc{\Z/2} &\cdot& \conc{\Z/2} &\cdot& \conc{\Z/2} &\cdot& \conc{\Z/2}  \\
  & \cdot &\cdot& \Mackey\Z &\cdot& \conc{\Z/2} &\cdot& \conc{\Z/2} &\cdot& \conc{\Z/2}  \\
  & \cdot &\cdot& \cdot &\Z_{\mathrlap{-}}& \conc{\Z/2} &\cdot& \conc{\Z/2} &\cdot& \conc{\Z/2}  \\
  \ar@{-}'[rrrrr]'[rrrrrrr]'[rrrrrrrrr][rrrrrrrrrr]
   &  &  &  &  & \Mackey\Z &  & \conc{\Z/2} &  & \conc{\Z/2} & a\\
  & \cdot &\cdot& \cdot &\cdot& \cdot &\Mackey \Z_{\mathrlap{-}}& \conc{\Z/2} & \cdot & \conc{\Z/2}  \\
  & \cdot &\cdot& \cdot &\cdot& \cdot &\cdot& \Mackey \Z & \cdot & \conc{\Z/2}   \\
  & \cdot &\cdot& \cdot &\cdot& \cdot &\cdot& \cdot & \Mackey \Z_{\mathrlap{-}} & \conc{\Z/2}   \\
  & \cdot &\cdot& \cdot &\cdot& \cdot &\cdot& \cdot & \cdot & \Mackey \Z   \\
  & & & & & \ar@{-}'[u]'[uu]'[uuu]'[uuuu]'[uuuuu]'[uuuuuu]'[uuuuuuu]'[uuuuuuuu]'[uuuuuuuuu]'[uuuuuuuuuu]'[uuuuuuuuuuu]
 }
\]
\caption{$\protect\Mackey H_\GG^{RO(\GG)}((E\GG)_+)$}\label{fig:EvenEG}
\end{figure}
\begin{figure}
\[\def\objectstyle{\scriptstyle}
 \xymatrix@!0@R=5ex@C=2.5em{
  & & & & &\phantom{a}\ b\sigma & & & & \\
  & e\xi^2 &\cdot& e^3\xi &\cdot& e^5 &\cdot& e^7\xi^{-1} &\cdot& e^9\xi^{-2}  \\
  & \xi^2 &\cdot& e^2\xi &\cdot& e^4 &\cdot& e^6\xi^{-1} &\cdot& e^8\xi^{-2} \\
  & \cdot & (\iota^3) & e\xi &\cdot& e^3 &\cdot& e^5\xi^{-1} &\cdot& e^7\xi^{-2}  \\
  & \cdot &\cdot& \xi &\cdot& e^2 &\cdot& e^4\xi^{-1} &\cdot& e^6\xi^{-2}  \\
  & \cdot &\cdot& \cdot & (\iota) & e &\cdot& e^3\xi^{-1} &\cdot& e^5\xi^{-2}  \\
  \ar@{-}'[rrrrr]'[rrrrrrr]'[rrrrrrrrr][rrrrrrrrrr]
   &  &  &  &  & 1 &  & e^2\xi^{-1} &  & e^4\xi^{-2} & a\\
  & \cdot &\cdot& \cdot &\cdot& \cdot & (\iota^{-1}) & e\xi^{-1} & \cdot & e^3\xi^{-2}  \\
  & \cdot &\cdot& \cdot &\cdot& \cdot &\cdot& \xi^{-1} & \cdot & e^2\xi^{-2}   \\
  & \cdot &\cdot& \cdot &\cdot& \cdot &\cdot& \cdot & (\iota^{-3}) & e\xi^{-2}   \\
  & \cdot &\cdot& \cdot &\cdot& \cdot &\cdot& \cdot & \cdot & \xi^{-2}   \\
  & & & & & \ar@{-}'[u]'[uu]'[uuu]'[uuuu]'[uuuuu]'[uuuuuu]'[uuuuuuu]'[uuuuuuuu]'[uuuuuuuuu]'[uuuuuuuuuu]'[uuuuuuuuuuu]
 }
\]
\caption{The generators of $\protect\Mackey H_\GG^{RO(\GG)}((E\GG)_+)$}\label{fig:EvenEGGens}
\end{figure}

\begin{proof}[Proof of Theorem~\ref{thm:EvenEG}]
Because $E\GG$ is free, we have in integer gradings that
\begin{align*}
 H_\GG^n((E\GG)_+) &\iso H^n((E\GG/\GG)_+;\Z) = H^n((B\GG)_+;\Z) \\
  &= \begin{cases}
  		\Z & \text{if $n = 0$} \\
		\Z/2 & \text{if $n>0$ is even} \\
		0 & \text{otherwise.}
  	 \end{cases}
\end{align*}
Because $E\GG$ is nonequivariantly contractible, we have
\[
 H^n((E\GG)_+;\Z) \iso H^n(S^0;\Z)
  = \begin{cases}
  		\Z & \text{if $n=0$} \\
		0 & \text{otherwise.}
    \end{cases}
\]
The map $\GG\times E\GG \to E\GG$ induces the restriction map from
$H_\GG^n((E\GG)_+)$ to $H^n((E\GG)_+;\Z) \iso H^n(S^0;\Z)$; on taking orbits by $\GG$ this
is the map induced by the nonequivariant map
$S^0 \hmtpc (E\GG)_+\to (B\GG)_+$, which leads to the calculation
\[
 \Mackey H_\GG^n((E\GG)_+)
  \iso \begin{cases}
  			\Mackey \Z & \text{if $n = 0$} \\
			\conc{\Z/2} & \text{if $n>0$ is even} \\
			0 & \text{otherwise.}
  	   \end{cases}
\]

Now consider, again for $n$ an integer,
\begin{align*}
 H_\GG^{n-\sigma+1}((E\GG)_+) &\iso H^{n+1}((E\GG)_+\smsh_\GG S^\sigma;\Z) \iso H^{n}((B\GG)_+;\tilde\Z) \\
  &= \begin{cases}
  		\Z/2 & \text{if $n > 0$ is odd} \\
		0 & \text{otherwise,}
     \end{cases}
\end{align*}
where $\tilde\Z$ is the nontrivial twisted coefficient system on $B\GG$.
On the other hand,
\[
 H_\GG^{n-\sigma+1}((\GG\times E\GG)_+) \iso H^{n+1}((\GG)_+\smsh_\GG S^\sigma;\Z)
  = \begin{cases}
  		\Z & \text{if $n=0$} \\
		0 & \text{otherwise}
    \end{cases}
\]
with the induced action of $\Z[\GG]$ on $H^{1}((\GG)_+\smsh_\GG S^\sigma;\Z) = \Z$ 
being the nontrivial one. This gives the calculation
\[
 \Mackey H_\GG^{n-\sigma+1}
 \iso \begin{cases}
 		\Mackey\Z_- & \text{if $n=0$} \\
		\conc{\Z/2} & \text{if $n>0$ is odd.}
 	  \end{cases}
\]

To calculate the groups in the remaining gradings, 
consider the equivariant bundle $q\colon E\GG\times \R^{2}\to E\GG$.
Because $E\GG$ is free and the action of $\GG$ on $\R^{2\sigma}$ preserves nonequivariant orientation,
$q$ can be considered an $\R^{2\sigma}$-dimensional bundle in the sense of \cite{CMW:orientation}.
By \cite{CW:ordinaryhomology}, it has a Thom class
\[
 \xi \in H_\GG^{2\sigma}(\susp^2(E\GG)_+) \iso H_\GG^{2\sigma-2}((E\GG)_+),
\]
the choice of which is determined by our nonequivariant identification of $\R^{2\sigma}$ and $\R^2$,
and multiplication by $\xi$ gives the Thom isomorphism
\[
 -\cup \xi\colon H_\GG^{\alpha}((E\GG)_+) \xrightarrow{\iso} H_\GG^{\alpha+2\sigma}(\susp^2(E\GG)_+)
  \iso H_G^{\alpha+2\sigma-2}((E\GG)_+).
\]
Because multiplication by $\xi$ is an isomorphism, $\xi$ is an invertible element of
$H_\GG^{RO(\GG)}((E\GG)_+)$.
This allows us to extend the calculations above to all gradings in $RO(\GG)$,
giving the additive part of the theorem.

To complete the multiplicative structure, we have just constructed the element $\xi$
and we let $e$ and $\iota$ be the images of the elements of the same name
from the cohomology of a point, restricted along the map $E\GG\to *$.
Because $E\GG\to *$ is a nonequivariant equivalence, the (positive and negative)
powers of $\iota$ generate $\Mackey H_\GG^{RO(\GG)}((E\GG)_+)(\GG/E)$
just as they do in the cohomology of a point.
Because we constructed $\xi$ using our chosen nonequivariant identification of $\R^{2\sigma}$
with $\R^2$, we have that $\rho(\xi) = \iota^2$.
Since the same is true in the cohomology of a point, we must have that
the element named $\xi$ in the cohomology of a point maps to the one of the same name
in the cohomology of $E\GG$. 
Because $e$ is the Euler class of $\R^\sigma$,
our calculation of $H_\GG^\sigma((E\GG)_+)$ shows that $e$ must be the nonzero element there.
The various relations can now be seen as
inherited from the cohomology of a point.
\end{proof}

Turning to $\Mackey H_\GG^{RO(\GG)}(\tE\GG)$, its ring structure is not that useful to us, partly because
it is a ring without an identity. We will, instead, describe it as a module
over $\Mackey H_\GG^{RO(\GG)}(S^0)$. In fact, it is a module over the following localization.

\begin{proposition}\label{prop:localizedPoint}
On inverting $e$ in $\Mackey H_\GG^{RO(\GG)}(S^0)$ we get
\[
 \Mackey H_\GG^{RO(\GG)}(S^0)[e^{-1}] \iso
  \conc{\Z}[e,e^{-1},\xi]/\langle 2\xi\rangle.
\]
\end{proposition}

\begin{proof}
Because $\rho(e) = 0$, inverting $e$ kills $\Mackey H_\GG^{RO(\GG)}(S^0)(\GG/e)$,
hence the result is a module over $\conc{\Z}$.
In $\Mackey H_\GG^{RO(\GG)}(S^0)$, every element is annihilated by a high enough power of $e$
except the terms $e^{-m}\kappa$ and $e^m\xi^n$ for $m\geq 0$ and $n\geq 0$.
We have $e^{m+1}\cdot e^{-m}\kappa = 2e$, so $e^{-m}\kappa = 2e^{-m}$
in $\Mackey H_\GG^{RO(\GG)}(S^0)[e^{-1}]$.
We have $2e\xi = 0$, so $2\xi = 0$ in $\Mackey H_\GG^{RO(\GG)}(S^0)[e^{-1}]$.
\end{proof}

\begin{theorem}\label{thm:EvenTEG}
Additively, for $\alpha\in RO(\GG)$,
\[
 \Mackey H_\GG^\alpha(\tE\GG) \iso
  \begin{cases}
    \conc{\Z} & \text{if $\alpha^\GG = 0$} \\
    \conc{\Z/2} & \text{if $\alpha^\GG\geq 3$ is odd} \\
    0 & \text{otherwise.}
  \end{cases}
\]
Multiplication by $e\in\Mackey H_\GG^{RO(\GG)}(S^0)$ is an isomorphism on $\Mackey H_\GG^{RO(\GG)}(\tE\GG)$,
so the latter is a module over $\Mackey H_\GG^{RO(\GG)}(S^0)[e^{-1}]$.
As such, it is generated by elements
\begin{align*}
 \kappa &\in \Mackey H_\GG^0(\tE\GG)(\GG/\GG) \quad\text{and}\\
 \delta\xi^{-k} &\in \Mackey H_\GG^{1 - 2k(\sigma-1)}(\tE\GG)(\GG/\GG) \quad k \geq 1
\end{align*}
such that 
\begin{align*}
 \xi\cdot\kappa &= 0 \\
 \xi\cdot\delta\xi^{-k} &= \delta\xi^{-(k-1)}\quad\text{$k>1$, and} \\
 \xi\cdot\delta\xi^{-1} &= 0.
\end{align*}
\qed
\end{theorem}

Figures~\ref{fig:EvenTEG} and~\ref{fig:EvenTEGGens} show $\Mackey H_\GG^{RO(\GG)}(\tE\GG)$ and its generators.
\begin{figure}
\[\def\objectstyle{\scriptstyle}
 \xymatrix@!0@R=5ex@C=2.5em{
  & \phantom{a}\ b\sigma \ar@{-}'[d]'[dd]'[ddd]'[dddd]'[ddddd]'[dddddd]'[ddddddd][dddddddd] \\
  & \conc\Z & \cdot & \cdot & \conc{\Z/2} & \cdot & \conc{\Z/2} & \cdot & \conc{\Z/2} \\
  & \conc\Z & \cdot & \cdot & \conc{\Z/2} & \cdot & \conc{\Z/2} & \cdot & \conc{\Z/2} \\
  & \conc\Z & \cdot & \cdot & \conc{\Z/2} &  & \conc{\Z/2} &  & \conc{\Z/2} \\
  \ar@{-}'[r]'[rrrr]'[rrrrrr]'[rrrrrrrr][rrrrrrrrr]
   & \conc\Z & & & \conc{\Z/2} & & \conc{\Z/2} & & \conc{\Z/2} &  a \\
  & \conc\Z & \cdot & \cdot & \conc{\Z/2} & \cdot & \conc{\Z/2} & \cdot & \conc{\Z/2} \\
  & \conc\Z & \cdot & \cdot & \conc{\Z/2} & \cdot & \conc{\Z/2} & \cdot & \conc{\Z/2} \\
  & \conc\Z & \cdot & \cdot & \conc{\Z/2} & \cdot & \conc{\Z/2} & \cdot & \conc{\Z/2} \\
  & &  
 }
\]
\caption{$\protect\Mackey H_\GG^{RO(\GG)}(\tE\GG)$}\label{fig:EvenTEG}
\end{figure}
\begin{figure}
\[\def\objectstyle{\scriptstyle}
 \xymatrix@!0@R=5ex@C=2.5em{
  & \phantom{a}\ b\sigma \ar@{-}'[d]'[dd]'[ddd]'[dddd]'[ddddd]'[dddddd]'[ddddddd][dddddddd] \\
  & e^3\kappa & \cdot & \cdot & e^5\delta\xi^{-1} & \cdot & e^7\delta\xi^{-2} & \cdot & e^9\delta\xi^{-3} \\
  & e^2\kappa & \cdot & \cdot & e^4\delta\xi^{-1} & \cdot & e^6\delta\xi^{-2} & \cdot & e^8\delta\xi^{-3} \\
  & e\kappa & \cdot & \cdot & e^3\delta\xi^{-1} & \cdot & e^5\delta\xi^{-2} & \cdot & e^7\delta\xi^{-3} \\
  \ar@{-}'[r]'[rrrr]'[rrrrrr]'[rrrrrrrr][rrrrrrrrr]
   & \kappa & & & e^{2}\delta\xi^{-1} & & e^{4}\delta\xi^{-2} & & e^{6}\delta\xi^{-3} & a  \\
  & e^{-1}\kappa & \cdot & \cdot & e\delta\xi^{-1} & \cdot & e^{3}\delta\xi^{-2} & \cdot & e^{5}\delta\xi^{-3} \\
  & e^{-2}\kappa & \cdot & \cdot & \delta\xi^{-1} & \cdot & e^{2}\delta\xi^{-2} & \cdot & e^{4}\delta\xi^{-3} \\
  & e^{-3}\kappa & \cdot & \cdot & e^{-1}\delta\xi^{-1} & \cdot & e\delta\xi^{-2} & \cdot & e^{3}\delta\xi^{-3} \\
  & &  
 }
\]
\caption{The generators of $\protect\Mackey H_\GG^{RO(\GG)}(\tE\GG)$}\label{fig:EvenTEGGens}
\end{figure}

\begin{proof}[Proof of Theorem~\ref{thm:EvenTEG}]
Consider the long exact sequence
\[
 \cdots \to \Mackey H_\GG^{n-1}((E\GG)_+)
  \xrightarrow{\delta} \Mackey H_\GG^{n}(\tE\GG) \\
  \xrightarrow{\psi} \Mackey H_\GG^{n}(S^0)
  \xrightarrow{\phi} \Mackey H_\GG^{n}((E\GG)_+)
  \to \cdots
\]
in integer grading. From the dimension axiom, we get that
$\delta\colon \Mackey H_\GG^{n-1}((E\GG)_+) \to \Mackey H_\GG^{n}(\tE\GG)$
is an isomorphism for $n\neq 0$ or $1$.
On the other hand, we know that $\Mackey H_\GG^{-1}((E\GG)_+) = 0$
(which is true for any space) and 
$\phi\colon \Mackey H_\GG^{0}(S^0) \to H_\GG^{0}((E\GG)_+)$
is the epimorphism $\epsilon\colon \Mackey A\to \Mackey\Z$. Therefore, we have the short exact sequence
\[
 \xymatrix@R1em{
 	0 \ar[r] & \Mackey H_\GG^{0}(\tE\GG) \ar[r]
		& \Mackey H_\GG^{0}(S^0) \ar@{=}[d] \ar[r] & \Mackey H_\GG^{0}((E\GG)_+) \ar@{=}[d] \ar[r]
		& 0 \\
	& & \Mackey A \ar[r]^{\epsilon} & \Mackey\Z
 }
\]
The kernel of $\epsilon$ is the copy of $\conc\Z$ generated by $\kappa$.
This gives us the additive calculation of $\Mackey H_\GG^{n}(\tE\GG)$ for $n\in\Z$ and also
identifies the generators as $\kappa$ when $n=0$ and
$\delta(e^{n-1}\xi^{(n-1)/2}) = e^{n-1}\delta\xi^{(n-1)/2}$ when $n\geq 3$ is odd.

Now $S^0\to S^\sigma$ is an equivalence on taking fixed points by $\GG$, hence
the map $\tE\GG \to \tE\GG\smsh S^\sigma$ obtained by smashing with $\tE\GG$ is a $\GG$-equivalence.
The induced map in cohomology is given by multiplication by $e$, so
multiplication by $e$ is an isomorphism on $\Mackey H_\GG^{RO(\GG)}(\tE\GG)$.
The theorem follows.
\end{proof}

To complete the picture, we need to describe the maps in the $RO(\GG)$-graded
long exact sequence
\begin{multline*}
 \cdots \to \susp\Mackey H_\GG^{RO(\GG)}((E\GG)_+)
  \xrightarrow{\delta} \Mackey H_\GG^{RO(\GG)}(\tE\GG) \\
  \xrightarrow{\psi} \Mackey H_\GG^{RO(\GG)}(S^0)
  \xrightarrow{\phi} \Mackey H_\GG^{RO(\GG)}((E\GG)_+)
  \to \cdots.
\end{multline*}
It helps to look at Figures~\ref{fig:EvenCohomPoint}--\ref{fig:EvenTEGGens}
when reading the following result, which follows from the proofs of the
computations done above.

\begin{proposition}\label{prop:EvenMaps}
$\delta\colon \susp\Mackey H_\GG^{RO(\GG)}((E\GG)_+)\to \Mackey H_\GG^{RO(\GG)}(\tE\GG)$ is given by
\begin{align*}
  \delta(\iota^k) &= 0 \text{ for all $k$ and} \\
  \delta(e^m\xi^n) &=
   \begin{cases}
    e^m\delta\xi^n & \text{if $n\leq -1$} \\
    0 & \text{otherwise.}
   \end{cases}
\end{align*}
$\psi\colon \Mackey H_\GG^{RO(\GG)}(\tE\GG) \to \Mackey H_\GG^{RO(\GG)}(S^0)$ is given by
\begin{align*}
  \psi(e^m\kappa) &= e^m\kappa \\
  \psi(e^m\delta\xi^{-n}) &= 
     \begin{cases}
      e^m\delta\xi^{-n} & \text{if $m\leq -1$} \\
      0 & \text{otherwise.}
     \end{cases}
\end{align*}
$\phi\colon \Mackey H_\GG^{RO(\GG)}(S^0) \to \Mackey H_\GG^{RO(\GG)}((E\GG)_+)$ is given by
\begin{align*}
 \phi(\iota^k) &= \iota^k \\
 \phi(e^m\xi^n) &= e^m\xi^n \text{ for $m\geq 0$ and $n\geq 0$} \\
 \phi(e^{-m}\kappa) &= 0 \text{ for $m\geq 1$} \\
 \phi(e^{-m}\delta\xi^{-n}) &= 0.
\end{align*}
\qed
\end{proposition}

\section{The cohomology of a point with other coefficient systems}\label{sec:cohomPointOtherCoeffs}

As mentioned, the cohomology of a point with constant $\Z$ coefficients has been used
frequently in the recent literature. We review its calculation and point out
that several other calculations follow from it easily.

First, let $C$ be any abelian group and consider $\conc{C}$, the Mackey functor with
$\conc{C}(\GG/\GG) = C$ and $\conc{C}(\GG/e) = 0$.

\begin{proposition}\label{prop:concCcohomology}
There is a natural isomorphism
\[
	H_\GG^\alpha(X;\conc{C}) \iso H^{\alpha^\GG}(X^\GG;C),
\]
for any $\alpha\in RO(\GG)$.
\end{proposition}

\begin{proof}
This is a special case of \cite[1.13.22]{CW:ordinaryhomology}, or can be seen directly as follows:
$H^{\alpha^\GG}(X^\GG;C)$ is an $RO(\GG)$-graded cohomology theory in based $\GG$-spaces $X$.
It obeys a dimension axiom in integer grading, with
$H^n((\GG/\GG)^\GG_+;C) \iso C$ if $n=0$ but equal to $0$ if $n\neq 0$, and
$H^n((\GG/e)^\GG_+;C) = 0$ for all $n$.
This is precisely the dimension axiom satisfied by $H_\GG^{RO(\GG)}(X;\conc{C})$,
so the two theories must be naturally isomorphic by the uniqueness of equivariant
ordinary cohomology.
\end{proof}

In particular, consider $\conc\Z$ and the short exact sequence
\[
 0 \to \conc\Z \xrightarrow{\kappa} \Mackey A \xrightarrow{\epsilon} \Mackey\Z \to 0,
\]
where the map $\kappa$ takes $1$ to $\kappa\in A(\GG)$.
We use this to think of $\conc\Z$ as the submodule of $\Mackey A$
consisting of the multiples of $\kappa$.

\begin{proposition}\label{prop:concZcohomologypoint}
Additively,
\[
 \Mackey H_\GG^\alpha(S^0;\conc\Z) \iso
  \begin{cases}
    \conc\Z & \text{if $\alpha^\GG = 0$} \\
    0 & \text{otherwise.}
  \end{cases}
\]
The map $\Mackey H_\GG^{RO(\GG)}(S^0;\conc\Z) \to \Mackey H_\GG^{RO(\GG)}(S^0;\Mackey A)$
is injective with image the
ideal
\[
  \langle e^n\kappa \mid n\in\Z \rangle 
  \subset \Mackey H_\GG^{RO(\GG)}(S^0;\Mackey A),
\]
so we identify $\Mackey H_\GG^{RO(\GG)}(S^0;\conc\Z)$ with this ideal.
(Recall that $e^n\kappa = 2e^n$ if $n\geq 1$.)
The inclusion factors through
$\Mackey H_\GG^{RO(\GG)}(\tE\GG;\Mackey A)$, with image there
the submodule of elements in gradings
$\alpha$ with $\alpha^\GG = 0$.
\end{proposition}

\begin{proof}
The additive calculation is immediate from Proposition~\ref{prop:concCcohomology}.
The dimension axiom tells us that  $\Mackey H_\GG^0(S^0;\conc\Z)\to \Mackey H_\GG^0(S^0;\Mackey A)$
is the inclusion $\conc\Z\to \Mackey A$.
Write $\kappa\in \Mackey H_\GG^0(S^0;\conc\Z)$ for the generator that maps
to $\kappa\in \Mackey H_\GG^0(S^0;\Mackey A)$.

Consider the inclusion $S^0\to S^{\sigma}$. Again by Proposition~\ref{prop:concCcohomology},
we have that the induced map
$\Mackey H_\GG^{RO(\GG)}(S^{\sigma};\conc\Z)\to \Mackey H_\GG^{RO(\GG)}(S^0;\conc\Z)$ is an
isomorphism. But this map is multiplication by $e$, so multiplication by $e$
is an isomorphism on $\Mackey H_\GG^{RO(\GG)}(S^0;\conc\Z)$. This implies that, for $n\in\Z$,
$\Mackey H_\GG^{n\sigma}(S^0;\conc\Z)$
is generated by $e^n\kappa$, an element that maps to
$e^n\kappa\in H_\GG^\alpha(S^0;\Mackey A)$.

We can now see that
the map $\Mackey H_\GG^{RO(\GG)}(S^0;\conc\Z) \to \Mackey H_\GG^{RO(\GG)}(S^0;\Mackey A)$
is injective with image the ideal $\langle e^n\kappa \mid n\in\Z \rangle$.

For the last statement, Proposition~\ref{prop:concCcohomology} implies
that
\[
 \Mackey H_\GG^{RO(\GG)}(S^0;\conc\Z) \iso \Mackey H_\GG^{RO(\GG)}(\tE\GG;\conc\Z)
\]
because $\tE\GG^\GG \hmtpc S^0$. Thus, we get a factorization
\begin{multline*}
 \Mackey H_\GG^{RO(\GG)}(S^0;\conc\Z) \iso \Mackey H_\GG^{RO(\GG)}(\tE\GG;\conc\Z) \\
  \to \Mackey H_\GG^{RO(\GG)}(\tE\GG;\Mackey A)
  \to \Mackey H_\GG^{RO(\GG)}(S^0;\Mackey A).
\end{multline*}
The image in $\Mackey H_\GG^{RO(\GG)}(\tE\GG;\Mackey A)$ is clear from our computations.
\end{proof}

It's useful to notice that the image of
$\Mackey H_\GG^{RO(\GG)}(S^0;\conc\Z) \to \Mackey H_\GG^{RO(\GG)}(\tE\GG;\Mackey A)$
is a direct summand. The other summand is given by the groups in gradings $\alpha$
with $\alpha^\GG\neq 0$; it's straightforward to check that this is a submodule.

\begin{theorem}\label{thm:pointEvenRZcoeffs}
Additively,
\[
 \Mackey H_\GG^\alpha(S^0;\Mackey\Z) \iso
  \begin{cases}
    \Mackey\Z & \text{if $|\alpha| = 0$ and $\alpha^\GG \leq 0$ is even} \\
    \Mackey\Z_- & \text{if $|\alpha| = 0$ and $\alpha^\GG \leq 1$ is odd} \\
    \Mackey\Z' & \text{if $|\alpha| = 0$ and $\alpha^\GG > 0$ is even} \\
    \Mackey\Z'_- & \text{if $|\alpha| = 0$ and $\alpha^\GG \geq 3$ is odd} \\
    \conc{\Z/2} & \text{if $|\alpha| > 0$ and $\alpha^\GG \leq 0$ is even} \\
    \conc{\Z/2} & \text{if $|\alpha| < 0$ and $\alpha^\GG \geq 3$ is odd} \\
    0 & \text{otherwise}.
  \end{cases}
\]
$\Mackey H_\GG^{RO(\GG)}(S^0;\Mackey\Z)$ is a strictly commutative $RO(\GG)$-graded 
algebra over $\Mackey\Z$,
generated multiplicatively by elements
\begin{align*}
 \iota &\in \Mackey H_\GG^{\sigma-1}(S^0;\Mackey\Z)(\GG/e) \\
 \iota^{-1} &\in \Mackey H_\GG^{1-\sigma}(S^0;\Mackey\Z)(\GG/e) \\
 \xi &\in \Mackey H_\GG^{2(\sigma-1)}(S^0;\Mackey\Z)(\GG/\GG) \\
 e &\in \Mackey H_\GG^{\sigma}(S^0;\Mackey\Z)(\GG/\GG) \\
 e^{-m}\delta\xi^{-n} &\in \Mackey H_\GG^{1 - m\sigma - 2n(\sigma-1)}(S^0;\Mackey\Z)(\GG/\GG)
   & & m, n \geq 1.
\end{align*}
These generators satisfy the following {\em structural} relations:
\begin{align*}
 \tau(\iota^{-1}) &= 0 \\
 \rho(\xi) &= \iota^2 \\
 \rho(e) &= 0 \\
 e^{-1}\delta\xi^{-n} &= \tau(\iota^{-2n-1}) & & \text{for $n\geq 1$} \\
 \rho(e^{-m}\delta\xi^{-n}) &= 0 & & \text{for $m\geq 2$ and $n\geq 1$}\\
\intertext{and the following {\em multiplicative} relations:}
  \iota\cdot \iota^{-1} &= \rho(1) \\
 e\cdot e^{-m}\delta\xi^{-n} &= e^{-m+1}\delta\xi^{-n} & &\text{for $m\geq 2$ and $n\geq 1$} \\
 \xi\cdot e^{-m}\delta\xi^{-n} &= e^{-m}\delta\xi^{-n+1} & & 
   \text{for $m\geq 1$ and $n\geq 2$} \\
 \xi\cdot e^{-m}\delta\xi^{-1} &= 0 & & \text{for $m\geq 2$} \\
 e^{-m}\delta\xi^{-k}\cdot e^{-n}\delta\xi^{-\ell} &= 0 & & 
   \text{if $m, n, k, \ell\geq 1$}
\end{align*}
The following relations are implied by the preceding ones:
\begin{align*}
 2e^m\xi^n &= 0 & & \text{if $m > 0$ and $n \geq 0$} \\
 2e^{-m}\delta\xi^{-n} &= 0 & & \text{if $m\geq 2$ and $n\geq 1$} \\
 t\iota^k &= (-1)^k\iota^k & & \text{for all $k$} \\
 \xi\cdot\tau(\iota^k) &= \tau(\iota^{k+2}) & &\text{for all $k$} \\
 e\cdot\tau(\iota^k) &= 0 & &\text{for all $k$} \\
 e^{-m}\delta\xi^{-n}\cdot \tau(\iota^k) &= 0 & & \text{for all $m, n\geq 1$ and $k$} \\
 \tau(\iota^k)\cdot\tau(\iota^\ell) &= 0 & &\text{if $k$ or $\ell$ is odd} \\
 \tau(\iota^{2k})\cdot\tau(\iota^{2\ell}) &= 2\tau(\iota^{2(k+\ell)}) & &\text{for all $k$ and $\ell$}\\
 \tau(\iota^{2k+1}) &= 0 & & \text{if $k\geq 0$} \\
 e\cdot e^{-1}\delta\xi^{-n} &= 0 & & \text{if $n\geq 1$}
\end{align*}
\end{theorem}

\begin{proof}
By Proposition~\ref{prop:concZcohomologypoint}, we have a short exact sequence
\[
 0 \to \Mackey H_\GG^{RO(\GG)}(S^0;\conc\Z) \to \Mackey H_\GG^{RO(\GG)}(S^0;\Mackey A)
  \to \Mackey H_\GG^{RO(\GG)}(S^0;\Mackey\Z) \to 0
\]
exhibiting $\Mackey H_\GG^{RO(\GG)}(S^0;\Mackey\Z)$ as a quotient ring of
$\Mackey H_\GG^{RO(\GG)}(S^0;\Mackey A)$.
The additive calculation follows by noticing that the elements killed in the quotient
are all of the $e^{-m}\kappa$ for $m\geq 0$ and the elements $2e^m$ for $m\geq 1$.
The rest is just seeing what relations are still necessary from those in
Theorem~\ref{thm:evenCohomPoint}.
We also simplify a bit by noticing that, for a $\Mackey\Z$-module
generated by an element $x$ at level $\GG/\GG$, $\rho(x) = 0$ implies that $2x = 0$.
\end{proof}

\begin{figure}
\[\def\objectstyle{\scriptstyle}
 \xymatrix@!0@R=5ex@C=2.5em{
  & & & & &\phantom{a}\ b\sigma & & & & & \\
  & \conc{\Z/2} &\cdot& \conc{\Z/2} &\cdot& \conc\Z &\cdot& \cdot &\cdot& \cdot & \cdot \\
  & \Mackey\Z &\cdot& \conc{\Z/2} &\cdot& \conc\Z &\cdot& \cdot &\cdot& \cdot & \cdot \\
  & \cdot &\Mackey\Z_{\mathrlap{-}}& \conc{\Z/2} &\cdot& \conc\Z &\cdot& \cdot &\cdot& \cdot & \cdot \\
  & \cdot &\cdot& \Mackey\Z &\cdot& \conc\Z &\cdot& \cdot &\cdot& \cdot & \cdot \\
  & \cdot &\cdot& \cdot &\Z_{\mathrlap{-}}& \conc\Z &\cdot& \cdot &\cdot& \cdot & \cdot \\
  \ar@{-}'[rrrrr][rrrrrrrrrrr]
   &  &  &  &  & \Mackey \Z &  &  &  &  &  & a\\
  & \cdot &\cdot& \cdot &\cdot&  &\Mackey \Z_{\mathrlap{-}}& \cdot & \cdot & \cdot & \cdot  \\
  & \cdot &\cdot& \cdot &\cdot&  &\cdot& \Mackey \Z' & \cdot & \cdot & \cdot  \\
  & \cdot &\cdot& \cdot &\cdot&  &\cdot& \cdot & \Mackey \Z'_{\mathrlap{-}} & \cdot & \cdot  \\
  & \cdot &\cdot& \cdot &\cdot&  &\cdot& \cdot & \conc{\Z/2} & \Mackey \Z' &  \cdot  \\
  & \cdot &\cdot& \cdot &\cdot&  &\cdot& \cdot & \conc{\Z/2} & \cdot &  \Mackey \Z'_{\mathrlap{-}}  \\
  & \cdot &\cdot& \cdot &\cdot&  &\cdot& \cdot & \conc{\Z/2} & \cdot &  \conc{\Z/2}  \\
  & & & & & \ar@{-}'[uuuuuuu]'[uuuuuuuu]'[uuuuuuuuu]'[uuuuuuuuuu]'[uuuuuuuuuuu]'[uuuuuuuuuuuu][uuuuuuuuuuuuu]
 }
\]
\caption{$\protect\Mackey H_\GG^{RO(\GG)}(S^0;\protect\Mackey\Z)$}\label{fig:EvenCohomPointRZ}
\end{figure}
\begin{figure}
\[\def\objectstyle{\scriptstyle}
 \xymatrix@!0@R=5ex@C=3em{
  & & & & & \phantom{a}\  b\sigma & & & & & \\
  & e\xi^2 &\cdot& e^3\xi &\cdot& e^5 &\cdot& \cdot &\cdot& \cdot & \cdot \\
  & \xi^2 &\cdot& e^2\xi &\cdot& e^4 &\cdot& \cdot &\cdot& \cdot & \cdot \\
  & \cdot &(\iota^3) & e\xi &\cdot& e^3 &\cdot& \cdot &\cdot& \cdot & \cdot \\
  & \cdot &\cdot& \xi &\cdot& e^2 &\cdot& \cdot &\cdot& \cdot & \cdot \\
  & \cdot &\cdot& \cdot &(\iota)& e &\cdot& \cdot &\cdot& \cdot & \cdot \\
  \ar@{-}'[rrrrr][rrrrrrrrrrr]
   &  &  &  &  & 1 &  &  &  &  &  & a\\
  & \cdot &\cdot& \cdot &\cdot&  &(\iota^{-1})& \cdot & \cdot & \cdot &  \cdot  \\
  & \cdot &\cdot& \cdot &\cdot&  &\cdot& (\iota^{-2}) & \cdot & \cdot & \cdot  \\
  & \cdot &\cdot& \cdot &\cdot&  &\cdot& \cdot & (\iota^{-3}) & \cdot &  \cdot  \\
  & \cdot &\cdot& \cdot &\cdot&  &\cdot& \cdot & e^{-2}\delta\xi^{-1} & (\iota^{-4}) &  \cdot  \\
  & \cdot &\cdot& \cdot &\cdot&  &\cdot& \cdot & e^{-3}\delta\xi^{-1} & \cdot &  (\iota^{-5})  \\
  & \cdot &\cdot& \cdot &\cdot&  &\cdot& \cdot & e^{-4}\delta\xi^{-1} & \cdot &  e^{-2}\delta\xi^{-2}  \\
  & & & & & \ar@{-}'[uuuuuuu]'[uuuuuuuu]'[uuuuuuuuu]'[uuuuuuuuuu]'[uuuuuuuuuuu]'[uuuuuuuuuuuu]'[uuuuuuuuuuuuu]
 }
\]
\caption{The generators of $\protect\Mackey H_\GG^{RO(\GG)}(S^0;\protect\Mackey\Z)$}\label{fig:EvenCohomPointRZGens}
\end{figure}
Figures~\ref{fig:EvenCohomPointRZ} and~\ref{fig:EvenCohomPointRZGens} 
show $\Mackey H_\GG^{RO(\GG)}(S^0;\Mackey\Z)$ and its
generators.
You can see something interesting in these figures.
Recall that the integer-graded part lies along the horizontal axis. 
As it must, because of the dimension axiom,
it contains only one nonzero group, which is $\Mackey\Z$.
Now look at the horizontal line just below the axis. It also contains only one nonzero group,
which is $\Mackey\Z_-$. So, if we were to shift the $\Mackey\Z$ cohomology up and to the left by one,
uniqueness shows that we get cohomology with coefficients in $\Mackey\Z_-$.
This continues for the two horizontal lines below that and gives us the following result.

\begin{corollary}\label{cor:manyIsoCoeffs}
There are natural isomorphisms
\begingroup\allowdisplaybreaks[0]
\begin{align*}
 \Mackey H_\GG^\alpha(X;\Mackey\Z)
  &\iso \Mackey H_\GG^{\alpha+(\sigma-1)}(X;\Mackey\Z_-) \\
  &\iso \Mackey H_\GG^{\alpha+2(\sigma-1)}(X;\Mackey\Z') \\
  &\iso \Mackey H_\GG^{\alpha+3(\sigma-1)}(X;\Mackey\Z'_-).
\end{align*}
\endgroup
Moreover, these are isomorphisms of modules over
$\Mackey H_\GG^{RO(\GG)}(S^0;\Mackey\Z)$.
\end{corollary}

\begin{proof}
The proof of the isomorphisms was given above. That all these theories are modules
over $\Mackey H_\GG^{RO(\GG)}(S^0;\Mackey\Z)$ follows from the fact that
$\Mackey\Z_-$, $\Mackey\Z'$, and $\Mackey\Z'_-$ are all modules over $\Mackey\Z$,
as noted in \S\ref{sec:MackeyFunctors}, or as exhibited in the fact that all
appear as groups within $\Mackey H_\GG^{RO(\GG)}(S^0;\Mackey\Z)$.
\end{proof}

\part{The cohomology of $B_\GG U(1)$}\label{part:BU1}

\section{The topology of $B_\GG U(1)$}

We introduce an explicit model for the classifying space $B_\GG U(1)$ and discuss its topology
and its fundamental groupoid.

Recall from Definition~\ref{def:irrRepresentations}
that we write $\C$ for the trivial complex representation of $G$ and
$\C^\sigma$ for the nontrivial irreducible complex representation.
As a model for $B_\GG U(1)$ we take
\[
 B_\GG U(1) = \CP_\GG^\infty = \CP(\C^\infty\dirsum(\C^\sigma)^\infty),
\]
the projective space of complex lines in $\C^\infty\dirsum(\C^\sigma)^\infty$.
Nonequivariantly, this is $\CP^\infty$.
Its fixed sets are
\[
 B_\GG U(1)^\GG = \CP(\C^\infty) \disjunion \CP((\C^\sigma)^\infty),
\]
the disjoint union of two copies of $\CP^\infty$.
For notational simplicity we shall write $B = B_\GG U(1)$, $B^0 = \CP(\C^\infty)$,
and $B^1 = \CP((\C^\sigma)^\infty)$,
so $B^\GG = B^0 \disjunion B^1$.

Let $\omega$ denote the tautological complex line bundle over $B$.
Nonequivariantly, it is the usual tautological line bundle over the infinite complex projective space.
Its restriction to $B^0$ is the nonequivariant tautological bundle with $\GG$ acting trivially on the fibers;
its restriction to $B^1$ is the nonequivariant tautological bundle with $\GG$
acting on each fiber by negation, as it does on $\C^\sigma$.

There is an equivariant involution of $B$ we will want to take into account in our calculations.
Write elements of $B = \CP(\C^\infty\dirsum(\C^\sigma)^\infty)$ as equivalence classes of pairs
$[z_0 : z_1]$, with $z_0\in \C^\infty$ and $z_1\in (\C^\sigma)^\infty$ not both 0. 
We define a $G$-map $\chi\colon B\to B$ by
\[
 \chi[z_0 : z_1] = [z_1 : z_0],
\]
using our chosen nonequivariant identification of $\C$ with $\C^\sigma$.
It is straightforward to check that $\chi$ is equivariant, and clearly $\chi^2$ is the identity.
Further, $\chi^\GG$ swaps $B^0$ and $B^1$ via a homeomorphism between them.
Notice also that $\chi^*\omega \iso \omega\tensor_\C \C^\sigma$.
More generally, if $f\colon X\to B$ classifies the complex line bundle $\theta$ over $X$,
then $\chi f$ classifies $\theta\tensor_\C \C^\sigma$.

The fundamental groupoid $\pi\colon \Pi B\to \orb \GG$
is relatively simple because $B$ and 
the components of its fixed set are all simply connected.
It is equivalent to a category over $\orb \GG$ having one object
$b$ over $\GG/e$ and two objects, $b_0$ and $b_1$, over $\GG/\GG$,
corresponding to the two components of $B^\GG$.
The self-maps of $b$ map isomorphically to the self-maps of $\GG/e$;
there is one map $b\to b_k$ for each $k$, over the map $\GG/e\to \GG/\GG$;
and there are no non-identity self-maps of $b_k$.
We can picture the category as follows:
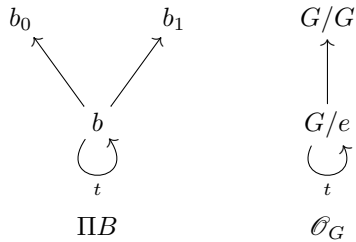
\[
 \xymatrix@!C=.5em{
  b_0 && b_1 & \qquad\qquad & \GG/\GG \\
  &b \ar[ul]^{\rho_0} \ar[ur]_{\rho_1} \ar@(dl,dr)[]_t &&& \GG/e \ar[u]_\rho \ar@(dl,dr)[]_t \\
  &\Pi B &&& \orb{\GG}
 }
\]
There is also an action of $\chi$, fixing $b$ and exchanging the objects $b_0$ and $b_1$.
Note that $\rho_k t = \rho_k$ for $k = 0$ or $1$.

We now want to compute the representation ring $RO(\Pi B)$.
Let $\alpha$ be a real virtual representation
of $\Pi B$.
We have $\alpha(b) = (\GG\times \R^i)\ominus (\GG\times \R^j)$ and, by an abuse of notation,
we shall write this as $\alpha(b) = \GG\times\R^{i-j} = \GG\times\R^n$ where
$n\in\Z$. (For this discussion, we care only about the isomorphism class of $\alpha$.)
The map $t$ acts on this bundle by its nontrivial action on $\GG$ and the
homotopy class of a linear map on $\R^n$ whose square is homotopically trivial.
There are, therefore, two possible
actions of $t$ on $\alpha(b)$: the map $t_+$ in which $t$ acts (homotopically) trivially on $\R^n$
and the map $t_-$ in which $t$ acts on $\R^n$ by any orientation-reversing linear map.
If $\alpha(b_0) = \R^{n_0} \dirsum (\R^\sigma)^{n_1}$ with $n_0 + n_1 = n$, when is there a $\GG$-map
$\alpha(\rho_0)\colon \GG\times\R^n\to \R^{n_0} \dirsum (\R^\sigma)^{n_1}$ such that $\alpha(\rho_0)\alpha(t) = \alpha(\rho_0)$,
so that $\alpha$ is a functor?
If the action of $t$ on $\GG\times\R^n$ is by $t_+$, then we must have $n_1$ even, whereas, if
the action of $t$ is by $t_-$, we must have $n_1$ odd. The same applies to $\alpha(b_1)$.
Thus, $\alpha$ is entirely determined by its values $\alpha(b_0)$ and $\alpha(b_1)$,
with the restriction that the parity of their nontrivial parts must agree.
This gives us the following.

\begin{proposition}\label{prop:repringB}
\[
 RO(\Pi B)
  \iso \{ (\alpha_0, \alpha_1) \mid \alpha_k\in RO(\GG), |\alpha_0| = |\alpha_1|,
      \text{ and }\alpha_0^\GG \equiv \alpha_1^\GG \pmod 2 \}
\]
is a free abelian group of rank 3. It has as a basis the elements
\begin{align*}
 1 &= (1, 1) \\
 \sigma &= (\sigma, \sigma) \\
 \mathllap{\text{and}\qquad}\Omega &= (1-\sigma, \sigma-1).
\end{align*}
\end{proposition}

\begin{proof}
The argument that $RO(\Pi B)$ is the indicated subgroup of $RO(\GG)^2$
was given just before the statement of the proposition,
where we are now writing $\alpha_0 = \alpha(b_0)$ and $\alpha_1 = \alpha(b_1)$.
Note that, if $|\alpha_0| = |\alpha_1|$, then $\alpha_0^\GG \equiv \alpha_1^\GG \pmod 2$
if and only if $|\alpha_0 - \alpha_0^\GG| \equiv |\alpha_1 - \alpha_1^\GG| \pmod 2$.

Now, for any $(\alpha_0,\alpha_1) \in RO(\Pi B)$, let
$n = (\alpha_0^\GG-\alpha_1^\GG)/2$. Then
\[
 (\alpha_0, \alpha_1) = 
  (\alpha_0^\GG - n)\cdot 1 + (|\alpha_0|-\alpha_0^\GG + n)\cdot\sigma
    + n\cdot\Omega.
\]
On the other hand, if $a\cdot 1 + b\cdot\sigma + c\cdot\Omega = 0$, 
then, looking at components, it is easy to see that we must have
$a = b = c = 0$. Hence, $\{1, \sigma, \Omega\}$ is a basis.
\end{proof}

It will be useful to introduce the following elements:
\begin{align*}
 \Omega_{0} &= (2\sigma - 2, 0) = -1 + \sigma - \Omega \\
 \Omega_{1} &= (0, 2\sigma - 2) = -1 + \sigma + \Omega.
\end{align*}
We can then write
\[
 RO(\Pi B) \iso \Z\{1,\sigma,\Omega_0,\Omega_1\}/\langle \Omega_0 + \Omega_1 = 2\sigma - 2\rangle
\]
where $\Z\{1,\sigma,\Omega_0,\Omega_1\}$ indicates the free abelian group
on the generators $1$, $\sigma$, $\Omega_0$ and $\Omega_1$.
This is the presentation of $RO(\Pi B)$ we will use most often.

The involution $\chi$ acts on $RO(\Pi B)$ by 
$\chi(\alpha_0,\alpha_1) = (\alpha_1,\alpha_0)$, so we have
\begin{align*}
 \chi(1) &= 1 & \chi(\sigma) &= \sigma \\
 \chi(\Omega_{0}) &= \Omega_{1} & \chi(\Omega_{1}) &= \Omega_{0} \\
 \chi(\Omega) &= -\Omega.
\end{align*}

The tautological line bundle $\omega$ induces a representation 
$\omega^* \in RO(\Pi B)$, with
\[
 \omega^* = (2, 2\sigma) = 2 + \Omega_{1}.
\]
Similarly,
\[
 (\chi\omega)^* = (2\sigma, 2) = 2 + \Omega_0.
\]
For simplicity of notation, we will often write $\omega$ for $\omega^*$
and $\chi\omega$ for $(\chi\omega)^*$ when the context
makes clear that we are speaking of the associated representation of the bundle, not the bundle itself.

\begin{lemma}\label{lem:kernelproj}
The inclusion $B^0 \to B$ induces the map $RO(\Pi B)\to RO(\Pi B^0) \iso RO(\GG)$ given by
$(\alpha_0, \alpha_1) \mapsto \alpha_0$. Its kernel is the free abelian subgroup
generated by $\Omega_1$.

Similarly, the inclusion $B^1 \to B$ induces the map $RO(\Pi B)\to RO(\GG)$
given by $(\alpha_0, \alpha_1) \mapsto \alpha_1$, and the kernel of
this map is generated by $\Omega_0$.
\end{lemma}

\begin{proof}
That the inclusions $B^k \to B$ induce the maps claimed follows from the way we identified
elements of $RO(\Pi B)$ with pairs $(\alpha_0, \alpha_1)$.

To identify the kernels, let $k = 0$ or $1$.
If $\alpha\in RO(\Pi B)$ and $\alpha_k = 0$, then $|\alpha_{1-k}| = 0$
because the dimensions are equal, and $\alpha_{1-k}^\GG$ is even
because the fixed sets have the same parity.
The set of elements $\beta \in RO(\GG)$ with $|\beta| = 0$ and $\beta^\GG$ even
is the subgroup generated by  $2\sigma-2$.
This proves the lemma.
\end{proof}

Our calculation of the cohomology of $B$ will use a combination of separation of isotropy and
restrictions to fixed points. 
These give us the following diagram,
in which we write $R = RO(\Pi B)$ for brevity,
whose rows and columns are parts of long exact sequences.
\begin{equation}\label{diag:main}\def\objectstyle{\scriptstyle}
 \xymatrix@C-0.5em{
   \Mackey H_\GG^R(B_+\smsh_B\tE\GG) \ar[r]^\eta \ar[d]_{\psi}
     & \Mackey H_\GG^R(B_+^\GG\smsh_B\tE\GG) \ar[d]^\psi \\
   \Mackey H_\GG^R(B_+) \ar[r]^\eta \ar[d]_\phi
     & \Mackey H_\GG^R(B_+^\GG) \ar[r]^-\theta \ar[d]^\phi
     & \susp^{-1}\Mackey H_\GG^{R}(B/_B B^\GG) \ar[d]^\phi \\
   \Mackey H_\GG^R(B_+\smsh_B (E\GG)_+) \ar[r]^\eta \ar[d]_\delta
     & \Mackey H_\GG^R(B_+^\GG\smsh_B (E\GG)_+) \ar[r]^-\theta \ar[d]^\delta
     & \susp^{-1}\Mackey H_\GG^{R}(B/_B B^\GG\smsh_B (E\GG)_+) \\
   \susp^{-1}\Mackey H_\GG^{R}(B_+\smsh_B\tE\GG) \ar[r]^\eta
     & \susp^{-1}\Mackey H_\GG^{R}(B_+^\GG\smsh_B\tE\GG)
  }
\end{equation}
(The missing corners are 0, as we shall see.)
With that in mind, in the next two sections we calculate the cohomologies that appear
in this diagram other than the cohomology of $B$ itself.

\section{The cohomology of $B_\GG U(1)^\GG$}

In this section we calculate the $RO(\Pi B)$-graded cohomologies of the spaces
$B_+^\GG$, $B_+^\GG\smsh_{B} (E\GG)_+$, and $B^\GG_+\smsh_{B}\tE\GG$, as ex-spaces over $B$.
Recall that $B^\GG = B^0 \disjunion B^1$, where each $B^k$ is a copy of $\CP^\infty$.
So, we begin by considering the equivariant cohomology of this nonequivariant space,
for which we need some results that hold for all groups $G$.

\begin{proposition}\label{prop:cohomologyTrivialAction}
Let $X$ be a based space with trivial $G$-action and let $\Mackey T$ be a Mackey functor.
Then, in integer grading,
\[
 \Mackey H_G^n(X;\Mackey T)(G/K) \iso H^n(X;\Mackey T(G/K)),
\]
naturally in $X$ and $K$.
\end{proposition}

\begin{proof}
Both sides may be considered as integer-graded nonequivariant cohomology theories in $X$.
They obey the same dimension axiom, so, by the uniqueness of ordinary cohomology,
they are naturally isomorphic. Naturality in $K$ follows similarly. 
\end{proof}

\begin{proposition}
Let $X$ be a based space with trivial $G$-action, let $A$ be a $G$-space,
and let $Y\to A$ be an ex-space over $A$.
Let $R$ be a commutative ring and suppose that the nonequivariant cohomology
groups $H^n(X;R)$ are flat $R$-modules for all $n$.
Further, suppose that $\Mackey T$ is a Mackey $R$-module, that is, a functor
from $\orb G$ to $R$-modules.
Then
\[
 \Mackey H_G^{RO(\Pi A)} (X\smsh Y; \Mackey T)
 	\iso H^{\Z}(X;R) \tensor_R \Mackey H_G^{RO(\Pi A)}(Y; \Mackey T).
\]
\end{proposition}

\begin{proof}
Extending the usual notation, by the tensor product on the right we mean the $RO(\Pi A)$-graded 
Mackey functor defined by
\[
 \left(H^{\Z}(X;R) \tensor_R \Mackey H_G^{RO(\Pi A)}(Y; \Mackey T)\right)^\alpha
  = \Dirsum_{n\in\Z} H^n(X;R) \tensor_R \Mackey H_G^{\alpha-n}(Y;\Mackey T).
\]

Now, both sides of the claimed isomorphism
are cohomology theories in ex-spaces $Y$ over $A$, the right side because of 
the flatness assumption.
There is a pairing carrying the right side to the left,
which is an isomorphism in integer grading for $Y = G/K_+\to A$,
using Proposition~\ref{prop:cohomologyTrivialAction}.
There are then Atiyah-Hirzebruch spectral sequences from \cite[3.6.1]{CW:ordinaryhomology}
of the form
\[
 \Mackey H_G^{\alpha}(Y;\Mackey H_G^q(X;\Mackey T)) \convto
   \Mackey H_G^{\alpha+q}(X\smsh Y;\Mackey T)
\]
and
\begin{multline*}
 \Mackey H_G^{\alpha}\left(Y;\left(H^{\Z}(X;R) \tensor_R \Mackey H_G^{RO(\Pi A)}(S^0; \Mackey T)\right)^q\right) \\
   \convto \left(H^{\Z}(X;R) \tensor_R \Mackey H_G^{RO(\Pi A)}(Y; \Mackey T)\right)^{\alpha+q}
\end{multline*}
whose $E_2$ pages are isomorphic via the pairing,
showing that the two cohomology theories are isomorphic for all $Y$.
\end{proof}

We now return to our specific case of $G = \GG$.

\begin{proposition}\label{prop:integerCPinfty}
Let $\CP^\infty$ be the infinite complex projective space considered as a $\GG$-space
with trivial $G$-action. We have
\begin{align*}
 \Mackey H_\GG^{RO(\GG)}(\CP^\infty_+) &\iso \Mackey H_\GG^{RO(\GG)}(S^0)[c], \\
 \Mackey H_\GG^{RO(\GG)}(\CP^\infty_+\smsh (E\GG)_+) &\iso 
        \Mackey H_\GG^{RO(\GG)}(S^0)[c]\tensor_{\Mackey H_\GG^{RO(\GG)}(S^0)} \Mackey H_\GG^{RO(\GG)}((E\GG)_+) \\
        &\iso \Mackey H_\GG^{RO(\GG)}((E\GG)_+)[c],  \\
\intertext{and}
 \Mackey H_\GG^{RO(\GG)}(\CP^\infty_+\smsh\tE\GG) &\iso 
        \Mackey H_\GG^{RO(\GG)}(S^0)[c]\tensor_{\Mackey H_\GG^{RO(\GG)}(S^0)} \Mackey H_\GG^{RO(\GG)}(\tE\GG),
\end{align*}
where $c$ is the Euler class of the tautological complex line bundle $\omega$
over $\CP^\infty$ with trivial $\GG$-action, so $|c| = 2$. 
Moreover, the long exact sequence coming from the cofibration sequence
\[
 \CP^\infty_+\smsh (E\GG)_+ \to \CP^\infty_+ \to \CP^\infty_+\smsh\tE\GG
\]
is given by tensoring $\Mackey H_\GG^{RO(\GG)}(S^0)[c]$ with the long exact sequence
coming from the cofibration sequence
$(E\GG)_+\to S^0 \to \tE\GG$.
\end{proposition}

\begin{proof}
This follows from the preceding proposition, taking $X = \CP^\infty_+$ and $Y = S^0$, $(E\GG)_+$, or
$\tE\GG$,
using the nonequivariant calculation
$H^\Z(\CP^\infty_+;\Z) \iso \Z[c]$, where $c$ is the nonequivariant Euler class
of the tautological line bundle.
\end{proof}

Now consider $B^0$ and $B^1$ as spaces over $B$. The following result
calculates their $RO(\Pi B)$-graded cohomologies.

\begin{proposition}\label{prop:calcFixedPoints}
\begin{align*}
 \Mackey H_\GG^{RO(\Pi B)}(B^0_+) 
   &\iso \Mackey H_\GG^{RO(\GG)}(S^0)[c, \zeta_1, \zeta_1^{-1}] \\
 \Mackey H_\GG^{RO(\Pi B)}(B^0_+\smsh_B (E\GG)_+) 
   &\iso \Mackey H_\GG^{RO(\Pi B)}(B^0_+)\tensor_{\Mackey H_\GG^{RO(\GG)}(S^0)} \Mackey H_\GG^{RO(\GG)}((E\GG)_+) \\
 \Mackey H_\GG^{RO(\Pi B)}(B^0_+\smsh_B \tE\GG) 
   &\iso \Mackey H_\GG^{RO(\Pi B)}(B^0_+)\tensor_{\Mackey H_\GG^{RO(\GG)}(S^0)} \Mackey H_\GG^{RO(\GG)}(\tE\GG),
\end{align*}
where $|c| = 2$ and $|\zeta_1| = \Omega_{1}$.
Similarly,
\begin{align*}
 \Mackey H_\GG^{RO(\Pi B)}(B^1_+) 
   &\iso \Mackey H_\GG^{RO(\GG)}(S^0)[c, \zeta_0, \zeta_0^{-1}] \\
 \Mackey H_\GG^{RO(\Pi B)}(B^1_+\smsh_B (E\GG)_+) 
   &\iso \Mackey H_\GG^{RO(\Pi B)}(B^1_+)\tensor_{\Mackey H_\GG^{RO(\GG)}(S^0)} \Mackey H_\GG^{RO(\GG)}((E\GG)_+) \\
 \Mackey H_\GG^{RO(\Pi B)}(B^1_+\smsh_B \tE\GG) 
   &\iso \Mackey H_\GG^{RO(\Pi B)}(B^1_+)\tensor_{\Mackey H_\GG^{RO(\GG)}(S^0)} \Mackey H_\GG^{RO(\GG)}(\tE\GG),
\end{align*}
where $|c| = 2$ and $|\zeta_0| = \Omega_{0}$.
This gives the following calculations of the cohomology of $B^\GG = B^0\disjunion B^1$:
\begin{align*}
 \Mackey H_\GG^{RO(\Pi B)}(B^\GG_+)
   &\iso \Mackey H_\GG^{RO(\Pi B)}(B^0_+) \dirsum \Mackey H_\GG^{RO(\Pi B)}(B^1_+)  \\
 \Mackey H_\GG^{RO(\Pi B)}(B^\GG_+\smsh_B (E\GG)_+) 
   &\iso \Mackey H_\GG^{RO(\Pi B)}(B^\GG_+)\tensor_{\Mackey H_\GG^{RO(\GG)}(S^0)} \Mackey H_\GG^{RO(\GG)}((E\GG)_+) \\
 \Mackey H_\GG^{RO(\Pi B)}(B^\GG_+\smsh_B \tE\GG) 
   &\iso \Mackey H_\GG^{RO(\Pi B)}(B^\GG_+)\tensor_{\Mackey H_\GG^{RO(\GG)}(S^0)} \Mackey H_\GG^{RO(\GG)}(\tE\GG).
\end{align*}
\end{proposition}

\begin{proof}
Consider the case of $B^0$.
Because $B^0$ is simply connected and has trivial $\GG$-action,
$RO(\Pi B^0) \iso RO(\GG)$.
By Lemma~\ref{lem:kernelproj}, the kernel of the induced map
$RO(\Pi B) \to RO(\Pi B^0)$ is the free abelian subgroup generated
by the element $\Omega_{1}$.
Let $\zeta_1\in \Mackey H_\GG^{\Omega_1}(B^0_+)$ be the unit given
by Proposition~\ref{prop:kerunits}.
Because $RO(\Pi B)$ is generated by $RO(G)$ and $\Omega_1$, we see that
the $RO(\Pi B)$-graded cohomology of $B^0$ is completely determined by its $RO(G)$-graded part together
with the invertible element $\zeta_{1}$.
Together with the preceding proposition, this gives the calculations stated.
The case of $B^1$ is similar, and the case of $B^\GG= B^0\disjunion B^1$ follows.
\end{proof}

The next calculation is key to much of what follows.
Note that $\omega$ is an $\omega^*$-dimensional bundle
in the sense of \cite{CMW:orientation} and \cite{CW:ordinaryhomology},
so we have the Euler class
$c_{\omega} = e(\omega)\in \Mackey H_\GG^{\omega}(B_+)$.

\begin{proposition}\label{prop:calcEulerClass}
The restriction of $c_{\omega}$ to $B^0$ is
\begin{align*}
 c_{\omega}|B^0 &= c \zeta_1 \in \Mackey H_\GG^\omega(B^0_+). \\
\intertext{The restriction of $c_{\omega}$ to $B^1$ is}
 c_{\omega}|B^1 &= (e^2 + \xi c)\zeta_0^{-1} \in \Mackey H_G^\omega(B^1_+). \\
\intertext{Put another way,}
 c_{\omega}|B^\GG &= (c \zeta_1, (e^2 + \xi c)\zeta_0^{-1}).
\end{align*}
\end{proposition}

\begin{proof}
Consider first the case of $c_{\omega}|B^0$. 
This is the Euler class of the tautological line bundle over $B^0$ with trivial $\GG$-action.
Grading over $RO(\GG)$, this is how $c\in\Mackey H_\GG^2(B^0_+)$ was defined.
When we grade on $RO(\Pi B)$ via the restriction $RO(\Pi B)\to RO(\Pi B^0) = RO(\GG)$,
the corresponding element in grading $\omega = 2 + \Omega_1$ is $c\zeta_1$,
so we have the first equality claimed.

For $c_{\omega}|B^1$, we use the preceding proposition and the calculation of the cohomology
of a point to see that
\[
 \Mackey H_\GG^{\omega}(B^1_+) \iso \conc\Z \dirsum \Mackey\Z
\]
with the first summand generated by $e^2\zeta_0^{-1}$ and the second by
$\xi c\zeta_0^{-1}$.
Thus, $c_{\omega}|B^1 = ae^2\zeta_0^{-1} + b\xi c\zeta_0^{-1}$ for some integers $a$ and $b$.
Restricting to nonequivariant cohomology, $c_{\omega}|B^1$ restricts to the first nonequivariant
Chern class of the tautological bundle, $e^2$ restricts to 0, and $\xi c\zeta_0^{-1}$
also restricts to the first nonequivariant Chern class of the tautological bundle.
Thus, $b=1$.
To determine $a$, consider the (equivariant) restriction to
a single point in $B^1$. This time, $c_{\omega}$
must restrict to the Euler class of the fiber over that point, which
is a copy of $\C^\sigma$, whose Euler class is $e^2\zeta_0^{-1}$.
Because $c$ restricts to 0 at any point, being the Euler class of
the tautological bundle with trivial $\GG$-action,
$ae^2\zeta_0^{-1} + \xi c\zeta_0^{-1}$ restricts to $ae^2\zeta_0^{-1}$,
so we must have $a = 1$, giving the second equality claimed.
\end{proof}

Because $\chi$ exchanges the fixed-point components, we get
the following immediate corollary, where we write $c_{\chi\omega}$ for the
Euler class of $\chi\omega = \omega\tensor_\C \C^\sigma$.
Note that $c_{\chi\omega} = \chi^* c_{\omega}$ because $\chi$ classifies
the operation $-\tensor_\C \C^\sigma$.

\begin{corollary}
$c_{\chi\omega}|B^\GG = ((e^2 + \xi c)\zeta_1^{-1}, c\zeta_0)
   \in \Mackey H_\GG^{\chi\omega}(B^\GG_+)$.
\qed
\end{corollary}

\section{The cohomologies of $B_\GG U(1)_+\smsh (E\GG)_+$ and $B_\GG U(1)_+\smsh\tE\GG$}
\label{sec:egcohomologies}

We start with an easy result.

\begin{proposition}\label{prop:BEGtilde}
The inclusion $B^\GG\to B$ induces an isomorphism
\[
 \Mackey H_\GG^{RO(\Pi B)}(B_+\smsh_B \tE\GG) \iso \Mackey H_\GG^{RO(\Pi B)}(B_+^\GG\smsh_B \tE\GG)
\]
where the latter is calculated in Proposition~\ref{prop:calcFixedPoints}.
\end{proposition}

\begin{proof}
For any ex-$\GG$-space $X$ over $B$, the inclusion $X^\GG\to X$ induces a weak equivalence
$X^\GG\smsh_B\tE\GG\to X\smsh_B\tE\GG$, from which the isomorphism follows.
\end{proof}

In one sense, the calculation of the cohomology of $B_+\smsh_B (E\GG)_+$ is just as easy,
as we have the following result.

\begin{proposition}\label{prop:BEGplus}
Each inclusion $B^k\to B$, $k=0$ or $1$, induces an isomorphism
\[
 \Mackey H_\GG^{RO(\Pi B)}(B_+\smsh_B (E\GG)_+) \iso \Mackey H_\GG^{RO(\Pi B)}(B^k_+\smsh_B (E\GG)_+)
\]
Consequently, the inclusion $B^\GG\to B$ induces a split short exact sequence
\begin{multline*}
 0
 \to \Mackey H_\GG^{RO(\Pi B)}(B_+\smsh_B (E\GG)_+)
 \xrightarrow{\eta} \Mackey H_\GG^{RO(\Pi B)}(B_+^\GG\smsh_B (E\GG)_+) \\
 \xrightarrow{\theta} \susp^{-1}\Mackey H_\GG^{RO(\Pi B)}(B/_B B^\GG\smsh_B (E\GG)_+)
 \to 0.
\end{multline*}
For the last group we have an isomorphism
\[
 \susp^{-1}\Mackey H_\GG^{RO(\Pi B)}(B/_B B^\GG\smsh_B (E\GG)_+)
  \iso \Mackey H_\GG^{RO(\Pi B)}(B^1_+\smsh_B (E\GG)_+).
\]
\end{proposition}

\begin{proof}
For each $k$, the inclusion $B^k\to B$ is a nonequivariant equivalence, hence
$B^k\times E\GG \to B\times E\GG$ is an equivariant equivalence,
giving the first isomorphism.

We can take as a splitting of $\eta$ the map
\begin{multline*}
 \eta'\colon \Mackey H_\GG^{RO(\Pi B)}(B_+^\GG\smsh_B (E\GG)_+) \\
  \iso \Mackey H_\GG^{RO(\Pi B)}(B^0_+\smsh_B (E\GG)_+) \dirsum \Mackey H_\GG^{RO(\Pi B)}(B^1_+\smsh_B (E\GG)_+) \\
  \to \Mackey H_\GG^{RO(\Pi B)}(B^0_+\smsh_B (E\GG)_+) \iso \Mackey H_\GG^{RO(\Pi B)}(B_+\smsh_B (E\GG)_+)
\end{multline*}
given by projection to the $B^0$ summand.
We then have
\[
 \Mackey H_\GG^{RO(\Pi B)}(B/_B B^\GG\smsh_B (E\GG)_+) \iso \susp\ker\eta',
\]
giving the rest of the proposition.
\end{proof}

We will need to know exactly what the inclusion
\[
 \eta\colon \Mackey H_\GG^{RO(\Pi B)}(B_+\smsh_B (E\GG)_+)
  \to \Mackey H_\GG^{RO(\Pi B)}(B^\GG_+\smsh_B (E\GG)_+)
\]
does. First, we write the source in the following way.

\begin{definition}
Let $c_{\omega}\in \Mackey H_\GG^{\omega}(B_+\smsh (E\GG)_+)$ denote the image of
$c_{\omega}\in \Mackey H_\GG^{\omega}(B_+)$ under pullback along the projection $B\times E\GG\to B$.
Let $\zeta_1 \in \Mackey H_\GG^{\Omega_1}(B_+\smsh (E\GG)_+)$
be the element that restricts
to $\zeta_1\in \Mackey H_\GG^{\Omega_1}(B^0_+\smsh (E\GG)_+)$.
\end{definition}

\begin{corollary}\label{cor:BtimesEG}
\[
 \Mackey H_\GG^{RO(\Pi B)}(B_+\smsh_B (E\GG)_+)
 \iso \Mackey H_\GG^{RO(\GG)}((E\GG)_+)[c_{\omega}, \zeta_1, \zeta_1^{-1}].
\]
\end{corollary}

\begin{proof}
This is clear from the isomorphism 
\[
 \Mackey H_\GG^{RO(\Pi B)}(B_+\smsh_B (E\GG)_+) \iso
	\Mackey H_\GG^{RO(\Pi B)}(B^0_+\smsh_B (E\GG)_+)
\]
given by the
preceding proposition, together with the calculation of
Proposition~\ref{prop:calcFixedPoints},
except possibly for the use of $c_{\omega}$ as one
of the generators. This is justified by Proposition~\ref{prop:calcEulerClass},
which tells us that
$c_{\omega}$ maps to $c\zeta_1 \in \Mackey H_\GG^{\omega}(B^0_+\smsh_B (E\GG)_+)$.
\end{proof}

We also define the following element.

\begin{definition}
Let
\[
 \zeta_0 = \xi\zeta_1^{-1} \in \Mackey H_\GG^{\Omega_0}(B_+\smsh_B (E\GG)_+).
\]
\end{definition}

Recalling that $\xi$ is invertible in $\Mackey H_\GG^{RO(\GG)}((E\GG)_+)$, we see that
$\zeta_0$ is invertible in $\Mackey H_\GG^{RO(\Pi B)}(B_+\smsh_B (E\GG)_+)$ and that
\[
 \zeta_0 \cdot \zeta_1 = \xi.
\]

\begin{proposition}\label{prop:imageeta}
Under the inclusion 
\[
 \eta\colon \Mackey H_\GG^{RO(\Pi B)}(B_+\smsh_B (E\GG)_+) \includesin \Mackey H_\GG^{RO(\Pi B)}(B^\GG_+\smsh_B (E\GG)_+)
\]
we have
\begin{align*}
  \eta(c_{\omega}) &= ( c\zeta_1, (e^2 + \xi c)\zeta_0^{-1} ), \\
  \eta(\zeta_0) &= ( \xi\zeta_1^{-1}, \zeta_0 ), \\
 \mathllap{\text{and}\qquad} \eta(\zeta_1) &= ( \zeta_1, \xi\zeta_0^{-1} ).
\end{align*}
\end{proposition}

\begin{proof}
The calculation of $\eta(c_{\omega})$ was done in Proposition~\ref{prop:calcEulerClass}.
We have that 
\begin{align*}
 \zeta_1|B^0 &= \zeta_1 \\
\intertext{and}
 \zeta_0|B^0 &= \xi\zeta_1^{-1}
\end{align*}
by definition.
Now $\zeta_1|B^1$ is a unit in $\Mackey H_\GG^{\Omega_1}(B^1_+\smsh_B (E\GG)_+)$, but this is
a copy of $\Mackey\Z$ generated by $\xi\zeta_0^{-1}$.
Both $\zeta_1|B^1$ and $\xi\zeta_0^{-1}$ restrict to $\iota^2\zeta_0^{-1}$ nonequivariantly, so we must have
$\zeta_1|B^1 = \xi\zeta_0^{-1}$.
The equation $\zeta_0|B^1 = \zeta_0$ then follows from the definition
that $\zeta_0 = \xi\zeta_1^{-1}$.
\end{proof}


\section{Preliminary results on the cohomology of $B_\GG U(1)$}\label{sec:prelimresults}

Many of our arguments in the remainder of the calculation will be based on
diagram~(\ref{diag:main}), which we repeat here, with $R = RO(\Pi B)$ again.
\[\def\objectstyle{\scriptstyle}
 \xymatrix@C-0.5em{
   \Mackey H_\GG^R(B_+\smsh_B\tE\GG) \ar[r]^\eta_\iso \ar[d]_{\psi}
     & \Mackey H_\GG^R(B_+^\GG\smsh_B\tE\GG) \ar[d]^\psi \\
   \Mackey H_\GG^R(B_+) \ar[r]^\eta \ar[d]_\phi
     & \Mackey H_\GG^R(B_+^\GG) \ar[r]^-\theta \ar[d]^\phi
     & \susp^{-1}\Mackey H_\GG^{R}(B/_B B^\GG) \ar[d]^\phi_\iso \\
   \Mackey H_\GG^R(B_+\smsh_B (E\GG)_+) \ar[r]^\eta \ar[d]_\delta
     & \Mackey H_\GG^R(B_+^\GG\smsh_B (E\GG)_+) \ar[r]^-\theta \ar[d]^\delta
     & \susp^{-1}\Mackey H_\GG^{R}(B/_B B^\GG\smsh_B (E\GG)_+) \\
   \susp^{-1}\Mackey H_\GG^{R}(B_+\smsh_B\tE\GG) \ar[r]^\eta_\iso
     & \susp^{-1}\Mackey H_\GG^{R}(B_+^\GG\smsh_B\tE\GG)
  }
\]
The rows and columns are all parts of long exact sequences.
That the map $\eta$ in the top and bottom rows of the diagram is an isomorphism
is Proposition~\ref{prop:BEGtilde}.
That implies that $\Mackey H_\GG^{RO(\Pi B)}(B/_B B^\GG\smsh_B\tE\GG)$,
the group that would appear in the top and bottom right corners, is zero;
alternatively, this vanishing follows from the fact that
$(B/_B B^\GG)^\GG$ is the trivial ex-space over $B^\GG$.
In turn, this vanishing implies that
the map $\phi$ on the far right is an isomorphism.
In the third row, Proposition~\ref{prop:BEGplus} shows that $\eta$ is
a monomorphism and $\theta$ is an epimorphism.

\begin{proposition}\label{prop:etamono}
$\eta\colon \Mackey H_\GG^\alpha(B_+)\to \Mackey H_\GG^\alpha(B_+^\GG)$ is a monomorphism
for $|\alpha| \leq 0$.
It is also a monomorphism at level $\GG/\GG$ if $|\alpha| > 0$ and 
$\alpha_0^\GG$ and $\alpha_1^\GG$ are even.
\end{proposition}

\begin{proof}
From the long exact sequence, we see that $\eta$ is a monomorphism if the group
$\Mackey H_\GG^{\alpha}(B/_B B^\GG) = 0$.
In the diagram above we noted that
\[
	\Mackey H_\GG^{RO(\Pi B)}(B/_B B^\GG) \iso \Mackey H_\GG^{RO(\Pi B)}(B/_B B^\GG\smsh_B (E\GG)_+)
\]
and we calculated the latter in
Proposition~\ref{prop:BEGplus}.
As a module over the ring $\Mackey H_\GG^{RO(\GG)}((E\GG)_+)$, 
it is the suspension of an algebra with multiplicative generators
in gradings $2$ and $\Omega_0$.
Further, $\Mackey H_\GG^\beta((E\GG)_+)$ is zero
if $|\beta|< 0$; at level $\GG/\GG$ it is also zero for
those $|\beta|\geq 0$ for which $\beta^\GG$ is odd.
It follows that
\[
 \Mackey H_\GG^{\alpha}(B/_B B^\GG) \iso \Mackey H_\GG^{\alpha}(B/_B B^\GG\smsh_B (E\GG)_+) = 0
\]
for the $\alpha$ specified in the statement of the proposition.
\end{proof}

The following is then a diagram chase.

\begin{corollary}\label{cor:pullback}
The diagram
\[
 \xymatrix{
  H_\GG^\alpha(B_+) \ar[r]^\eta \ar[d]_\phi
  	& H_\GG^\alpha(B^\GG_+) \ar[d]^\phi \\
  H_\GG^\alpha(B_+\smsh_B (E\GG)_+) \ar@{>->}[r]^\eta
  	& H_\GG^\alpha(B^\GG_+\smsh_B(E\GG)_+)
 }
\]
is a pullback square for $|\alpha| < 0$ and also for
those $|\alpha|>0$
for which $\alpha_0^\GG$ and $\alpha_1^\GG$ are even.
\qed
\end{corollary}

\begin{corollary}\label{cor:xibar}
There are unique elements $\zeta_0\in \Mackey H_\GG^{\Omega_0}(B_+)$
and $\zeta_1\in \Mackey H_\GG^{\Omega_1}(B_+)$ such that
\begin{align*}
	\eta(\zeta_0) &= (\xi\zeta_1^{-1}, \zeta_0), & \phi(\zeta_0) &= \zeta_0, \\
	\eta(\zeta_1) &= (\zeta_1, \xi\zeta_0^{-1}), & \phi(\zeta_1) &= \zeta_1.
\end{align*}
\end{corollary}

\begin{proof}
By Proposition~\ref{prop:imageeta},
the element $\zeta_0 \in \Mackey H_\GG^{\Omega_0}(B_+\smsh_B (E\GG)_+)$
satisfies $\eta(\zeta_0) = (\xi\zeta_1^{-1}, \zeta_0)$.
On the other hand, we have the element $(\xi\zeta_1^{-1}, \zeta_0) \in \Mackey H_\GG^{\Omega_0}(B^\GG_+)$,
with $\phi(\xi\zeta_1^{-1}, \zeta_0) = \eta(\zeta_0)$.
It follows from Corollary~\ref{cor:pullback} that there is a unique element
$\zeta_0\in \Mackey H_\GG^{\Omega_0}(B_+)$ with the claimed images under $\eta$ and $\phi$.
The existence and uniqueness of $\zeta_1$ follow in the same way.
\end{proof}

We emphasize that $\zeta_0$ and $\zeta_1$ are not invertible as elements
of $\Mackey H_\GG^{RO(\Pi B)}(B_+)$, unlike their images in 
$\Mackey H_\GG^{RO(\Pi B)}(B_+\smsh_B (E\GG)_+)$---their images in
$\Mackey H_\GG^{RO(\Pi B)}(B^\GG_+)$ are not invertible.

The action of $\chi$ on these elements is given as follows.

\begin{proposition}
In $\Mackey H_\GG^{RO(\Pi B)}(B_+)$ we have
\begin{align*}
	\chi(\zeta_0) &= \zeta_1, & \chi(\zeta_1) &= \zeta_0, \\
	\chi(c_{\omega}) &= c_{\chi\omega}, & \chi(c_{\chi\omega}) &= c_{\omega}.
\end{align*}
\end{proposition}

\begin{proof}
These elements all live in gradings in which $\eta$ is a monomorphism, so it suffices to
check that the equalities are true after applying $\eta$. We have
\[
 \eta(\chi\zeta_0) = \chi\eta(\zeta_0)
 	= \chi(\xi\zeta_1^{-1}, \zeta_0)
	= (\zeta_1, \xi\zeta_0^{-1})
	= \eta(\zeta_1),
\]
hence $\chi\zeta_0 = \zeta_1$. That $\chi\zeta_1 = \zeta_0$ follows because $\chi^2 = 1$.

We have already noted that $\chi c_\omega = c_{\chi\omega}$ because $\chi$ classifies
$-\tensor_\C \C^\sigma$.
\end{proof}

\section{The proposed ring structure}\label{sec:ringstructure}

In \S\ref{sec:structure} we shall show that 
$\Mackey H_\GG^{RO(\Pi B)}(B_+)$ is generated multiplicatively by the elements
$\zeta_0$, $\zeta_1$, $c_{\omega}$, and $c_{\chi\omega}$, subject to the two relations given in the following result.

\begin{proposition}\label{prop:relations}
In $\Mackey H_\GG^{RO(\Pi B)}(B_+)$ we have the relations
\[
  \zeta_0 \zeta_1 = \xi
\]
and
\[
  \zeta_1 c_{\chi\omega} = (1-\kappa)\zeta_0 c_{\omega} + e^2.
\]
\end{proposition}

\begin{proof}
To show that the relations hold, we note that both take place in gradings where
Proposition~\ref{prop:etamono} says that 
$\eta\colon \Mackey H_G^{RO(\Pi B)}(B_+)\to \Mackey H_G^{RO(\Pi B)}(B^G_+)$ is a monomorphism,
so it suffices to show that the relations hold after applying $\eta$.
Recall that
\begin{align*}
 \eta(\zeta_0) &= (\xi\zeta_1^{-1}, \zeta_0), \\
 \eta(\zeta_1) &= (\zeta_1, \xi\zeta_0^{-1}), \\
 \eta(c_{\omega}) &= ( c\zeta_1, (e^2 + \xi c)\zeta_0^{-1}), \quad\text{and} \\
 \eta(c_{\chi\omega}) &= ( (e^2+\xi c)\zeta_1^{-1}, c\zeta_0).
\end{align*}
We see immediately that
\[
 \eta(\zeta_0\zeta_1) = (\xi,\xi) = \eta(\xi),
\]
hence $\zeta_0\zeta_1 = \xi$. For the second relation, we have
\begin{align*}
 \eta(\zeta_1 c_{\chi\omega} - (1-\kappa)\zeta_0 c_{\omega})
  &= (e^2 + \xi c, \xi c) - (1-\kappa)(\xi c, e^2 + \xi c) \\
  &= (e^2 + \xi c, \xi c) - (\xi c, -e^2 + \xi c) \\
  &= (e^2, e^2) \\
  &= \eta(e^2),
\end{align*}
using that $(1-\kappa)\xi = \xi$ and $(1-\kappa)e = -e$, hence the second relation also holds.
\end{proof}

Multiplying the second relation by $(1-\kappa)$ and rearranging gives the similar relation
\[
 \zeta_0 c_{\omega} = (1-\kappa) \zeta_1 c_{\chi\omega} + e^2.
\]

\begin{definition}
Let $\Mackey P^{RO(\Pi B)}$ be the algebra defined by
\[
 \Mackey P^{RO(\Pi B)}
  = \Mackey H_\GG^{RO(\GG)}(S^0)[\zeta_0, \zeta_1, c_{\omega}, c_{\chi\omega}]/
   \langle \zeta_0\zeta_1 - \xi, \zeta_1 c_{\chi\omega} - (1-\kappa)\zeta_0 c_{\omega} - e^2 \rangle.
\]
\end{definition}

In the following section we shall show that this is isomorphic to the cohomology of $B$.
The main result of this section is the following.

\begin{theorem}\label{thm:freeness}
$\Mackey P^{RO(\Pi B)}$ is a free $\Mackey H_\GG^{RO(\GG)}(S^0)$-module on a basis consisting of the images of
those monomials $\zeta_0^k \zeta_1^\ell c_{\omega}^m c_{\chi\omega}^n$ that are {\em not} multiples of
\begin{itemize}
\item $\zeta_0\zeta_1$,
\item $\zeta_1 c_{\chi\omega}$, or
\item $\zeta_0^2 c_{\omega}$.
\end{itemize}
\end{theorem}

\begin{proof}
Write $I = \langle \zeta_0\zeta_1 - \xi, \zeta_1 c_{\chi\omega} - (1-\kappa)\zeta_0 c_{\omega} - e^2 \rangle$.
Rather than ``write out tortuous verifications'' \cite[Introduction]{Berg:diamondLemma},
we use Bergman's diamond lemma, \cite[Theorem~1.2]{Berg:diamondLemma}, as modified
for the commutative case by the comments in his \S10.3.
The diamond lemma will show that, using the relations, every element of $\Mackey P^{RO(\Pi B)}$ can be reduced
in a unique way to a linear combination of the basis elements specified in
the statement of the theorem.
(Note: We should presumably be working with {\em anti-}commutative rings and algebras, but
Theorem~\ref{thm:evenCohomPoint} tells us that the cohomology of a point is strictly commutative,
and our polynomial generators $\zeta_0$, etc., are in ``even'' gradings where 
no signs will be introduced.)
To use the diamond lemma we need several things:
First, we need a partial ordering $\preccurlyeq$ on the set of all monomials
in $\Mackey H_\GG^{RO(\GG)}(S^0)[\zeta_0, \zeta_1, c_{\omega}, c_{\chi\omega}]$
such that, if $A \preccurlyeq B$ then $CA \preccurlyeq CB$ for any monomial $C$.
(By a monomial we shall always mean a product $\zeta_0^k \zeta_1^\ell c_{\omega}^m c_{\chi\omega}^n$
with no coefficient.)
We define $\preccurlyeq$ by ordering by total degree and, within the same degree, saying that 
\[
 \zeta_0^{k_1} \zeta_1^{\ell_1} c_{\omega}^{m_1} c_{\chi\omega}^{n_1}
  \prec \zeta_0^{k_2} \zeta_1^{\ell_2} c_{\omega}^{m_2} c_{\chi\omega}^{n_2}
\]
if $n_1 < n_2$. Note that, as required for the diamond lemma, this partial order 
satisfies the descending chain condition.

Next, we need a {\em reduction system} $S$, consisting of pairs $(W,f)$, 
where $W$ is a monomial, $f$ is a polynomial,
and $I$ is generated by the collection of elements $W-f$, $(W,f)\in S$.
We take $S$ to be the set containing the three pairs
\begin{align*}
	\sigma_1 &= (W_1, f_1) = (\zeta_0\zeta_1, \xi ) \\
	\sigma_2 &= (W_2, f_2) = (\zeta_1 c_{\chi\omega}, (1-\kappa)\zeta_0 c_{\omega} + e^2) \qquad\text{and} \\
	\sigma_3 &= (W_3, f_3) = (\zeta_0^2 c_{\omega}, \xi c_{\chi\omega} + e^2\zeta_0).
\end{align*}
The first two satisfy $W-f \in I$ by the definition of $I$; for the third, note that
\begin{align*}
 \zeta_0^2 c_\omega &= \zeta_0(\zeta_0 c_{\omega}) \\
  &\equiv \zeta_0((1-\kappa)\zeta_1 c_{\chi\omega} + e^2) \pmod I \\
  &\equiv \xi c_{\chi\omega} + e^2\zeta_0 \pmod I.
\end{align*}
Note that $W_1-f_1$ and $W_2-f_2$ already generate $I$, so it does no harm to add $\sigma_3$.
(We will see in a moment why we need it.)
We also need that the partial ordering on the monomials is compatible with the reduction system, meaning that,
for each pair $(W_i,f_i)$, $f_i$ is a linear combination of monomials $\prec W_i$,
which the reader can quickly check.

As discussed in \cite{Berg:diamondLemma}, $S$ defines reduction maps on 
$\Mackey H_\GG^{RO(\GG)}(S^0)[\zeta_0, \zeta_1, c_{\omega}, c_{\chi\omega}]$: Given $(W,f)$,
we define the reduction $r$ on monomials by defining $r(AW) = Af$,
$r(A) = A$ if $A$ is not a multiple of $W$, and then extending linearly to polynomials.
Write $r_1$, $r_2$, and $r_3$ for the reduction maps defined by $\sigma_1$, $\sigma_2$, and $\sigma_3$,
respectively.
The last thing we need to verify in order to use the diamond lemma
is that {\em all ambiguities of the reduction system are resolvable,}
meaning that,
if $(W_i, f_i)$ and $(W_j, f_j)$ are pairs in our reduction system, with corresponding reductions
$r_i$ and $r_j$, and
$W$ is the least common multiple of $W_i$ and $W_j$, then $r_i(W)$ and $r_j(W)$
can be further reduced to give the same polynomial. There are three cases we need to check:
\begin{enumerate}
\item $\sigma_1$ and $\sigma_2$: $\zeta_0\zeta_1 c_{\chi\omega}$ is the least common multiple of $\zeta_0\zeta_1$ and
$\zeta_1 c_{\chi\omega}$. In this case, we have
\[
 r_1(\zeta_0\zeta_1 c_{\chi\omega}) = \xi c_{\chi\omega}
\]
and
\[
 r_2(\zeta_0\zeta_1 c_{\chi\omega}) = (1-\kappa)\zeta_0^2 c_{\omega} + e^2\zeta_0.
\]
Applying $r_3$ to the latter gives
\[
 r_3 r_2(\zeta_0\zeta_1 c_{\chi\omega}) = (1-\kappa)\xi c_{\chi\omega} + (1-\kappa)e^2\zeta_0 + e^2\zeta_0
  = \xi c_{\chi\omega},
\]
resolving the ambiguity. (This is why we need to include $\sigma_3$ in our reduction system.)

\item $\sigma_1$ and $\sigma_3$: $\zeta_0^2\zeta_1 c_{\omega}$ is the least common multiple in question, and
\[
 r_1(\zeta_0^2\zeta_1 c_{\omega}) = \xi\zeta_0 c_{\omega}
\]
while
\[
 r_3(\zeta_0^2\zeta_1 c_{\omega}) = \xi\zeta_1 c_{\chi\omega} + e^2\zeta_0\zeta_1.
\]
Applying $r_1 r_2$ to the second gives
\begin{align*}
 r_1 r_2 r_3(\zeta_0^2\zeta_1 c_{\omega})
  &= r_1( \xi\zeta_0 c_{\omega} + e^2\xi + e^2\zeta_0\zeta_1) \\
  &= \xi\zeta_0 c_{\omega} + 2e^2\xi \\
  &= \xi\zeta_0 c_{\omega},
\end{align*}
resolving this ambiguity.

\item $\sigma_2$ and $\sigma_3$: $\zeta_0^2\zeta_1 c_{\omega}c_{\chi\omega}$ is the least common multiple,
\[
 r_2(\zeta_0^2\zeta_1 c_{\omega}c_{\chi\omega}) = (1-\kappa)\zeta_0^3 c_{\omega}^2 + e^2\zeta_0^2 c_{\omega},
\]
and
\[
 r_3(\zeta_0^2\zeta_1 c_{\omega}c_{\chi\omega}) = \xi\zeta_1 c_{\chi\omega}^2 + e^2\zeta_0\zeta_1 c_{\chi\omega}.
\]
In this case, we find that
\begin{align*}
 r_3 r_2(\zeta_0^2\zeta_1 c_{\omega}c_{\chi\omega})
  &= \xi\zeta_0 c_{\omega}c_{\chi\omega} - e^2\zeta_0^3 c_{\omega} + e^2\xi c_{\chi\omega} + e^4\zeta_0 \\
  &= r_2 r_3(\zeta_0^2\zeta_1 c_{\omega}c_{\chi\omega}).
\end{align*}
\end{enumerate}
Thus, we can resolve all the ambiguities, which completes our verification of the hypotheses of the diamond lemma.

The conclusion of the diamond lemma is that, as a module over $\Mackey H_\GG^{RO(\GG)}(S^0)$,
$\Mackey H_\GG^{RO(\GG)}(S^0)[\zeta_0, \zeta_1, c_{\omega}, c_{\chi\omega}]$ is the direct sum of $I$
and the submodule generated by the {\em irreducible} monomials, which are precisely those not divisible
by $W_1$, $W_2$, or $W_3$. Thus, the quotient ring $\Mackey P^{RO(\Pi B)}$ is free on the images of the
irreducible monomials, proving the theorem.
\end{proof}

\begin{definition}\label{def:basic}
We call a monomial in $\zeta_0$, $\zeta_1$, $c_{\omega}$ and $c_{\chi\omega}$ {\em basic}
if it is not a multiple of
\begin{itemize}
\item $\zeta_0\zeta_1$,
\item $\zeta_1 c_{\chi\omega}$, or
\item $\zeta_0^2 c_{\omega}$.
\end{itemize}
\end{definition}

The preceding theorem says that the basic monomials form a basis for $\Mackey P^{RO(\Pi B)}$
over the cohomology of a point.
In fact, the diamond lemma provides an algorithm for reducing any element to ``normal form,'' that is,
to write any element in terms of basic monomials: Simply apply the reductions $r_1$, $r_2$, and $r_3$
in any order until no further reductions can be achieved. This process is guaranteed to stop after
finitely many steps, with the end result not depending on the order in which the reductions are applied.

The following result lists the basic monomials more explicitly, organized by cosets of $RO(\GG)$ in $RO(\Pi B)$.

\begin{corollary}\label{cor:basic}
If $n > 0$, then the basic monomials with gradings in $n\omega+RO(\GG)$ are
\[
 \{ \zeta_1^n, \zeta_1^{n-1} c_{\omega}, \ldots, \zeta_1 c_{\omega}^{n-1}, c_{\omega}^n,
 	\zeta_0 c_{\omega}^{n+1}, c_{\omega}^{n+1} c_{\chi\omega}, \zeta_0 c_{\omega}^{n+2} c_{\chi\omega},
	c_{\omega}^{n+2} c_{\chi\omega}^2, \ldots \}.
\]
The basic monomials with gradings in $RO(\GG)$ are
\[
 \{ 1, \zeta_0 c_{\omega}, c_{\omega}c_{\chi\omega}, \zeta_0 c_{\omega}^2 c_{\chi\omega},
 	c_{\omega}^2 c_{\chi\omega}^2, \zeta_0 c_{\omega}^3 c_{\chi\omega}^2, \ldots \}.
\]
If $n < 0$, the basic monomials with gradings in $n\omega + RO(\GG)$ are
\[
 \{ \zeta_0^{|n|}, \zeta_0^{|n|-1} c_{\chi\omega}, \ldots, \zeta_0 c_{\chi\omega}^{|n|-1},
 	c_{\chi\omega}^{|n|}, \zeta_0 c_{\omega} c_{\chi\omega}^{|n|},
	c_{\omega}c_{\chi\omega}^{|n|+1}, \zeta_0 c_{\omega}^2 c_{\chi\omega}^{|n|+1},
	c_{\omega}^2 c_{\chi\omega}^{|n|+2}, \ldots \}.
\]
In each case, the monomials are listed in order so that the grading of each
is either $2$ or $2\sigma$ larger than the grading of the one preceding it.
\qed
\end{corollary}

The proof is an exercise in computing the gradings of the basic monomials and
arranging them in order, and is omitted.

\section{The cohomology of $B_\GG U(1)$}\label{sec:structure}

There is an algebra map
\[
 f\colon \Mackey P^{RO(\Pi B)} \to \Mackey H_\GG^{RO(\Pi B)}(B_+)
\]
defined by taking each of the generators $\zeta_0$, $\zeta_1$, $c_{\omega}$, and $c_{\chi\omega}$ to the
element of the same name.
Proposition~\ref{prop:relations} shows that this does define an algebra map.
The goal of this section is to show that $f$ is an isomorphism.

There are several ways the proof could proceed.
One would be to take the approach used in \cite{CHT:projspaces}
and verify that the basis found in Theorem~\ref{thm:freeness}
maps to a basis for the nonequivariant cohomology of $B$ and to a basis
for the nonequivariant cohomology of $B^\GG$.
By the argument given in that paper, that suffices to show that these elements
form a basis for the equivariant cohomology of $B$.
This could be carried out using the explicit description of the basis
given in Corollary~\ref{cor:basic}.

We will give a different argument here, partly for variety but, more importantly,
because this argument will generalize to allow a similar computation
in Part~\ref{part:basespace}, where the cohomology being computed
is not free, hence the preceding argument would not work.

What we shall do is the following.
Because $\Mackey P^{RO(\Pi B)}$ is a free $\Mackey H_\GG^{RO(G)}(S^0)$-module, 
tensoring with the long exact sequence of the cofibration
$(E\GG)_+\to S^0\to \tE\GG$ gives us a long exact sequence
\begin{multline*}
 \cdots\to \Mackey P^{RO(\Pi B)}\tensorS \Mackey H_\GG^{RO(\GG)}(\tE\GG)
  \to \Mackey P^{RO(\Pi B)} \\
  \to \Mackey P^{RO(\Pi B)}\tensorS \Mackey H_\GG^{RO(\GG)}((E\GG)_+)
  \to \cdots
\end{multline*}
The map $f$ induces the following map of long exact sequences:
\[
 \xymatrix{
 \vdots \ar[d] & \vdots \ar[d] \\
 \Mackey P^{RO(\Pi B)}\tensorS \Mackey H_\GG^{RO(\GG)}(\tE\GG) \ar[r]^-{f_{\tE\GG}} \ar[d]
 	& \Mackey H_\GG^{RO(\Pi B)}(B_+\smsh_B \tE\GG) \ar[d] \\
 \Mackey P^{RO(\Pi B)} \ar[r]^f \ar[d]
	& \Mackey H_\GG^{RO(\Pi B)}(B_+) \ar[d] \\
 \Mackey P^{RO(\Pi B)}\tensorS \Mackey H_\GG^{RO(\GG)}((E\GG)_+) \ar[r]^-{f_{(E\GG)_+}} \ar[d]
	& \Mackey H_\GG^{RO(\Pi B)}(B_+\smsh_B (E\GG)_+) \ar[d] \\
 \vdots & \vdots
 }
\]
We will show that both $f_{(E\GG)_+}$ and $f_{\tE\GG}$ are isomorphisms, 
so conclude that $f$ is an isomorphism.
The first isomorphism is straightforward, the second is less so.

\begin{proposition}\label{prop:fEG+}
\[
 f_{(E\GG)_+}\colon \Mackey P^{RO(\Pi B)}\tensorS \Mackey H_\GG^{RO(\GG)}((E\GG)_+)
\to \Mackey H_\GG^{RO(\Pi B)}(B_+\smsh_B (E\GG)_+)
\]
is an isomorphism.
\end{proposition}

\begin{proof}
Write
\begin{multline*}
  \Mackey P^{RO(\Pi B)}\tensorS \Mackey H_\GG^{RO(\GG)}((E\GG)_+) \\
   \iso \Mackey H_\GG^{RO(\GG)}((E\GG)_+)[\zeta_0, \zeta_1, c_{\omega}, c_{\chi\omega}] /
   \langle \zeta_0\zeta_1 - \xi, \zeta_1 c_{\chi\omega} - \zeta_0 c_{\omega} - e^2 \rangle
\end{multline*}
and remember that $\Mackey H_\GG^{RO(\GG)}((E\GG)_+) \iso \Mackey H_\GG^{RO(\GG)}(S^0)[\xi^{-1}]$.
(Note that $\kappa = 0$ and $2e^2 = 0$ in this ring.)
The relation
$\zeta_0\zeta_1 = \xi$ tells us that both $\zeta_0$ and $\zeta_1$ are invertible in this ring
and that
\[
 \zeta_0 = \xi\zeta_1^{-1}.
\]
Further, the second relation can be written as
\[
 c_{\chi\omega} = \zeta_0\zeta_1^{-1} c_{\omega} + \zeta_1^{-1}e^2
  = \xi\zeta_1^{-2} c_{\omega} + \zeta_1^{-1}e^2.
\]
Thus, we have that
\[
 \Mackey P^{RO(\Pi B)}\tensorS \Mackey H_\GG^{RO(\GG)}((E\GG)_+)
  \iso \Mackey H_\GG^{RO(\GG)}((E\GG)_+)[c_{\omega}, \zeta_1, \zeta_1^{-1}].
\]
The isomorphism now follows from Corollary~\ref{cor:BtimesEG}.
\end{proof}

\begin{proposition}\label{prop:ftEG}
\[
 f_{\tE\GG}\colon \Mackey P^{RO(\Pi B)}\tensorS \Mackey H_G^{RO(G)}(\tE\GG)
	\to \Mackey H_G^{RO(\Pi B)}(B_+\smsh_B \tE\GG)
\]
is an isomorphism.
\end{proposition}

\begin{proof}
By Proposition~\ref{prop:BEGtilde}, we can write $f_{\tE\GG}$ as
\begin{multline*}
  f_{\tE\GG}\colon \Mackey P^{RO(\Pi B)}\tensorS \Mackey H_\GG^{RO(\GG)}(\tE\GG) \\
  \to \Mackey H_\GG^{RO(\GG)}(S^0)[c,\zeta_1,\zeta_1^{-1}]
 	\tensorS \Mackey H_\GG^{RO(\GG)}(\tE\GG) \\
 \dirsum
 \Mackey H_\GG^{RO(\GG)}(S^0)[c,\zeta_0,\zeta_0^{-1}]
 	\tensorS \Mackey H_G^{RO(G)}(\tE\GG)
\end{multline*}
and, by previous calculations, it is induced by the map
\[
 \Mackey P^{RO(\Pi B)}
 \to 
 \Mackey H_\GG^{RO(\GG)}(S^0)[c,\zeta_1,\zeta_1^{-1}]
  \dirsum \Mackey H_\GG^{RO(\GG)}(S^0)][c,\zeta_0,\zeta_0^{-1}]
\]
given by
\begin{align*}
  \zeta_0 &\mapsto ( \xi\zeta_1^{-1}, \zeta_0 ) \\
  \zeta_1 &\mapsto (\zeta_1, \xi\zeta_0^{-1} ) \\
  c_{\omega} &\mapsto ( c\zeta_1, (e^2 + \xi c)\zeta_0^{-1} ) \\
  c_{\chi\omega} &\mapsto ( (e^2 + \xi c)\zeta_1^{-1}, c\zeta_0 ).
\end{align*}
Because $e$ acts by isomorphisms on $\Mackey H_\GG^{RO(\GG)}(\tE\GG)$, we may as well invert $e$ on both
sides of this map and consider the resulting map of 
modules over 
\[
 \Mackey H_\GG^{RO(\GG)}(S^0)[e^{-1}]
  \iso \conc\Z [\xi, e, e^{-1}] / \langle 2\xi \rangle.
\]
Note that $\kappa = 2$ after inverting $e$, so the second relation in the definition of
$\Mackey P^{RO(\Pi B)}$ can be written
\[
 \zeta_0 c_{\omega} + \zeta_1 c_{\chi\omega} = e^2.
\]

Now, because every element of $\Mackey H_\GG^{RO(\GG)}(\tE\GG)$ is annihilated by a sufficiently
high power of $\xi$, we have that
\begin{multline*}
 \Mackey P^{RO(\Pi B)}\tensorS \Mackey H_\GG^{RO(\GG)}(\tE\GG) \\
  \iso \Mackey P^{RO(\Pi B)}[e^{-1}]^\wedge_\xi
   \tensor_{\Mackey H_\GG^{RO(\GG)}(S^0)[e^{-1}]} \Mackey H_\GG^{RO(\GG)}(\tE\GG),
\end{multline*}
where ${}^\wedge_\xi$ indicates completion at $\xi$,
and similarly for the right-hand side of $f_{\tE\GG}$. We therefore consider the map
\begin{multline*}
  \bar f_{\tE\GG}\colon \Mackey P^{RO(\Pi B)}[e^{-1}]^\wedge_\xi \\
  \to \Mackey H_\GG^{RO(\GG)}(S^0)[e^{-1}][c,\zeta_1,\zeta_1^{-1}]^\wedge_\xi
 \dirsum
 \Mackey H_\GG^{RO(\GG)}(S^0)[e^{-1}][c,\zeta_0,\zeta_0^{-1}]^\wedge_\xi.
\end{multline*}
To show that $f_{\tE\GG}$ is an isomorphism, it suffices to show that $\bar f_{\tE\GG}$ is,
and to do this we find an explicit inverse.

To that end, define
\[
 c = \sum_{n=0}^\infty (e^{-2}\xi)^{2^n-1}(e^{-2}c_{\omega}c_{\chi\omega})^{2^n}
  \in  \Mackey P^{RO(\Pi B)}[e^{-1}]^\wedge_\xi.
\]
We claim that this element satisfies $\bar f_{\tE\GG}(c) = (c, c)$.
First calculate
\begin{align*}
 \bar f_{\tE\GG}(e^{-2}c_{\omega}c_{\chi\omega}) 
  &= e^{-2}(c(e^2 + \xi c), c(e^2+\xi c)) \\
  &= (c, c) + e^{-2}\xi(c, c)^2,
\end{align*}
so
\begin{align*}
 \bar f_{\tE\GG}((e^{-2}c_{\omega}c_{\chi\omega})^{2^n})
  &= ((c, c) + e^{-2}\xi(c, c)^2)^{2^n} \\
  &= (c, c)^{2^n} + (e^{-2}\xi)^{2^n}(c, c)^{2^{n+1}},
\end{align*}
using that $2\xi = 0$. It follows that, on applying $\bar f_{\tE\GG}$, the infinite
sum in the definition of $c$ telescopes and converges to $(c, c)$.

We then have
\begin{align*}
 \bar f_{\tE\GG}(c_{\omega} - \zeta_1 c) &= (0, e^2\zeta_0^{-1}) \\
\intertext{and}
 \bar f_{\tE\GG}(c_{\chi\omega} - \zeta_0 c) &= (e^2\zeta_1^{-1},0).
\end{align*}
Define
\begin{multline*}
  \bar g \colon \Mackey H_\GG^{RO(\GG)}(S^0)[e^{-1}][c,\zeta_1,\zeta_1^{-1}]^\wedge_\xi
 \dirsum
 \Mackey H_\GG^{RO(\GG)}(S^0)[e^{-1}][c,\zeta_0,\zeta_0^{-1}]^\wedge_\xi \\
 \to \Mackey P^{RO(\Pi B)}[e^{-1}]^\wedge_\xi
\end{multline*}
as the continuous algebra map determined by
\begin{align*}
 \bar g(1,0) &= e^{-2}\zeta_1(c_{\chi\omega} - \zeta_0 c) \\
 \bar g(\zeta_1,0) &= e^{-2}\zeta_1^2(c_{\chi\omega} - \zeta_0 c) \\
 \bar g(\zeta_1^{-1},0) &= e^{-2}(c_{\chi\omega} - \zeta_0 c) \\
 \bar g(c,0) &= e^{-2}\zeta_1 c(c_{\chi\omega} - \zeta_0 c) \\
 \bar g(0,1) &= e^{-2}\zeta_0(c_{\omega} - \zeta_1 c) \\
 \bar g(0,\zeta_0) &= e^{-2}\zeta_0^2(c_{\omega} - \zeta_1 c) \\
 \bar g(0,\zeta_0^{-1}) &= e^{-2}(c_{\omega} - \zeta_1 c) \\
 \bar g(0,c) &= e^{-2}\zeta_0 c(c_{\omega} - \zeta_1 c).
\end{align*}
It is now a matter of calculation to check that $\bar g$ is the inverse of $\bar f_{\tE\GG}$.
\end{proof}

Finally, we get to our main result.

\begin{theorem}\label{thm:multstructure}
$\Mackey H_\GG^{RO(\Pi B)}(B_+)$ is the commutative $\Mackey H_\GG^{RO(\GG)}(S^0)$-algebra generated by
the elements $\zeta_0$, $\zeta_1$, $c_{\omega}$, and $c_{\chi\omega}$ subject to the two relations
\[
  \zeta_0 \zeta_1 = \xi
\]
and
\[
  \zeta_1 c_{\chi\omega} = (1-\kappa)\zeta_0 c_{\omega} + e^2.
\]
It is free as a module over $\Mackey H_\GG^{RO(\GG)}(S^0)$ on the basic monomials
given in Definition~\ref{def:basic} and enumerated in Corollary~\ref{cor:basic}.
\end{theorem}

\begin{proof}
It follows from Propositions~\ref{prop:fEG+} and~\ref{prop:ftEG}
and the comparison of long exact sequences outlined at the beginning of this section that
$f\colon \Mackey P^{RO(\Pi B)} \to \Mackey H_\GG^{RO(\Pi B)}(B_+)$ is a ring isomorphism.
The theorem then follows from Theorem~\ref{thm:freeness}.
\end{proof}

\begin{corollary}
The pairings
\[
 \Mackey H_\GG^{RO(\Pi B)}(B_+)\tensorS \Mackey H_\GG^{RO(\GG)}((E\GG)_+) \to
  \Mackey H_\GG^{RO(\Pi B)}(B_+\smsh_B (E\GG)_+)
\]
and
\[
 \Mackey H_\GG^{RO(\Pi B)}(B_+)\tensorS \Mackey H_\GG^{RO(\GG)}(\tE\GG) \to
  \Mackey H_\GG^{RO(\Pi B)}(B_+\smsh_B \tE\GG)
\]
are isomorphisms. The long exact sequence for the cofibration
$B_+\smsh_B (E\GG)_+ \to B_+ \to B_+\smsh_B \tE\GG$ is isomorphic to the long exact sequence
for the cofibration $(E\GG)_+\to S^0\to \tE\GG$ tensored with $\Mackey H_\GG^{RO(\Pi B)}(B_+)$.
\end{corollary}

\begin{proof}
This follows from the isomorphism $\Mackey P^{RO(\Pi B)} \iso \Mackey H_\GG^{RO(\Pi B)}(B_+)$
and the isomorphism of the long exact sequences displayed above Proposition~\ref{prop:fEG+}.
\end{proof}

The cohomology group $\Mackey H_\GG^0(B_+)$ has an interesting structure. Additively, it is
$\Mackey A \dirsum \conc\Z$, where the second summand is generated by
\[
 \epsilon = e^{-2}\kappa \zeta_0 c_{\omega}.
\]
The following result is needed in the computation to follow of the Euler classes of the duals
of $\omega$ and $\chi\omega$.

\begin{proposition}
The units in $\Mackey H_\GG^0(B_+)$ are
\[
 \pm 1,\ \pm (1-\kappa),\ \pm (1-\epsilon), \text{ and } \pm (1-\kappa + \epsilon) = \pm(1-\kappa)(1-\epsilon).
\]
Each of these elements squares to $1$.
\end{proposition}

\begin{proof}
We know that $\kappa^2 = 2\kappa$ and we also have
\begin{align*}
 \epsilon^2 &= (e^{-2}\kappa \zeta_0 c_{\omega})^2 \\
 	&= 2e^{-4}\kappa (\zeta_0 c_{\omega})^2 \\
	&= 2e^{-4}\kappa \zeta_0 c_{\omega}(e^2 + (1-\kappa)\zeta_1 c_{\chi\omega}) \\
	&= 2e^{-2}\kappa \zeta_0 c_{\omega} + 2(1-\kappa)e^{-4}\kappa\xi c_{\omega}c_{\chi\omega} \\
	&= 2e^{-2}\kappa \zeta_0 c_{\omega} \\
	&= 2\epsilon,
\end{align*}
using the fact that $e^{-4}\kappa \cdot \xi = 0$. From this it follows that
$(1-\kappa)^2 = 1$ and $(1-\epsilon)^2 = 1$.
We also then have that $[(1-\kappa)(1-\epsilon)]^2 = 1$, but
$(1-\kappa)(1-\epsilon) = 1-\kappa+\epsilon$ because $\kappa\epsilon = 2\epsilon$.
Therefore, the elements listed all square to 1 and are all units.

On the other hand, an arbitrary element of $\Mackey H_G^0(B_+)$ can be written as a sum
$a+b\kappa+c\epsilon$ for integers $a$, $b$, and $c$,
and consideration of when two such elements can multiply to give 1 leads to the conclusion
that only the eight elements above can be units.
\end{proof}

Notice that
\begin{align*}
 \kappa - \epsilon
 	&= \kappa - e^{-2}\kappa\zeta_0 c_{\omega} \\
	&= \kappa - e^{-2}\kappa((1-\kappa)\zeta_1 c_{\chi\omega} + e^2) \\
	&= e^{-2}\kappa\zeta_1 c_{\chi\omega} \\
	&= \chi\epsilon,
\end{align*}
so, where $\epsilon$ appears above, we can also expect to see $\kappa-\epsilon$, and we do.

Although the tautological bundle $\omega$ is a natural thing to consider,
in some applications, and in algebraic geometry, the dual of $\omega$
is important as well.
An example is when considering finite complex projective spaces as in \cite{CHT:projspaces}.
We write
\[
 \omega\dual = \Hom(\omega,\C).
\]
Note that
\[
 (\chi\omega)\dual = \chi(\omega\dual) = \Hom(\omega,\C^\sigma),
\]
so we can write $\chi\omega\dual$ unambiguously.
For applications we need to know the Euler classes of these bundles.

\begin{proposition}
\begin{align*}
	e(\omega\dual) &= -(1-\epsilon)c_{\omega} \\
\intertext{and}
	e(\chi\omega\dual) &= -(1-\kappa)(1-\epsilon)c_{\chi\omega}.
\end{align*}
\end{proposition}

\begin{proof}
The bundle $\omega\dual$ gives the same representation of $\Pi B$ as does $\omega$,
so $e(\omega\dual)$ lies in grading $\omega$.
In this grading, $\eta\colon \Mackey H_\GG^\omega(B_+)\to \Mackey H_\GG^\omega(B^\GG_+)$
is monomorphic, so it suffices to show the desired equality in $\Mackey H_\GG^\omega(B^\GG_+)$.
Consider $\omega\dual|B^0 = (\omega|B^0)\dual$.
As is true nonequivariantly, $\omega\dual$ is the same underlying real bundle as $\omega$ but
with the complex structure conjugated.
The bundle $\omega|B^0$ is the nonequivariant tautological bundle over $\CP^\infty$ with
trivial $\GG$-action on its fiber, so conjugating the $\GG$-action on the fibers 
has the same effect on the Euler class as locally
applying an orientation-reversing map $\R^2\to \R^2$ on each fiber.
Thinking of the effect on the Thom class, because the resulting map $S^2\to S^2$ represents
$-1$ in $A(\GG)$, the Thom class of $\omega\dual|B^0$ is the negative
of the Thom class of $\omega|B^0$, hence the same is true for the Euler classes.
That is,
\[
 e(\omega\dual)|B^0 = e(\omega\dual|B^0) = -e(\omega|B^0).
\]
On the other hand, $\omega|B^1$ is the nonequivariant tautological bundle with each fiber
isomorphic to $\C^\sigma = \R^{2\sigma}$ as a representation of $\GG$. Conjugating the complex structure
amounts locally to applying an orientation-reversing map $\R^{2\sigma}\to\R^{2\sigma}$ on each fiber.
The resulting map $S^{2\sigma}\to S^{2\sigma}$ represents the unit $\kappa-1$ in $A(\GG)$
rather than $-1$ (see the discussion at the end of \S\ref{sec:RepBurnsideRings}), so
\[
 e(\omega\dual)|B^1 = (\kappa-1)e(\omega|B^1).
\]
Putting together $e(\omega\dual)|B^0$ and $e(\omega\dual)|B^1$, we get
\[
 \eta(e(\omega\dual)) = (-1, \kappa-1)\eta(c_{\omega}).
\]

We also calculate 
\[
 \eta(1-\epsilon)
 	= \eta(1 - e^{-2}\kappa\zeta_0 c_{\omega}) 
	= (1, 1) - (0, \kappa)
	= (1, 1-\kappa),
\]
so
\[
 \eta(-(1-\epsilon)c_{\omega})
 	= (-1, \kappa-1)\eta(c_{\omega})
	= \eta(e(\omega\dual)),
\]
showing that $e(\omega\dual) = -(1-\epsilon)c_{\omega}$ as claimed.

The second equality follows on applying $\chi$:
\[
 e(\chi\omega\dual) = \chi e(\omega\dual) = \chi(-(1-\epsilon)c_{\omega})
  = -(1-\kappa+\epsilon)c_{\chi\omega}
  = -(1-\kappa)(1-\epsilon)c_{\chi\omega}. \qedhere
\]
\end{proof}

We can look at the calculation of the preceding proposition in the following way as well. Recalling that
\begin{align*}
 \eta(c_\omega) &= (c\zeta_1, (e^2+\xi c)\zeta_0^{-1}) \\
 \mathllap{\text{and}\qquad} \eta(c_{\chi\omega}) &= ((e^2 + \xi c)\zeta_1^{-1}, c\zeta_0),
\end{align*}
we have
\begin{align*}
 \eta(e(\omega\dual)) &= -(1-\epsilon)\eta(c_\omega) = (-c\zeta_1, (e^2 - \xi c)\zeta_0^{-1}) \\
 \mathllap{\text{and}\qquad} \eta(e(\chi\omega\dual)) 
 	&= -(1-\kappa)(1-\epsilon)\eta(c_{\chi\omega}) = ((e^2 - \xi c)\zeta_1^{-1}, -c\zeta_0).
\end{align*}
Dualizing negates the nonequivariant Chern class $c$, but multiplies the Euler class $e^2$
of $\C^\sigma$ by $\kappa-1$, which leaves it unchanged.

\section{Comparison to Lewis's calculation in $RO(\GG)$-grading}\label{sec:Lewis}

In \cite{Le:projectivespaces}, Gaunce Lewis calculated $\Mackey H_\GG^{RO(\GG)}(B_+)$,
the $RO(\GG)$-graded part 
of the cohomology of $B$. We here recover his results from ours.

Part of Lewis's Theorem~5.1 shows that $\Mackey H_\GG^{RO(\GG)}(B_+)$ is generated by
two elements we shall call $\gamma \in \Mackey H_\GG^{2\sigma}(B_+)$ (Lewis's $c$) and
$\Gamma\in\Mackey H_G^{2\sigma+2}(B_+)$ (Lewis's $C(1)$), with the single relation
\[
 \gamma^2 = e^2\gamma + \xi\Gamma.
\]
This results in an additive basis consisting of
\[
 1,\ \gamma,\ \Gamma,\ \gamma\Gamma,\ \Gamma^2,\ \gamma\Gamma^2,\ \ldots
\]
Our additive calculation gave a basis for 
$\Mackey H_\GG^{RO(\GG)}(B_+)$ consisting of
\[
 1,\ \zeta_0 c_{\omega},\ c_{\omega}c_{\chi\omega},\ \zeta_0 c_{\omega}^2 c_{\chi\omega},\ 
 	c_{\omega}^2 c_{\chi\omega}^2,\ \zeta_0 c_{\omega}^3 c_{\chi\omega}^2,\ \ldots
\]
(See Corollary~\ref{cor:basic}.)
If we let $\gamma = \zeta_0 c_{\omega}$ and $\Gamma = c_{\omega}c_{\chi\omega}$, we recover
Lewis's generators. 
Moreover, his relation between these generators can be seen as follows:
\begin{align*}
  \gamma^2
  	&= (\zeta_0 c_{\omega})^2 \\
	&= \zeta_0 c_{\omega}(e^2 + (1-\kappa)\zeta_1 c_{\chi\omega}) \\
	&= e^2\zeta_0 c_{\omega} + (1-\kappa)\xi c_{\omega} c_{\chi\omega} \\
	&= e^2 \gamma + \xi \Gamma,
\end{align*}
using the fact that $(1-\kappa)\zeta_0\zeta_1 = (1-\kappa)\xi = \xi$.

In his Remark~5.3, Lewis
introduces an element $\tilde\gamma = e^2 + (1-\kappa)\gamma$ (his $\tilde c$) and points
out that $\tilde\gamma$ could be used as a generator in place of $\gamma$.
In fact, $\tilde\gamma = \zeta_1 c_{\chi\omega}$, and his equation relating
$\gamma$ and $\tilde\gamma$ is our basic relation between $\zeta_0 c_{\omega}$
and $\zeta_1 c_{\chi\omega}$.

In \cite{CHT:projspaces} we also generalize Lewis's calculation for finite
projective spaces to the $RO(\Pi B)$-grading. In that case, 
our description is considerably simpler than the one he obtained.

\section{Cohomology with other coefficient systems}

The fact that the cohomology of $B$ with $\Mackey A$ coefficients is a free
$\Mackey H_\GG^{RO(\GG)}(S^0)$-module allows us to calculate its cohomology with any coefficient system.
We begin with a general result along the lines of Adams' splitting \cite{Ad:stablehomotopy}.
For the next several results, we let $G$ be any finite group, $\alpha\in RO(\Pi B)$, 
and let $X\to B$ be any ex-$G$-space.
By \cite[\S3.7]{CW:ordinaryhomology}, both $\Mackey H^G_{\alpha+RO(G)}(-;\Mackey A)$ and
$\Mackey H_G^{\alpha+RO(G)}(-;\Mackey A)$ are represented by a spectrum $H\Mackey A^\alpha$
parametrized by $B$. Let $\rho\colon B\to *$ denote the map to a $G$-fixed point.

\begin{proposition}\label{prop:spectrumsplitting}
Suppose that
$\Mackey H^G_{\alpha+RO(G)}(X;\Mackey A)$ is free as a module over
the ring $\Mackey H^G_{RO(G)}(S^0;\Mackey A)$.
Then the (nonparametrized) spectrum $\rho_!(X \smsh_B H\Mackey A^\alpha)$ is equivalent to a wedge of suspensions of
the spectrum $H\Mackey A$
representing $RO(G)$-graded cohomology.

Similarly, if $\Mackey H_G^{\alpha+RO(G)}(X;\Mackey A)$ is free as a module over
$\Mackey H_G^{\alpha+RO(G)}(S^0;\Mackey A)$,
then the spectrum $F_B(X,H\Mackey A^\alpha)$ is equivalent to a wedge of suspensions of
$H\Mackey A$.
\end{proposition}

\begin{proof}
Recall from \cite{CW:ordinaryhomology} that
\[
 \Mackey H^G_{\alpha+\beta}(X;\Mackey A)
  \iso [S^\beta, \rho_!(X \smsh_B H\Mackey A^\alpha)]_{G}
\]
for $\beta\in RO(G)$. Take a basis for $\Mackey H^G_{\alpha+RO(G)}(X;\Mackey A)$ 
over $\Mackey H^G_{RO(G)}(S^0;\Mackey A)$ and represent each basis element as a map
$S^\beta \to \rho_!(X \smsh_B H\Mackey A^\alpha)$. This gives a map
\[
 \bigvee S^\beta \to \rho_!(X \smsh_B H\Mackey A^\alpha),
\]
where the wedge is taken over the basis.
Using the fact that $H\Mackey A^\alpha$ is a module over $H\Mackey A$
(or $\rho^*H\Mackey A$), we then have the
composite
\[
 \bigvee S^\beta\smsh H\Mackey A 
  \to \rho_!(X \smsh_B H\Mackey A^\alpha) \smsh H\Mackey A
  \to \rho_!(X \smsh_B H\Mackey A^\alpha),
\]
which is a weak equivalence.

The proof for cohomology is similar, using
\[
 \Mackey H_G^{\alpha+\beta}(X;\Mackey A)
  \iso [S^\beta, F_B(X,H\Mackey A^\alpha)]_{G}. \qedhere
\]
\end{proof}

\begin{proposition}\label{prop:directsum}
Let $G$ be any finite group, let $X\to B$ be an ex-$G$-space, and suppose that
$\Mackey H^G_{\alpha+RO(G)}(X;\Mackey A)$ is a free module over $\Mackey H^G_{RO(G)}(S^0;\Mackey A)$
for an $\alpha\in RO(\Pi B)$. If $\Mackey T$ is any Mackey functor, then
\[
 \Mackey H^G_{\alpha+RO(G)}(X;\Mackey T) \iso
  \Mackey H^G_{\alpha+RO(G)}(X;\Mackey A)
  	\tensor_{\Mackey H^G_{RO(G)}(S^0;\Mackey A)} \Mackey H^G_{RO(G)}(S^0;\Mackey T).
\]
Similarly, 
if $\Mackey H_G^{\alpha+RO(G)}(X;\Mackey A)$ is a free module over $\Mackey H_G^{RO(G)}(S^0;\Mackey A)$,
then
\[
 \Mackey H_G^{\alpha+RO(G)}(X;\Mackey T) \iso
  \Mackey H_G^{\alpha+RO(G)}(X;\Mackey A)
  	\tensor_{\Mackey H_G^{RO(G)}(S^0;\Mackey A)} \Mackey H_G^{RO(G)}(S^0;\Mackey T).
\]
\end{proposition}

\begin{proof}
Consider homology, where we have the representing spectrum $H\Mackey T^\alpha$, meaning that
\[
 \Mackey H^G_{\alpha+\beta}(X;\Mackey T)
  \iso [S^\beta, \rho_!(X \smsh_B H\Mackey T^\alpha)]_{G}
\]
for $\beta\in RO(G)$.
Using the preceding proposition, we have
\begin{align*}
 \rho_!(X \smsh_B H\Mackey T^\alpha)
  &\hmtpc \rho_!(X\smsh_B H\Mackey A^\alpha \smsh_{H\Mackey A} \rho^*H\Mackey T) \\
  &\hmtpc \rho_!(X\smsh_B H\Mackey A^\alpha) \smsh_{H\Mackey A} H\Mackey T \\
  &\hmtpc \big(\bigvee S^\beta \smsh H\Mackey A\big) \smsh_{H\Mackey A} H\Mackey T \\
  &\hmtpc \bigvee S^\beta \smsh H\Mackey T,
\end{align*}
with the wedge taken over a basis for $\Mackey H^G_{\alpha+RO(G)}(X;\Mackey A)$ 
over $\Mackey H^G_{RO(G)}(S^0;\Mackey A)$.
This implies the algebraic claim of the proposition.

The proof for cohomology is similar, using $F_B(X,H\Mackey T^\alpha)$.
\end{proof}

Returning to the case $G = \GG$ and $B = \CP^\infty_G$, we can now easily compute the cohomology 
of $B$ with any coefficient system. The one that is probably of most interest is
the constant $\Z$ functor $\Mackey\Z$.
Refer to Theorem~\ref{thm:pointEvenRZcoeffs}
for the structure of $\Mackey H_\GG^{RO(\GG)}(S^0;\Mackey\Z)$.

\begin{theorem}\label{thm:BGU1cohomRZ}
$\Mackey H_\GG^{RO(\Pi B)}(B_+;\Mackey\Z)$ is the commutative $\Mackey H_\GG^{RO(\GG)}(S^0;\Mackey\Z)$-algebra generated by
the elements $\zeta_0$, $\zeta_1$, $c_{\omega}$, and $c_{\chi\omega}$ subject to the two relations
\[
  \zeta_0 \zeta_1 = \xi
\]
and
\[
  \zeta_1 c_{\chi\omega} = \zeta_0 c_{\omega} + e^2.
\]
It is free as a module over $\Mackey H_\GG^{RO(\GG)}(S^0;\Mackey\Z)$.
\end{theorem}

\begin{proof}
It follows from Proposition~\ref{prop:directsum} that
\[
 \Mackey H_\GG^{RO(\Pi B)}(B_+;\Mackey\Z)
  \iso \Mackey H_\GG^{RO(\Pi B)}(B_+) \tensorS \Mackey H_\GG^{RO(\GG)}(S^0;\Mackey\Z).
\]
By the proof of Theorem~\ref{thm:pointEvenRZcoeffs},
$\Mackey H_\GG^{RO(\GG)}(S^0;\Mackey\Z)$ is the quotient of $\Mackey H_\GG^{RO(\GG)}(S^0)$
by the ideal consisting of the elements $e^k\kappa$, $k\in \Z$.
(Recall that $e^k\kappa = 2e^k$ for $k > 0$.)
The result then follows readily from Theorem~\ref{thm:multstructure},
noting that the second relation simplifies because $\kappa = 0$ in the quotient.
\end{proof}

Corollary~\ref{cor:manyIsoCoeffs}
implies calculations of $\Mackey H_\GG^*(B_+;\Mackey T)$
with $\Mackey T$ equal to $\Mackey\Z'$, $\Mackey\Z_-$,
or $\Mackey\Z'_-$.
Each of these is probably best understood as 
a shifted version of $\Mackey H_\GG^*(B_+;\Mackey\Z)$, that is, as
a free module over $\Mackey H_\GG^*(B_+;\Mackey\Z)$
on a single generator.

\part{Invariants of the component structure}\label{part:basespace}

\section{The space $K(2)$}

The elements $\zeta_0$ and $\zeta_1$ that appear in the cohomology of
$B_\GG U(1)$ are examples of a much more general phenomenon.
They come from the cohomology of another space $K(2)$ whose
cohomology we will compute.
Its structure as a space is fairly simple.

\begin{definition}
Let $K = K(2)$ be a $\GG$-space of the homotopy type of a $\GG$-CW complex
such that $K$ is nonequivariantly contractible and $K^\GG$ consists
of two contractible components.
\end{definition}

Such a $\GG$-space exists and is unique up to $\GG$-homotopy equivalence
using Elmendorf's construction~\cite{El:fixedpoints}.
We could similarly construct spaces $K(n)$ with $K(n)^\GG$ having $n$ components,
but the case $n=2$ will suffice for our purposes here.
We could also construct $K$ as the classifying space $B\Pi_\GG B_\GG U(1)$
using the construction given in \cite[24.1]{CMW:orientation}.

\begin{remark}\label{rem:K2etilde}
Another model for $K(2)$ is $\tE\GG$ considered as an unbased space.
We computed the $RO(\GG)$-graded (reduced) cohomology of $\tE\GG$ in Theorem~\ref{thm:EvenTEG};
using a retraction $S^0\to (\tE\GG)_+\to S^0$, we see that the cohomology
of $(\tE\GG)_+$ is the direct sum of the cohomology of $\tE\GG$ and the cohomology of a point,
in the $RO(\GG)$ grading. However, this argument does not extend to other gradings
because the projection $(\tE\GG)_+\to S^0$ does not induce a map in cohomology in those other gradings.
Instead, we need an argument similar to the one we used in calculating the cohomology of
$B_\GG U(1)$.
We use the notation $K(2)$ to distinguish this unbased space from the based
space $\tE\GG$.
\end{remark}

Write $K^\GG = K^0 \disjunion K^1$. Elmendorf's construction also gives us
a $\GG$-map $B = B_\GG U(1)\to K$, unique up to $\GG$-homotopy, such that
$B^0$ maps to $K^0$ and $B^1$ maps to $K^1$. This map induces an isomorphism
$RO(\Pi K) \iso RO(\Pi B)$, so we will use the same notation for elements
of $RO(\Pi K)$ as we did for $RO(\Pi B)$.

The calculation of $\Mackey H_\GG^{RO(\Pi K)}(K_+)$ is similar to the calculation of the cohomology
of $B$, with some parts much simpler but some parts, including the final result, actually more complicated.
We begin with some of the simple calculations.

\begin{proposition}\label{prop:calcFixedPointsC}
For $k = 0$ or $1$, we have the following calculation:
\begin{align*}
 \Mackey H_\GG^{RO(\Pi K)}(K^k_+) 
   &\iso \Mackey H_\GG^{RO(\GG)}(S^0)[\zeta_{1-k}, \zeta_{1-k}^{-1}] \\
 \Mackey H_\GG^{RO(\Pi K)}(K^k_+\smsh_K (E\GG)_+) 
   &\iso \Mackey H_\GG^{RO(\Pi K)}(K^k_+)\tensor_{\Mackey H_\GG^{RO(\GG)}(S^0)} \Mackey H_\GG^{RO(\GG)}((E\GG)_+) \\
 \Mackey H_\GG^{RO(\Pi K)}(K^k_+\smsh_K \tE\GG) 
   &\iso \Mackey H_\GG^{RO(\Pi K)}(K^k_+)\tensor_{\Mackey H_\GG^{RO(\GG)}(S^0)} \Mackey H_\GG^{RO(\GG)}(\tE\GG),
\end{align*}
where $|\zeta_i| = \Omega_i$.
Further,
\begin{align*}
 \Mackey H_\GG^{RO(\Pi K)}(K^\GG_+)
   &\iso \Mackey H_\GG^{RO(\Pi K)}(K^0_+) \dirsum \Mackey H_\GG^{RO(\Pi K)}(K^1_+)  \\
 \Mackey H_\GG^{RO(\Pi K)}(K^\GG_+\smsh_K (E\GG)_+) 
   &\iso 
    \Mackey H_\GG^{RO(\Pi K)}(K^\GG_+)\tensor_{\Mackey H_\GG^{RO(\GG)}(S^0)} \Mackey H_\GG^{RO(\GG)}((E\GG)_+) \\
 \Mackey H_\GG^{RO(\Pi K)}(K^\GG_+\smsh_K \tE\GG) 
   &\iso 
    \Mackey H_\GG^{RO(\Pi K)}(K^\GG_+)\tensor_{\Mackey H_\GG^{RO(\GG)}(S^0)} \Mackey H_\GG^{RO(\GG)}(\tE\GG).
\end{align*}
\end{proposition}

\begin{proof}
The proof is essentially the same as that of Proposition~\ref{prop:calcFixedPoints},
simplified by the fact that each $K^k$ is contractible.
\end{proof}

Also simple are the cohomologies of $K_+\smsh_K (E\GG)_+$ and $K_+\smsh_K \tE\GG$.

\begin{proposition}\label{prop:CEGplus}
We have
\[
 \Mackey H_\GG^{RO(\Pi K)}(K_+\smsh_K (E\GG)_+) \iso 
  \Mackey H_\GG^{RO(\GG)}((E\GG)_+)[\zeta_0,\zeta_1] / 
  		\langle \zeta_0\zeta_1 = \xi \rangle.
\]
The long exact sequence of the pair $(K,K^\GG)$
reduces to a split short exact sequence
\begin{multline*}
 0 \to \Mackey H_\GG^{RO(\Pi K)}(K_+\smsh_K (E\GG)_+) \xrightarrow{\eta}
  \Mackey H_\GG^{RO(Pi K)}(K^\GG_+\smsh_K (E\GG)_+) \\*
  \to
  \susp^{-1}\Mackey H_\GG^{RO(\Pi K)}(K/_K K^\GG \smsh_K (E\GG)_+) \to 0
\end{multline*}
with
\begin{align*}
	\eta(\zeta_0) &= (\xi\zeta_1^{-1}, \zeta_0) \qquad\text{and} \\
	\eta(\zeta_1) &= (\zeta_1, \xi\zeta_0^{-1}).
\end{align*}
\qed
\end{proposition}

The proof is the similar to the proofs of the corresponding statements about $B_+\smsh_B (E\GG)_+$
in \S\ref{sec:egcohomologies}, but simpler.

Note that, because $\xi$ is invertible in $\Mackey H_\GG^{RO(\GG)}((E\GG)_+)$, the relation
$\zeta_0\zeta_1 = \xi$ implies that both $\zeta_0$ and $\zeta_1$ are invertible.

\begin{proposition}\label{prop:KEGtilde}
The inclusion $K^\GG\to K$ induces an isomorphism
\begin{align*}
 \Mackey H_\GG^{RO(\Pi K)}(K_+&{}\smsh_K \tE\GG) \\
  &\iso \Mackey H_\GG^{RO(\Pi K)}(K^\GG_+\smsh_K \tE\GG) \\
  &\iso \left(\Mackey H_\GG^{RO(\GG)}(S^0)[\zeta_1,\zeta_1^{-1}]\dirsum
  				\Mackey H_\GG^{RO(\GG)}(S^0)[\zeta_0,\zeta_0^{-1}] \right) \\
     &\qquad\qquad \tensorS \Mackey H_\GG^{RO(\GG)}(\tE\GG).
\end{align*}
\end{proposition}

\begin{proof}
As in Proposition~\ref{prop:BEGtilde}, this is just the observation that
the inclusion induces a weak equivalence
$K^\GG_+\smsh \tE\GG \to K_+\smsh \tE\GG$.
\end{proof}

\section{The cohomology of $K(2)$}\label{sec:basespace}

Unlike $B$, the cohomology of $K$ is not a free module over $\Mackey H_\GG^{RO(\GG)}(S^0)$.
This is where the argument gets more complicated.

Diagram~(\ref{diag:main}) works equally well with
$K$ in place of $B$.
By the same argument used for the cohomology of $B$, we can find elements
$\zeta_i \in \Mackey H_\GG^{\Omega_i}(K_+)$, $k = 0$ and $1$, characterized by
\begin{align*}
	\eta(\zeta_0) &= (\xi\zeta_1^{-1}, \zeta_0) \qquad\text{and} \\
	\eta(\zeta_1) &= (\zeta_1, \xi\zeta_0^{-1}).
\end{align*}
These map to the elements of $\Mackey H_\GG^{RO(\Pi K)}(K_+\smsh_K (E\GG)_+)$ of the same names that 
we saw in the preceding section.
We have again that $\zeta_0\zeta_1 = \xi$.

We define
\[
 \Mackey Q^{RO(\Pi K)} = \Mackey H_\GG^{RO(\GG)}(S^0)[\zeta_0, \zeta_1]/
  \langle \zeta_0\zeta_1 = \xi \rangle.
\]
The following is easily proved.

\begin{proposition}
$\Mackey Q^{RO(\Pi K)}$ is a free $\Mackey H_\GG^{RO(\GG)}(S^0)$-module on a basis consisting
of the elements $1$, $\zeta_0^n$ for $n\geq 1$, and $\zeta_1^n$ for $n\geq 1$.
\qed
\end{proposition}

To have a simple way of referring to this basis, we make the following definition.

\begin{definition}\label{def:basicC}
We say that a monomial $\zeta_0^m\zeta_1^n$ is {\em basic} if $m=0$ or $n=0$.
\end{definition}

These basis elements are distributed nicely:

\begin{lemma}
Let $\alpha\in RO(\Pi K)$. Then there is exactly one basic
monomial with grading in the coset $\alpha + RO(\GG) \subset RO(\Pi K)$.
\end{lemma}

\begin{proof}
Write $\alpha = m\Omega_0 + n\Omega_1 + \gamma$ where $\gamma\in RO(\GG)$.
If $m\geq n$, then
\[
 m\Omega_0 + n\Omega_1 + \gamma = (m-n)\Omega_0 + \gamma + n(2\sigma-2),
\]
and the only basic monomial having grading in $(m-n)\Omega_0 + RO(\GG)$ is $\zeta_0^{m-n}$.
Similarly, if $m < n$, the only basic monomial with grading in $\alpha+RO(\GG)$
is $\zeta_1^{n-m}$.
\end{proof}

Consider the algebra map 
$f\colon \Mackey Q^{RO(\Pi K)}\to \Mackey H_\GG^{RO(\Pi K)}(K_+)$ given by taking each
$\zeta_i$ to the element of the same name. As in the case of $B$, we have a map of
long exact sequences:
\[
 \xymatrix{
 \vdots \ar[d] & \vdots \ar[d] \\
 \Mackey Q^{RO(\Pi K)}\tensorS \Mackey H_\GG^{RO(\GG)}(\tE\GG) \ar[r]^-{f_{\tE\GG}} \ar[d]
 	& \Mackey H_\GG^{RO(\Pi K)}(K_+\smsh_K \tE\GG) \ar[d] \\
 \Mackey Q^{RO(\Pi K)} \ar[r]^f \ar[d]
	& \Mackey H_\GG^{RO(\Pi K)}(K_+) \ar[d] \\
 \Mackey Q^{RO(\Pi K)}\tensorS \Mackey H_\GG^{RO(\GG)}((E\GG)_+) \ar[r]^-{f_{(E\GG)_+}} \ar[d]
	& \Mackey H_\GG^{RO(\Pi K)}(K_+\smsh_K (E\GG)_+) \ar[d] \\
 \vdots & \vdots
 }
\]
We shall show that $f_{(E\GG)_+}$ is again an isomorphism, but in this case
$f_{\tE\GG}$ is {\em not}, and we shall see where that leads us.

\begin{proposition}
\begin{multline*}
 f_{(E\GG)_+}\colon \Mackey Q^{RO(\Pi K)}\tensorS \Mackey H_\GG^{RO(\GG)}((E\GG)_+) \\
	\to \Mackey H_\GG^{RO(\Pi K)}(K_+\smsh_K (E\GG)_+)
\end{multline*}
is an isomorphism.
\end{proposition}

\begin{proof}
We have that
\begin{multline*}
 \Mackey Q^{RO(\Pi K)}\tensorS \Mackey H_\GG^{RO(\GG)}((E\GG)_+)  \\ \iso
    \Mackey H_\GG^{RO(\GG)}((E\GG)_+)[\zeta_0, \zeta_1]/
  \langle \zeta_0\zeta_1 = \xi \rangle,
\end{multline*}
and the result now follows from Proposition~\ref{prop:CEGplus}.
\end{proof}

\begin{proposition}
\[
 f_{\tE\GG}\colon \Mackey Q^{RO(\Pi K)}\tensorS \Mackey H_\GG^{RO(\GG)}(\tE\GG) 
 \to \Mackey H_\GG^{RO(\Pi K)}(K_+\smsh_K \tE\GG)
\]
is a monomorphism, but not an isomorphism. It is split as a map of modules over
$\Mackey H_\GG^{RO(\GG)}(S^0)$.
\end{proposition}

\begin{proof}
Consider a grading $m\Omega_0 + RO(\GG)$, $m\geq 0$, and the corresponding
basic monomial $\zeta_0^m$ that lives in this grading.
(The case of $n\Omega_1 + RO(\GG)$ is similar.)

From the calculation in the preceding section,
in gradings $m\Omega_0+RO(\GG)$ the target
of $f_{\tE\GG}$ is
\[
 \left(\Mackey H_\GG^{RO(\GG)}(S^0)\{\zeta_1^{-m}\}
 	\dirsum \Mackey H_\GG^{RO(\GG)}(S^0)\{\zeta_0^{m}\}\right)
  \tensor_{\Mackey H_\GG^{RO(\GG)}(S^0)} \Mackey H_\GG^{RO(\GG)}(\tE\GG).
\]
If $x\in \Mackey H_\GG^{RO(\GG)}(\tE\GG)$, we have
\[
 f_{\tE\GG}\big(\zeta_0^m \tensor x\big)
  = (\xi^m\zeta_1^{-m}\tensor x, \zeta_0^m\tensor x).
\]
From this we can see that projection to the second summand is
a splitting of $f_{\tE\GG}$ as a map of $\Mackey H_\GG^{RO(\GG)}(S^0)$-modules in gradings $m\Omega_0+RO(\GG)$.

The cokernel of $f_{\tE\GG}$ in these gradings is isomorphic to the first summand,
so is nontrivial, hence $f_{\tE\GG}$ is not an isomorphism.
\end{proof}

We can now calculate the cohomology of $K$.

\begin{theorem}
As a module over $\Mackey H_\GG^{RO(\GG)}(S^0)$, $\Mackey H_\GG^{RO(\Pi K)}(K_+)$ is the pushout in the diagram
\[
 \xymatrix{
  \Mackey Q^{RO(\Pi K)}\tensorS \Mackey H_\GG^{RO(\GG)}(\tE\GG) \ar[r] \ar[d]_{f_{\tE\GG}}
    & \Mackey Q^{RO(\Pi K)} \ar[d]^f \\
  \Mackey H_\GG^{RO(\Pi K)}(K_+\smsh_K \tE\GG) \ar[r] & \Mackey H_\GG^{RO(\Pi K)}(K_+)
 }
\]
in which the vertical arrows are both split monomorphisms. So
$\Mackey H_\GG^{RO(\Pi K)}(K_+)$ splits as a direct sum with
$\Mackey Q^{RO(\Pi K)}$ as one summand, but the
other summand is isomorphic to the direct sum of countably infinitely many copies of 
$\Mackey H_\GG^{RO(\GG)}(\tE\GG)$, with one copy in each coset of gradings $\alpha+RO(\GG)\subset RO(\Pi K)$.

As an algebra over $\Mackey H_\GG^{RO(\GG)}(S^0)$, $\Mackey H_\GG^{RO(\Pi K)}(K_+)$ is generated by elements
\begin{align*}
 \zeta_0, \zeta_1, \\
\intertext{with $\deg\zeta_0 = \Omega_0$ and $\deg\zeta_1 = \Omega_1$, and}
 x\psi_0^m, x\psi_1^m &\qquad x\in \Mackey H_\GG^{RO(\GG)}(\tE\GG),\ m\in\Z
\end{align*}
with $\deg x\psi_0^m = \deg x + m\Omega_1$ and $\deg x\psi_1^m = \deg x + m\Omega_0$.
All relations among these generators are consequences of the following relations:
\begin{align*}
 \zeta_0\zeta_1 &= \xi \\
 y(x\psi_k^m) &= (yx)\psi_k^m && \text{if $y \in \Mackey H_\GG^{RO(\GG)}(S^0)$} \\
 (x\psi_0^m)(y\psi_0^n) &= xy\psi_0^{m+n} \\
 (x\psi_1^m)(y\psi_1^n) &= xy\psi_1^{m+n} \\
 (x\psi_0^m)(y\psi_1^n) &= 0 \\
 \zeta_0\cdot x\psi_0^m &= \xi x\psi_0^{m-1} \\
 \zeta_1\cdot x\psi_0^m &= x\psi_0^{m+1} \\
 \zeta_0\cdot x\psi_1^m &= x\psi_1^{m+1} \\
 \zeta_1\cdot x\psi_1^m &= \xi x\psi_1^{m-1} \\
 \psi(x)\zeta_0 &= \xi x\psi_0^{-1} + x\psi_1^1 && \text{for $x\in\Mackey H_\GG^{RO(\GG)}(\tE\GG)$} \\
 \psi(x)\zeta_1 &= x\psi_0^1 + \xi x\psi_1^{-1} && \text{for $x\in\Mackey H_\GG^{RO(\GG)}(\tE\GG)$}
\end{align*}
In the last two relations, $\psi(x)$ is the image of $x$ under the map
\[
	\psi\colon H_\GG^{RO(\GG)}(\tE\GG) \to H_\GG^{RO(\GG)}(S^0).
\]
\end{theorem}

\begin{proof}
In the map of long exact sequences displayed after Definition~\ref{def:basicC},
we now know that $f_{(E\GG)_+}$ is an isomorphism and $f_{\tE\GG}$ is a monomorphism.
A diagram chase shows that these imply that the square involving $f_{\tE\GG}$ and $f$ is
a pushout square, as claimed in the theorem. 
Moreover, because $f_{\tE\GG}$ is a split monomorphism, its pushout, $f$, is as well.

Note that the cokernel of $f$ is therefore isomorphic to the cokernel of $f_{\tE\GG}$,
which is a direct sum of countably infinitely many copies of $\Mackey H_\GG^{RO(\GG)}(\tE\GG)$,
with one copy in each coset of gradings $\alpha+RO(\GG) \subset RO(\Pi K)$.

The analogue of Proposition~\ref{prop:etamono} holds for $K$, by a similar proof as for $B$.
We have already defined the elements $\zeta_i \in \Mackey H_G^{\Omega_i}(K_+)$.
We define the element $x\psi_k^m$ to be the image in $\Mackey H_\GG^{RO(\Pi K)}(K_+)$
of $\psi_k^m\tensor x \in \Mackey H_\GG^{RO(\Pi K)}(K_+\smsh_K \tE\GG)$,
where
\begin{align*}
	\psi_0^m &= (\zeta_1^m, 0) \qquad\text{and} \\
	\psi_1^m &= (0, \zeta_0^m).
\end{align*}
(See Proposition~\ref{prop:KEGtilde}.)
It follows that
\begin{align*}
	\eta(x\psi_0^m) &= (\psi(x)\zeta_1^m, 0) \qquad\text{and} \\
	\eta(x\psi_1^m) &= (0, \psi(x)\zeta_0^m).
\end{align*}

From the definition of $\Mackey Q^{RO(\Pi K)}$ and our calculation of
$\Mackey H_\GG^{RO(\Pi K)}(K_+\smsh_K \tE\GG)$, we see that the elements $\zeta_i$ and
$x\psi_k^m$ generate $\Mackey H_\GG^{RO(\Pi K)}(K_+)$ multiplicatively.

The relation $\zeta_0\zeta_1 = \xi$ holds because it is true on applying $\eta$ and takes
place in a grading in which $\eta$ is a monomorphism.
The next four relations listed in the theorem hold because they do in 
$\Mackey H_\GG^{RO(\Pi K)}(K_+\smsh_K \tE\GG)$.

To verify the formulas for $\zeta_i \cdot x\psi_k^m$, consider the following
diagram:
\[\def\objectstyle{\scriptstyle}
 \xymatrix{
   \Mackey H_\GG^{RO(\Pi K)}(K_+) \tensor \Mackey H_\GG^{RO(\Pi K)}(K_+\smsh\tE\GG) \ar[d] \ar[r]^\eta
    & \Mackey H_\GG^{RO(\Pi K)}(K^\GG_+)\tensor \Mackey H_\GG^{RO(\Pi K)}(K^\GG_+\smsh\tE\GG) \ar[d] \\
   \Mackey H_\GG^{RO(\Pi K)}(K_+\smsh K_+\smsh \tE\GG) \ar[d] \ar[r]^\iso
    & \Mackey H_\GG^{RO(\Pi K)}(K^\GG_+\smsh K^\GG_+\smsh \tE\GG) \ar[d] \\
   \Mackey H_\GG^{RO(\Pi K)}(K_+\smsh\tE\GG) \ar[d] \ar[r]^\iso
    & \Mackey H_\GG^{RO(\Pi K)}(K^\GG_+\smsh\tE\GG) \\
   \Mackey H_\GG^{RO(\Pi K)}(K_+)
  }
\]
This diagram implies that it suffices to check that the formula holds in the cohomology group
$\Mackey H_\GG^{RO(\Pi K)}(K^\GG_+\smsh\tE\GG)$, where it is easy to check component-wise.

The relation $\psi(x)\zeta_0 = \xi x\psi_0^{-1} + x\psi_1^1$
follows from the commutativity of the diagram in the theorem, as
the two sides are the images of $\zeta_0\tensor x$ around the two sides of the diagram.
Similarly for the last relation.

That the listed relations imply all relations follows from the fact that they easily allow
us to write any product of generators as a linear combination of basic monomials from
$\Mackey Q^{RO(\Pi K)}$ and elements $x\psi_k^m$ from
$\Mackey H_\GG^{RO(\Pi K)}(K_+\smsh_K \tE\GG)$.
\end{proof}

\section{$\zeta_0$ and $\zeta_1$ as invariants of the component structure}\label{sec:components}

Consider any $\GG$-space $X$ and 
suppose that we have a decomposition $X^\GG = X^0 \disjunion X^1$,
so each of $X^0$ and $X^1$ is a union of components of $X^\GG$.
There is a $\GG$-map $k\colon X\to K(2)$, unique up to $\GG$-homotopy,
sending $X^0$ to $K^0$ and $X^1$ to $K^1$.
Then we have the elements $k^*\zeta_0 \in \Mackey H_\GG^{k^*\Omega_0}(X_+)$
and $k^*\zeta_1 \in \Mackey H_\GG^{k^*\Omega_1}(X_+)$,
where we also write $k^*\colon RO(\Pi K)\to RO(\Pi X)$ for the induced map.
We think of $k^*\zeta_0$ and $k^*\zeta_1$ as invariants of the component structure of $X$
and we write $\zeta_0 = k^*\zeta_0$ and $\zeta_1 = k^*\zeta_1$ for simplicity of notation.
This is exactly how the elements $\zeta_0$ and $\zeta_1$ in the
cohomology of $B_\GG U(1)$ arise.
Note that we will always have $\zeta_0\zeta_1 = \xi$ because that relation is
true in the cohomology of $K$.

One phenomenon that showed up
in the computation of the cohomology of finite projective spaces in \cite{CHT:projspaces}
that does not occur in the cohomology of $B_\GG U(1)$
was that some nonzero cohomology elements were divisible by $\zeta_0$ or $\zeta_1$.
A general result is the following.
For $x$ a cohomology element, let $x|X^0$ be the image of $x$ under the restriction
map $\Mackey H_\GG^{RO(\Pi X)}(X_+) \to \Mackey H_\GG^{RO(\Pi X)}(X^0_+)$.

\begin{proposition}\label{prop:divisibility}
Let $x\in H_\GG^{RO(\Pi X)}(X_+)$ be an element such that
$x|X^0 = 0$. Then $x$ is infinitely divisible by $\zeta_0$.
\end{proposition}

A similar statement is true for $\zeta_1$, by symmetry. We first need a lemma.

\begin{lemma}
Let $f\colon Y\to X$ be an ex-$\GG$-space over $X$ and suppose that $f^{-1}(X^0)$
is contractible to the base section.
Then $f^*\zeta_0$ is an invertible element of $\Mackey H_\GG^{RO(\Pi X)}(Y)$.
\end{lemma}

\begin{proof}
Let $i\colon *\to K$ be the inclusion of a fixed point in $K^1$.
By its construction, $i^*\zeta_0 \in \Mackey H_\GG^{RO(\Pi K)}(S^0)$ is invertible.

Using Elmendorf's construction, we can construct a $\GG$-map $h\colon \tilde X\to X$
such that $h$ is a nonequivariant equivalence and $h^\GG$ is equivalent to the inclusion 
$X^\GG\setminus X^0\includesin X^\GG$.
The composite $k\circ h\colon \tilde X\to K$ factors, up to homotopy, through $i$.
Therefore, $h^*\zeta_0$ is invertible
in $\Mackey H_\GG^{RO(\Pi X)}(\tilde X_+)$.

Let $Y^0 = (f^\GG)^{-1}(X^0)$, which we are assuming is contractible to the base section.
Using Elmendorf's construction again, we can construct an ex-$\GG$-space $h\colon\tilde Y\to \tilde X$
and a diagram
\[
 \xymatrix{
 \tilde Y \ar[r]^\ell \ar[d]_-{\tilde f} &  Y \ar[d]^-f \\
 \tilde X \ar[r]_h & X
 }
\]
such that the map $h_!\tilde Y\to  Y$ over $X$ is a nonequivariant equivalence and, on fixed points,
is equivalent to the inclusion $(Y^\GG\setminus Y^0)_+ \includesin Y^\GG$,
hence is an equivalence by the assumption on $Y$.
(Here, $(Y^\GG\setminus Y^0)_+$ is the pushout $(Y^\GG\setminus Y^0)\union_{\sigma(X)\setminus Y^0} X$.)
It follows that
\[
    \ell^*\colon H_\GG^{RO(\Pi X)}(Y) \iso H_\GG^{RO(\Pi X)}(h_! \tilde Y^0).
\]
But $\ell^*f^*\zeta_0 = \tilde f^*h^*\zeta_0$ is invertible, hence $f^*\zeta_0$ is invertible.
\end{proof}

\begin{proof}[Proof of Proposition~\ref{prop:divisibility}]
Let $i\colon X^0\to X$ be the inclusion,
let $Ci$ be the cofiber over $X$ of $i$, and let $j\colon X_+\to Ci$ be inclusion in the cofiber.
If $i^*x = 0$, then $x = j^*y$ for some $y$.
Now $\zeta_0$ is invertible in the cohomology of $Ci$ by the preceding lemma, because the part of $(Ci)^\GG$ lying over $X^0$
is contractible to the base section.
So, for each $n>0$, there is a $z$ such that
$y = \zeta_0^n z$. It follows that $x = \zeta_0^n j^*z$, so $x$ is divisible by $\zeta_0^n$.
\end{proof}

Using the spaces $K(n)$ for $n>2$ we can generalize this construction.
If we have a decomposition
$X^\GG = \Disjunion_{k=0}^{n-1} X^k$, we can define elements $\zeta_k$ associated to each $X^k$
satisfying the appropriate version of Proposition~\ref{prop:divisibility}
and with $\prod_{k=0}^{n-1} \zeta_k = \xi$.

Note that Proposition~\ref{prop:divisibility}
is not an ``if and only if'' result because $\zeta_0|X^0\neq 0$ in general
(for example, in $K$).
In the applications so far, the present proposition has been sufficient.

\section{$B_\GG U(1)$ does not represent cohomology}

We noted earlier that we have a map $B = B_\GG U(1) \to K = K(2)$ mapping $B^0$ to $K^0$ and $B^1$ to $K^1$
and inducing an isomorphism $RO(\Pi B) \iso RO(\Pi K)$.
In this way we consider $B$ as a space parametrized by $K$,
which allows us to consider $\Mackey H_\GG^{RO(\Pi B)}(B_+)$ 
as an algebra over $\Mackey H_\GG^{RO(\Pi K)}(K_+)$,
via a ring map
\[
 \Mackey H_\GG^{RO(\Pi K)}(K_+) \to \Mackey H_\GG^{RO(\Pi B)}(B_+).
\]
We know that $\zeta_0 \mapsto \zeta_0$ and $\zeta_1 \mapsto \zeta_1$,
but what about the elements $x\psi_k^m \in \Mackey H_\GG^{RO(\Pi K)}(K_+)$?
Recall that 
\begin{align*}
	\eta(x\psi_0^m) &= (\psi(x)\zeta_1^m, 0) \qquad\text{and} \\
	\eta(x\psi_1^m) &= (0, \psi(x)\zeta_0^m).
\end{align*}
From this and the formulas used in the proof of Theorem~\ref{thm:multstructure} 
it is possible to work out the image of $x\psi_k^m$
in $\Mackey H_\GG^{RO(\Pi B)}(B_+)$.
It suffices for our purposes here that $x\psi_k^m$ maps to $0$ in $\Mackey H_\GG^{RO(\Pi B)}(B_+)$
if $\psi(x) = 0$.

From our calculations we get the following result.

\begin{proposition}
There is no equivariant section of the map $B\to K$, even up to homotopy.
\end{proposition}

\begin{proof}
If there were a homotopy section, then the ring map 
$\Mackey H_\GG^{RO(\Pi K)}(K_+)\to \Mackey H_\GG^{RO(\Pi B)}(B_+)$
would be injective. However, the element
$(\delta\xi^{-1})\psi_0^0$ is nonzero in $\Mackey H_\GG^{RO(\Pi K)}(K_+)$, but maps
to 0 in $\Mackey H_\GG^{RO(\Pi B)}(B_+)$ because $\psi(\delta\xi^{-1}) = 0$.
\end{proof}

\begin{remark}
This result settles a point on which the author was confused for a while.
The $\GG$-space $B$ is nonequivariantly equivalent to $K(\Z,2)$ and $B^\GG$ is equivalent
to the disjoint union of two copies of $K(\Z,2)$. This suggests that
$B$ should represent the functor $\Mackey H_\GG^2(-;\Mackey\Z)$ on spaces over $K$,
but the cohomology of $B$ just seems wildly wrong for that: What would be the
universal cohomology element in $\Mackey H_\GG^2(B;\Mackey\Z)$?
But the space representing $\Mackey H_G^2(-;\Mackey\Z)$ must first be
an ex-space over $K$, and the proposition shows that $B$ cannot be made into one,
hence it is not a candidate to represent any cohomology group.

In fact, the ex-space representing $\Mackey H_\GG^2(-;\Mackey\Z)$ is simply
$K\times K(\Z,2)$ for a nonequivariant $K(\Z,2)$, say $\CP^\infty$,
taken with trivial $\GG$-action, and with a section given by any chosen
basepoint in $K(\Z,2)$. The proposition shows that $B$ is not equivariantly
equivalent to $K\times K(\Z,2)$.

Thus, unlike the nonequivariant case, $\GG$-equivariant complex line bundles 
are not classified by a cohomology class (the Euler class in the nonequivariant case),
as $B_\GG U(1)$ is not the representing space for any cohomology group.
\end{remark}

\begin{corollary}
There is no complex line bundle over $K$ who associated representation is
$\omega\in RO(\Pi K) = RO(\Pi B)$.
\end{corollary}

\begin{proof}
If there were such a line bundle, its classifying map
$K \to B_\GG U(1)$ would be a (homotopy) section of the projection $B\to K$.
But the preceding proposition shows that no such map exists.
\end{proof}

\begin{remark}
It is often convenient to know that a given representation
of a fundamental groupoid comes from a bundle over the space,
but this result shows that, at least for complex bundles, this is
not something that is guaranteed to happen.
\end{remark}

\bibliographystyle{amsplain}

\providecommand{\bysame}{\leavevmode\hbox to3em{\hrulefill}\thinspace}
\providecommand{\MR}{\relax\ifhmode\unskip\space\fi MR }
\providecommand{\MRhref}[2]{%
  \href{http://www.ams.org/mathscinet-getitem?mr=#1}{#2}
}
\providecommand{\href}[2]{#2}

\end{document}